\documentclass[a4paper]{amsart}
\usepackage{microtype}
\usepackage{amsmath,amsfonts,amssymb,mathrsfs,stmaryrd,mathtools}
\usepackage{mathbbol,wasysym}
\usepackage[backref=pages]{hyperref}
\usepackage[hyperpageref]{backref}
\hypersetup{pdftitle={Seminormal forms and cyclotomic quiver Hecke algebras of type A},
  pdfauthor={Jun Hu and Andrew Mathas},
  colorlinks = true,
  linkcolor =blue,
  anchorcolor = red,
  citecolor = blue,
  urlcolor = blue
}

\renewcommand*{\backref}[1]{}
\renewcommand*{\backrefalt}[4]{{\tiny[%
  \ifcase #1 No citations.
  \or Page #2.%
  \else Pages #2.%
  \fi%
  ]}%
}
\DeclareSymbolFontAlphabet{\mathbb}{AMSb}
\usepackage{cite}

\newcommand\Comment[2][Hu]{\space\par\medskip\noindent%
   \fbox{\begin{minipage}{\textwidth}\textbf{Comment\ifx\relax#1\else---#1\fi}\newline%
        #2\end{minipage}}\medskip
}
\usepackage{tikz}
\usetikzlibrary{matrix}

\newcount\tableauRow\newcount\tableauCol
\newcommand\Tableau[3][-4]{
  \begin{tikzpicture}[scale=#2,draw/.append style={thick,black},baseline=#1mm]
    \tableauRow=0
    \foreach \Row in {#3} {
       \tableauCol=1
       \foreach\k in \Row {
          \draw(\the\tableauCol,\the\tableauRow)+(-.5,-.5)rectangle++(.5,.5);
          \draw(\the\tableauCol,\the\tableauRow)node{\k};
          \global\advance\tableauCol by 1
       }
       \global\advance\tableauRow by -1
    }
  \end{tikzpicture}
}
\def\Bitab(#1|#2){\Bigg(%
  \hspace*{1mm}\Tableau{0.5}{#1}\hspace*{2mm}\Bigg|%
  \hspace*{2mm}\Tableau{0.5}{#2}\hspace*{1mm}\Bigg)
}
\def\Tritab(#1|#2|#3){\Bigg(%
  \hspace*{1mm}\Tableau{0.5}{#1}\hspace*{1mm}\Bigg|%
  \hspace*{1mm}\Tableau{0.5}{#2}\hspace*{1mm}\Bigg|%
  \hspace*{1mm}\Tableau{0.5}{#3}\hspace*{1mm}\Bigg)
}
\def\tritab(#1|#2|#3){\big(%
  \Tableau[-1]{0.3}{#1}\big|\Tableau[-1]{0.3}{#2}\big|\Tableau[-1]{0.3}{#3}\big)
}

\newcommand\YoungDiagram[2][\relax]{
  \begin{tikzpicture}[scale=0.5,draw/.append style={thick,black},baseline=-2mm]
    \ifx\relax#1\relax%
    \else 
    \foreach\box in {#1} {
      \filldraw[blue!30]\box rectangle ++(1,1);
    }
    \fi
    \newcount\tableauRow
    \tableauRow=0
    \foreach \diagramCol in {#2} {
       \draw(1,\the\tableauRow)grid ++(\diagramCol,1);
       \global\advance\tableauRow by -1
    }
  \end{tikzpicture}
}
\def\TriDiagram(#1|#2|#3){\Bigg(%
  \hspace*{1mm}\YoungDiagram{#1}\hspace*{1mm}\Bigg|\hspace*{1mm}\YoungDiagram{#2}%
  \hspace*{1mm}\Bigg|\hspace*{1mm}\YoungDiagram{#3}\hspace*{1mm}\Bigg)
}

\let\Sect=\S
\let\<=\langle
\let\>=\rangle
\let\[=\llbracket
\let\]=\rrbracket
\def\({\big(}
\def\){\big)}

\def\Sum{\displaystyle\sum}
\def\bijection{\overset{\sim}{\longrightarrow}}
\let\surjection\twoheadrightarrow

\newcommand{\C}{\mathbb{C}}

\newcommand{\N}{\mathbb{N}}

\newcommand{\Q}{\mathbb{Q}}
\newcommand{\Z}{\mathbb{Z}}

\newcommand{\Sym}{\mathfrak S}

\newcommand\SetBox[3][60]{\Big\{\ #2\ \Big| \ \vcenter{\hsize#1mm\centering #3}\Big\}}
\newcommand\Setbox[2][34]{\Big\{\vcenter{\hsize#1mm\centering#2}\Big\}}

\def\map#1#2{\,{:}\,#1\!\longrightarrow\!#2}
{\catcode`\|=\active
  \gdef\set#1{\mathinner{\lbrace\,{\mathcode`\|"8000%
                                   \let|\midvert #1}\,\rbrace}}
}
\def\midvert{\egroup\,\mid\,\bgroup}

\def\pmod#1{\text{ }(\text{mod } #1)\,}

\def\diag#1{\llbracket#1\rrbracket}

\def\K{\mathscr{K}}
\def\KK{\widehat{\mathscr{K}}}
\def\O{{\mathcal{O}}}
\def\OO{{\widehat{\mathcal{O}}}}

\newcommand\m{\mathfrak m}
\newcommand\p{\mathfrak p}

\DeclareMathOperator\noedge{\:\rlap{\hspace*{0.25em}/}\text{---}\:}
\DeclareMathOperator\Deg{Deg}

\DeclareMathAlphabet{\mathpzc}{OT1}{pzc}{m}{it}

\def\L{\mathscr{L}}
\newcommand\HH{\mathcal{H}}   
\newcommand\RR{\mathcal{R}}
\renewcommand\H[1][n]{\HH^\Lambda_{#1}}

\newcommand\HKK{\HH^\Lambda_n(\K)}
\newcommand\HK{\HH^\Lambda_n(\K)}
\newcommand\HO{\HH^\Lambda_n(\O)}

\newcommand\Hlam[1][]{\HH^{\gdom\blam}_n#1}

\newcommand\Gram[1][\blam]{\mathcal{G}^{#1}}

\def\psio{\psi^\O}

\def\Boo{B^\O}
\def\Do{D^\O}

\newcommand\yo{y^{\O}}
\newcommand\ylam[1][\blam]{y^{#1}_{\O}}
\def\dyo[#1]{y^{\<#1\>}}
\def\OBraid{\mathcal{B}^\O_r}
\def\eps{1+\rho_r(\bi)}
\def\spe{1-\rho_r(\bi)}
\newcommand\fo[1][\bi]{f_{#1}^\O}

\newcommand\R[1][n]{\RR^\Lambda_{#1}}

\newcommand\UnR[1][n]{\underline{\RR}^\Lambda_{#1}}
\newcommand\RO[1][n]{R_{#1}(\O)}

\def\SS{\mathcal{S}}           
\renewcommand\S[1][n]{\SS^\Lambda_{#1}}

\newcommand\UnD{\underline{D}}
\newcommand\UnGram{\underline{\mathcal G}^\blam}
\newcommand\UnS{\underline{S}}

\def\Pcal{\mathcal{P}}
\newcommand{\Klesh}[1][n]{\mathcal{K}^{\Lambda}_{#1}}
\newcommand{\Parts}[1][n]{\Pcal^{\Lambda}_{#1}}

\newcommand\blam{{\boldsymbol\lambda}}

\newcommand\balpha{{\boldsymbol\alpha}}
\newcommand\bmu{{\boldsymbol\mu}}

\newcommand\bnu{{\boldsymbol\nu}}

\def\i{\hat\imath}
\def\j{\hat\jmath}
\newcommand\bi{\mathbf{i}}
\newcommand\bj{\mathbf{j}}
\def\a{{\mathfrak a}}
\def\b{{\mathfrak b}}

\def\s{{\mathfrak s}}

\def\t{{\mathfrak t}}

\def\u{{\mathfrak u}}
\def\v{{\mathfrak v}}

\def\rest#1{{}_{\downarrow #1}}
\def\ilam{\bi^\blam}
\def\flam{\fo[\ilam]}

\def\tlam{\t^\blam}

\def\tllam{\t_\blam}
\def\tmu{\t^\bmu}

\def\Std{\mathop{\rm Std}\nolimits}
\def\DStd{\mathop{\rm DStd}\nolimits_e}
\def\SStd{\mathop{\rm Std}^2\nolimits}

\let\gedom=\trianglerighteq
\let\gdom=\vartriangleright
\let\Gdom=\blacktriangleright
\DeclareMathOperator\Gedom{\,{\underline{\kern-.1ex{\blacktriangleright}\kern-0.1ex}}\,}

\newcommand{\charge}{{\boldsymbol{\kappa}}}

\DeclareMathOperator\Shape{Shape}

\DeclareMathOperator\Dim{\dim_q}

\DeclareMathOperator\comp{comp}

\def\Cont{\mathscr{C}}
\newcommand\Dr[1][r]{\mathscr{D}_{#1}(\bi)}
\DeclareMathOperator\cont{c}

\DeclareMathOperator\rad{rad}

\DeclareMathOperator\res{res}

\def\Add{\mathscr A}

\usepackage{aliascnt}
\def\NewTheorem#1{%
  \newaliascnt{#1}{equation}
  \newtheorem{#1}[#1]{#1}
  \aliascntresetthe{#1}
  \expandafter\def\csname #1autorefname\endcsname{#1}
}
\def\equationautorefname~#1\null{(#1)\null}

\theoremstyle{plain}
\newcounter{MainTheore}
\newtheorem{MainTheorem}[MainTheore]{Theorem}

\swapnumbers
\numberwithin{equation}{section}
\NewTheorem{Assumption}
\NewTheorem{Proposition}
\NewTheorem{Theorem}
\NewTheorem{Corollary}
\NewTheorem{Conjecture}
\NewTheorem{Lemma}
\NewTheorem{Definition}
\theoremstyle{definition}
\theoremstyle{remark}
\NewTheorem{Remark}
\NewTheorem{Remarks}

\newaliascnt{Example}{equation}
\aliascntresetthe{Example}

\newenvironment{Example}%
  {\refstepcounter{Example}\trivlist
   \item[\hskip\labelsep\theequation.~\textbf{Example}\space]
   \ignorespaces
  }{\unskip\nobreak\hfil%
    \penalty50\hskip2em\hbox{}\nobreak\hfil$\Diamond$%
    \parfillskip=0pt\finalhyphendemerits=0\penalty-100\endtrivlist
}
\newaliascnt{Examples}{equation}
\aliascntresetthe{Examples}

\newenvironment{Examples}[1]%
  {\refstepcounter{Examples}\trivlist
   \item[\hskip\labelsep\theequation.~\textbf{Examples}\space]
   \ignorespaces #1\enumerate
   }{\endenumerate\unskip\nobreak\hfil%
    \penalty50\hskip2em\hbox{}\nobreak\hfil$\Diamond$%
    \parfillskip=0pt\finalhyphendemerits=0\penalty-100\endtrivlist
}

\newcounter{MyCase}
\numberwithin{MyCase}{equation}
\newcommand\Case[1]{\refstepcounter{MyCase}
  \medskip\noindent\textbf{Case \arabic{MyCase}.}
  \space\textit{#1.}\newline
}

\def\Step#1.{\medskip\noindent\textit{Step #1.}\space}

\let\bigast\ast

\def\Email#1{\email{\href{mailto:#1}{#1}}}

\title[Seminormal forms and quiver Hecke algebras]%
  {Seminormal forms and cyclotomic quiver Hecke algebras of type~$A$}
\subjclass[2000]{20G43, 20C08, 20C30}
\keywords{Cyclotomic Hecke algebras, quiver Hecke algebras}
\author{Jun Hu}
\address{School of Mathematics,
Beijing Institute of Technology,
Beijing, 100081, P.R.~China}
\address{School of Mathematics and Statistics F07, University of Sydney, NSW 2006, Australia}

\Email{junhu303@qq.com}
\author{Andrew Mathas}
\address{School of Mathematics and Statistics F07,
University of Sydney, NSW 2006, Australia}
\Email{andrew.mathas@sydney.edu.au}

\begin{document}

\begin{abstract}
This paper shows that the cyclotomic quiver Hecke algebras of type~$A$, and the
gradings on these algebras, are intimately related to the classical seminormal
forms.  We start by classifying all seminormal bases and then give an explicit
``integral'' closed formula for the Gram determinants of the Specht modules in
terms of the combinatorics associated with the KLR grading. We then use
seminormal forms to give a deformation of the KLR algebras of type~$A$. This
makes it possible to study the cyclotomic quiver Hecke algebras in terms of the
semisimple representation theory and seminormal forms. As an application we
construct a new distinguished graded cellular basis of the cyclotomic KLR
algebras of type $A$.
\end{abstract}

\maketitle
\setcounter{tocdepth}{1}
\tableofcontents

\section{Introduction}
  The quiver Hecke algebras are a remarkable family of algebras that were
  introduced independently by Khovanov and
  Lauda~\cite{KhovLaud:diagI,KhovLaud:diagII} and Rouquier~\cite{Rouq:2KM}.
  These algebras are attached to an arbitrary oriented quiver, they are
  $\Z$-graded and they categorify the negative part of the associated quantum
  group. Over a field, Brundan and Kleshchev showed that the cyclotomic quiver
  Hecke algebras of type~$A$, which are certain quotients of the quiver
  Hecke algebras of type~$A$, are isomorphic to the cyclotomic Hecke
  algebras of type~$A$.

  The quiver Hecke algebras have a homogeneous presentation by
  generators and relations. As a consequence they have well-defined
  integral forms. Unlike Hecke algebras, which are generically
  semisimple, the cyclotomic quiver Hecke algebras are typically not
  semisimple even over the rational field. As a result the cyclotomic
  quiver Hecke algebras are rarely isomorphic to the cyclotomic Hecke
  algebras over an arbitrary ring.

  The first main result of this paper shows that the cyclotomic quiver Hecke
  algebras of type~$A$ admit a one-parameter deformation. Moreover, this
  deformation is isomorphic to cyclotomic Hecke algebra defined over the
  corresponding ring. Before we can state this result we need
  some notation.

  Fix integers $n\ge0$ and $e>1$ and let $\Gamma_e$ be the oriented quiver with
  vertex set $I=\Z/e\Z$ and edges $i\to i+1$, for $i\in I$. Given $i\in I$ let
  $\i\ge0$ be the smallest non-negative integer such that $i=\i+e\Z$.  For each
  dominant weight $\Lambda$ for the corresponding Kac-Moody algebra
  $\mathfrak g(\Gamma_e)$, there exists a cyclotomic quiver Hecke algebra~$\R$
  and a cyclotomic Hecke algebra~$\H$. To each tuple $\bi\in I^n$ we associate
  the set of standard tableaux $\Std(\bi)$ with residue sequence $\bi$.
  All of these terms are defined in \autoref{S:Contents}.

  Like the cyclotomic quiver Hecke algebra, our deformation of~$\R$ is adapted
  to the choice of~$e$ through the choice of base ring~$\O$ that must be an
  \textit{$e$-idempotent subring} (\autoref{D:IdempotentSubring}).  This
  definition ensures that the cyclotomic Hecke algebras are semisimple over~$\K$, the
  field of fractions of~$\O$, and that $\HO\otimes_\O K$ is a cyclotomic quiver
  Hecke algebra whenever $K=\O/\m$, for $\m$ a maximal ideal of~$\O$. For
  $t\in\O$ and $d\in\Z$ let $[d]=[d]_t$ be the corresponding quantum integer, so
  that $[d]=(t^d-1)/(t-1)$ if $t\neq 1$ or $[d]=d$ if $t=1$.

  We can now state our first main result.

  \begin{MainTheorem}\label{Thm:KLRDeformation}
    Suppose that $1<e<\infty$ and that $(\O,t)$ is an $e$-idempotent
    subring of a field $\K$.  Then the algebra $\HO$ is generated as an
    $\O$-algebra by the elements $$\set{\fo|\bi\in
    I^n}\cup\set{\psio_r|1\le r<n}\cup\set{\yo_r|1\le r\le n}$$ subject
    only to the following relations:
    {\setlength{\abovedisplayskip}{2pt}
     \setlength{\belowdisplayskip}{1pt}
     \begin{align*}
         \prod_{\substack{1\le l\le\ell\\\kappa_l\equiv i_1\pmod e}}
              (\yo_1-[\kappa_l-\i_1])\fo=0,
     \end{align*}
     \begin{align*}
        \fo \fo[\bj] = \delta_{\bi\bj} \fo, \qquad
           {\textstyle\sum_{\bi \in I^n}} \fo= 1, \qquad
            \yo_r \fo = \fo \yo_r,
     \end{align*}
     \begin{align*}
      \psio_r \fo&= \fo[s_r{\cdot}\bi] \psio_r,
      &\yo_r \yo_s &= \yo_s \yo_r,\\
      \psio_r \yo_{r+1} \fo&=(\yo_r\psio_r+\delta_{i_ri_{r+1}})\fo,&
      \yo_{r+1}\psio_r\fo&=(\psio_r \yo_r+\delta_{i_ri_{r+1}})\fo,\\
      \psio_r \yo_s  &= \yo_s \psio_r,&&\text{if }s \neq r,r+1,\\
      \psio_r \psio_s &= \psio_s \psio_r,&&\text{if }|r-s|>1,
    \end{align*}
    \begin{align*}
      (\psio_r)^2\fo &= \begin{cases}
        (\dyo[\eps]_r-\yo_{r+1})(\dyo[\spe]_{r+1}-\yo_r)\fo,&
                   \text{if }i_r\leftrightarrows i_{r+1},\\
        (\dyo[\eps]_r-\yo_{r+1})\fo,& \text{if }i_r\rightarrow i_{r+1},\\
        (\dyo[\spe]_{r+1}-\yo_r)\fo,& \text{if }i_r\leftarrow i_{r+1},\\
           0,&\text{if }i_r=i_{r+1},\\
           \fo,&\text{otherwise,}
    \end{cases}
    \end{align*}
    and where
    $\big(\psio_r\psio_{r+1}\psio_r-\psio_{r+1}\psio_r\psio_{r+1}\big)\fo$
    is equal to
    \begin{align*}
        \begin{cases}
            (\dyo[\eps]_r+\dyo[\eps]_{r+2}-
            \dyo[\eps]_{r+1}-\dyo[\spe]_{r+1})\fo,
                    &\text{if }i_{r+2}=i_r\rightleftarrows i_{r+1},\\
            -t^{\eps}\fo,&\text{if }i_{r+2}=i_r\rightarrow i_{r+1},\\
            \fo,&\text{if }i_{r+2}=i_r\leftarrow i_{r+1},\\
            0,&\text{otherwise},
        \end{cases}
    \end{align*}
  }
  where $\rho_r(\bi)=\i_r-\i_{r+1}$ and $\dyo[d]_r=t^d\yo_r+[d]$, for $d\in\Z$.
  \end{MainTheorem}

  As we explain in \autoref{C:Largee} by taking $e$ large enough  this
  result also applies when $e=0$. The appendix gives a direct
  treatment of this case.

  \newcommand\Op{\mathbb{Z}_{(p)}}

  To help the reader interpret \autoref{Thm:KLRDeformation} we include
  the following special case of this result that gives a new
  presentation of the group algebra of the symmetric group over the
  ring~$\Op$, where $p$ is an integer prime and $\Op$ is the
  localisation of~$\Z$ at the prime ideal~$p\Z$.

  \begin{Corollary}\label{MainCorPe1}
    Suppose that $e=p$ be an odd prime
    number and let $I=\Z/p\Z$ and $\Lambda=\Lambda_0$. Then the group
    algebra $\Op\Sym_n$ is generated as an $\Op$-algebra by the elements
    $$\set{\fo|\bi\in I^n}\cup\set{\psio_r|1\le r<n}\cup\set{\yo_s|1\le s\le n}$$
    subject only to the relations:
    {\setlength{\abovedisplayskip}{2pt}
     \setlength{\belowdisplayskip}{1pt}
     \begin{align*}
          (\yo_1)^{(\Lambda,\alpha_{i_1})}e(\bi)=0,\qquad
          \fo \fo[\bj] = \delta_{\bi\bj} \fo, \qquad
          {\textstyle\sum_{\bi \in I^n}} \fo= 1, \qquad
          \yo_r \fo = \fo \yo_r,
     \end{align*}
     \begin{align*}
      \psio_r \fo&= \fo[s_r{\cdot}\bi] \psio_r,
      &\yo_r \yo_s &= \yo_s \yo_r,\\
      \psio_r \yo_{r+1} \fo&=(\yo_r\psio_r+\delta_{i_ri_{r+1}})\fo,&
      \yo_{r+1}\psio_r\fo&=(\psio_r \yo_r+\delta_{i_ri_{r+1}})\fo,\\
      \psio_r \yo_s  &= \yo_s \psio_r,&&\text{if }s \neq r,r+1,\\
      \psio_r \psio_s &= \psio_s \psio_r,&&\text{if }|r-s|>1,
    \end{align*}
    \begin{align*}
      (\psio_r)^2\fo &= \begin{cases}
        (\yo_r-\yo_{r+1})\fo,& \text{if }i_r\rightarrow i_{r+1}\neq0,\\
        (\yo_{r}+p-\yo_{r+1})\fo,& \text{if }i_r\rightarrow i_{r+1}=0,\\
        (\yo_{r+1}-\yo_r)\fo,& \text{if }0\neq i_r\leftarrow i_{r+1},\\
        (\yo_{r+1}+p-\yo_r)\fo,& \text{if } 0=i_r\leftarrow i_{r+1},\\
        0,&\text{if }i_r=i_{r+1},\\
           \fo,&\text{otherwise,}
    \end{cases}
    \end{align*}
    \begin{align*}
      \big(\psio_r\psio_{r+1}\psio_r-\psio_{r+1}\psio_r\psio_{r+1}\big)\fo=
        \begin{cases}
          -\fo,&\text{if }i_{r+2}=i_r\rightarrow i_{r+1},\\
          \fo,&\text{if }i_{r+2}=i_r\leftarrow i_{r+1},\\
          0,&\text{otherwise},
        \end{cases}
    \end{align*}
    }
    for all admissible $r,s$ and $\bi\in I^n$.
  \end{Corollary}

  Except for the cyclotomic relation and the last two relations (that
  is, the quadratic relations and the braid relations for
  $\psio_1,\dots\psio_{n-1}$), all of the relations in
  \autoref{Thm:KLRDeformation} coincide with the corresponding
  KLR-relations in~$\R$. Interestingly, only the ``Jucys-Murphy like
  elements'' $\yo_r$ need to be modified in order to define a
  deformation of~$\R$.  Over a field $K=\O/\m$, the presentation in
  \autoref{Thm:KLRDeformation} collapses to give the KLR algebra $\R$
  because the definition of an idempotent subring ensures that
  $t^{1+\rho_r(\bi)}\otimes1_K=1$ and
  $y_r^{\<1\pm\rho_r(\bi)\>}\otimes1_K=\yo_r\otimes1_K$, for $1\le r\le
  n$.

  As a first application of \autoref{Thm:KLRDeformation},
  \autoref{C:GeneralNilpotence} gives what appears to be tight upper bounds
  on the nilpotency indices of the elements $y_1,\dots,y_n$ in the
  cyclotomic quiver Hecke algebras of type~$A$. Previously such a result
  was known only in the special case of the linear quiver or, equivalently,
  when $e=0$.

  To prove \autoref{Thm:KLRDeformation} we work almost entirely inside the
  semisimple representation theory of the cyclotomic Hecke algebras~$\H$. We
  show that definition of the quiver Hecke algebra~$\R$, and its grading, is
  implicit in \textit{Young's seminormal form}. With hindsight, using the
  perspective afforded by this paper, it is not too much of an exaggeration to
  say that Murphy could have discovered the cyclotomic quiver Hecke algebras in
  1983 soon after writing his paper on the Nakayama conjecture~\cite{M:Nak}.

  Our proof of \autoref{Thm:KLRDeformation} gives another explanation
  for the KLR relations and a more conceptual proof of one direction in
  Brundan and Kleshchev's isomorphism theorem~\cite{BK:GradedKL} (see
  \autoref{T:BKiso}). In fact, we give a new proof of the
  Brundan-Kleshchev isomorphism theorem by using the
  Ariki-Brundan-Kleshchev categorification
  theorem~\cite{Ariki:class,BK:GradedDecomp} to bound the dimension of
  the algebras defined by the presentation in
  \autoref{Thm:KLRDeformation}.

  For the algebras of type~$A$ the authors have constructed a graded
  cellular basis $\set{\psi_{\s\t}|(\s,\t)\in\SStd(\Parts)}$
  for~$\R$~\cite{HuMathas:GradedCellular}. Here $\SStd(\Parts)$ is the
  set of all pairs of standard tableaux of the same shape, where the
  shape is a multipartition of~$n$.  The element~$\psi_{\s\t}$ is
  homogeneous of degree $\deg_e\s+\deg_e\t$, where
  $\deg_e\map{\Std(\Parts)}\Z$ is the combinatorial degree function
  introduced by Brundan, Kleshchev and Wang~\cite{BKW:GradedSpecht}.
  Li~\cite{Li:PhD} has shown that $\{\psi_{\s\t}\}$ is a graded cellular
  basis of~$\R$ over an arbitrary ring. In particular, the KLR algebra
  $\R$ is always free of rank $\dim\H(K)$, for~$K$ a field.

  One of the problems with the basis $\{\psi_{\s\t}\}$ is that, because
  the KLR generators~$\psi_r$, for $1\le r<n$, do not satisfy the braid
  relations, the basis elements $\psi_{\s\t}$ depend upon a choice of
  reduced expression for certain permutations $d(\s),d(\t)\in\Sym_n$
  associated with the tableaux~$\s$ and $\t$; see
  \autoref{S:tableaux}. As a consequence, the results
  of~\cite{HuMathas:GradedCellular} constructs different $\psi$-bases
  for different choices of reduced expressions for the elements
  of~$\Sym_n$. The different $\psi$-bases constructed in this way are
  closely related and it would be advantageous to be able to make a
  canonical choice of basis, however, until now it has not been clear
  how to do this.

  Fix a modular system $(\K,\O,K)$ as in \autoref{Chap:BBasis} and
  consider the corresponding cyclotomic Hecke algebras $(\HK,\HO,\H)$,
  where $\H=\H(K)$.  The algebra $\HK$ is semisimple and has a
  seminormal basis $\set{f_{\s\t}|(\s,\t)\in\SStd(\Parts)}$, $\HO$ is a
  free $\O$-subalgebra of $\HK$ and $\H\cong\HO\otimes_\O K$. Finally,
  we recall that the set~$\SStd(\Parts)$ comes equipped with a naturally
  partial order~$\Gdom$; see \autoref{S:tableaux}.

  \begin{MainTheorem}\label{Thm:BBasis}
    Suppose that $K$ is a field of characteristic zero and that
    $(\s,\t)\in\SStd(\Parts)$. Then there is a unique
    element $\Boo_{\s\t}\in\HO$ such that
    \begin{enumerate}
      \item $\Boo_{\s\t}=f_{\s\t}
            +\sum_{(\u,\v)\Gdom(\s,\t)}p^{\s\t}_{\u\v}(x^{-1})f_{\u\v},$
          where if $(\u,\v)\Gdom(\s,\t)$ then $p^{\s\t}_{\u\v}(x)\in xK[x]$
          and $\deg p^{\s\t}_{\u\v}(x)\le\tfrac12(\deg\u-\deg\s+\deg\v-\deg\t).$
      \item $\Boo_{\s\t}\otimes_\O1_K=B_{\s\t}'+C_{\s\t}$, where
        $B_{\s\t}'$ is homogeneous of degree $\deg\s+\deg\t$ and~$C_{\s\t}$
        is a sum of homogeneous terms of degree strictly larger than
        $\deg B_{\s\t}$.
    \end{enumerate}
    Moreover,
    $\set{B_{\s\t}'|(\s,\t)\in\Std(\Parts)}$ is a graded cellular basis
    of~$\H$ and
    $$B_{\s\t}'=\psi_{\s\t}+\sum_{(\u,\v)\Gdom(\s,\t)}r_{\u\v}\psi_{\u\v},$$
    for some $r_{\u\v}\in K$.
  \end{MainTheorem}

  There is a similar graded cellular basis of $\HO$ when K is a field of
  positive characteristic, however, its' description is more complicated
  because the corresponding polynomial $p^{\s\t}_{\u\v}(x)$ do not
  necessarily satisfy the degree bound in \autoref{Thm:BBasis}(a).  The
  construction of the $B$-basis is reminiscent of the Kazhdan-Lusztig
  basis~\cite{KL}.  The $B$-basis can depend on the choice of
  multicharge.

  As remarked above, the basis element $\psi_{\s\t}$ depends upon
  choices of reduced expression for the permutations
  $d(\s),d(\t)\in\Sym_n$, yet for any choice $B_{\s\t}$ is equal to
  $\psi_{\s\t}$ plus a linear combination of more dominant terms by
  \autoref{Thm:BBasis}. The $B$-basis elements depend only on the
  indexing tableaux, and not on choices of reduced expressions. In
  this sense, the $B$-basis corrects for a deficiency in the
  definition of the $\psi$-bases.

  To prove the two theorems above, we define a seminormal basis of a semisimple
  Hecke algebra to be a basis of~$\H$ of simultaneous eigenvectors for the
  Gelfand-Zetlin subalgebra of~$\H$. Seminormal bases are classical objects
  that are ubiquitous in the literature, having been rediscovered many times
  since were first introduced for the symmetric groups by Young
  in~1900~\cite{QSAI}.

  Seminormal basis elements are a basis of eigenvectors for the action
  of the Jucys-Murphy elements on the regular representation of~$\H(K)$.
  Eigenbases are, of course, only unique up to scalar multiplication.
  This paper starts by introducing \textit{seminormal coefficient
  systems} that gives a combinatorial framework for describing the
  structure constants of the algebra in terms of the choice of
  eigenvectors. The real surprise is that seminormal coefficient systems
  encode the KLR grading.

  The close connections between the semisimple representation theory and the KLR
  gradings is made even more explicit in the third main result of this paper
  that gives a closed formula for the Gram determinants of the semisimple
  Specht modules of these algebras. Closed formulas for these determinants
  already exist in the literature~\cite{AMR,JM:det,JM:Schaper,JM:cyc-Schaper},
  however, all of these formulas describe these determinants as rational functions
  (or rational numbers in the degenerate case).  The theorem below gives the
  first \textit{integral} formula for these determinants.

  In order to state the closed integral formulas for the Gram determinant of the
  Specht module $S^\blam$, for a multipartition $\blam$ define
  $$\deg_e(\blam)=\sum_{\t\in\Std(\blam)}\deg_e(\t)\in\Z,$$
  where $\Std(\blam)$ is the set of standard $\blam$-tableaux.  Let
  $\Phi_e(t)\in\Z[t]$ be the $e$th cyclotomic polynomial for $e>1$.
  We prove the following (see \autoref{T:GramDetFactorization} for a more
  precise statement).

  \begin{MainTheorem}\label{Thm:IntegralDet}
    Suppose that $\H$ is a semisimple cyclotomic Hecke algebra over~$\Q(t)$,
    with Hecke parameter~$t$. Let $\blam$ be a multipartition of $n$. Then
    the Gram determinant of the Specht module $S^\blam$ is equal to
    $$t^N\prod_{e>1}\Phi_e(t)^{\deg_e(\blam)},$$
    for a known integer $N$. In particular, $\deg_e(\blam)\ge0$, for all
    $e\in\{0,2,3,4,\dots\}$.
  \end{MainTheorem}

  As the integers $\deg_e(\blam)$ are defined combinatorially, it should be
  possible to give a purely combinatorial proof that $\deg_e(\blam)\ge0$. In
  \autoref{S:GramDet} we give two representation theoretic proofs of this
  result. The first proof is elementary but not very enlightening. The second
  proof uses deep positivity properties of the graded decomposition numbers
  of~$\H(\C)$ to show that the tableaux combinatorics of~$\H$ provides a
  framework for giving purely combinatorial formulas for the graded dimensions
  of the simple $\H$-modules and for the graded decomposition numbers of~$\H$.
  Interestingly, we show that there is a close connection between the graded
  dimensions of the simple $\H$-modules and the graded decomposition numbers
  for~$\H$.  Note that in characteristic zero, the graded decomposition numbers
  of~$\H$ are parabolic Kazhdan-Lusztig polynomials of
  type~$A$~\cite{BK:GradedDecomp}, so our results show that the tableaux
  combinatorics leads to combinatorial formulas for these polynomials.
  Unfortunately, we are only able to prove that such formulas exist and we are
  not able to make them explicit or to show that they are canonical in any way.

  The outline of this paper is as follows. \autoref{Chap:HeckeAlgebras} defines
  the cyclotomic Hecke algebras of type~$A$, giving a uniform presentation for
  the degenerate and non-degenerate algebras. Previously these algebras have
  been treated separately in the literature. We then recall the basic results
  about these algebras that we need from the literature, including Brundan and
  Kleshchev's isomorphism theorem~\cite{BK:GradedKL}.
  \autoref{Chap:Seminormal} develops the theory of seminormal bases for these
  algebras in full generality. We completely classify the seminormal bases
  of~$\H$ and then use them to prove \autoref{Thm:IntegralDet}, thus
  establishing a link between the semisimple representation theory of~$\H$ and
  the quiver Hecke algebra~$\R$. Using this we prove the existence of
  combinatorial formulas for the graded dimensions of the simple modules and the graded
  decomposition numbers of~$\H$. In \autoref{Chap:QuiverHeckeAlgebras} we use
  the theory of seminormal forms to construct a deformation of the cyclotomic
  quiver Hecke algebras of type~$A$, culminating with the proof of
  \autoref{Thm:KLRDeformation}.  \autoref{Chap:PsiBases} builds on
  \autoref{Thm:KLRDeformation} to give a quicker construction of the graded
  cellular basis of~$\H(K)$, over a field~$K$, which was one of the main results
  of~\cite{HuMathas:GradedCellular}. Finally, in \autoref{Chap:BBasis} we use
  \autoref{Thm:KLRDeformation} to show that $\H(K)$ has the distinguished graded
  cellular basis described in \autoref{Thm:BBasis}.

\section{Cyclotomic Hecke algebras}\label{Chap:HeckeAlgebras} This chapter
  defines the cyclotomic Hecke and quiver Hecke algebras of type~$A$ and it
  introduces some of the basic machinery that we need for understanding these
  algebras.  We give a new presentation for the \textit{cyclotomic Hecke
  algebras of type~$A$}, which simultaneously captures the degenerate and
  non-degenerate cyclotomic Hecke algebras that currently appear in the
  literature, and then we recall the results from the literature that we need,
  including Brundan and Kleshchev's isomorphism
  theorem~\cite{BK:GradedKL}.

  \subsection{Quiver combinatorics}\label{S:Quiver}
  Fix an integer $e\in\{0,2,3,4\dots\}$ and let $\Gamma_e$ be the oriented quiver with
  vertex set $I=\Z/e\Z$ and edges $i\longrightarrow i+1$, for $i\in I$. If
  $i,j\in I$ and~$i$ and~$j$ are not connected by an edge in~$\Gamma_e$ then we write
  $i\noedge j$.

  To the quiver
  $\Gamma_e$ we attach the Cartan matrix
  $(c_{ij})_{i,j\in I}$, where
  $$c_{i,j}=\begin{cases} 2,&\text{if } i=j,\\
    -1,&\text{if $i\rightarrow j$ or $i\leftarrow j$},\\
    -2,&\text{if }i\leftrightarrows j,\\
    0,&\text{otherwise},
  \end{cases}$$
  Let $\widehat{\mathfrak{sl}}_e$ be the corresponding Kac-Moody
  algebra~\cite{Kac} with fundamental weights $\set{\Lambda_i|i\in I}$, positive
  weight lattice $P^+_e=\sum_{i\in I}\N\Lambda_i$ and  positive root lattice
  $Q^+=\bigoplus_{i\in I}\N\alpha_i$. Let $(\cdot,\cdot)$ be the bilinear form
  determined by
  $$(\alpha_i,\alpha_j)=c_{ij}\qquad\text{and}\qquad
          (\Lambda_i,\alpha_j)=\delta_{ij},\qquad\text{for }i,j\in I.$$
  More details can be found, for example, in ~\cite[Chapter~1]{Kac}.

  Fix, once and for all, a \textbf{multicharge}
  $\charge=(\kappa_1,\dots,\kappa_\ell)\in\Z^\ell$ that is a sequence of integers such that
  if $e\ne0$ then $\kappa_l-\kappa_{l+1}\ge n$ for $1\le l<\ell$. Define
  $\Lambda=\Lambda_e(\charge)=\Lambda_{\bar\kappa_1}+\dots+\Lambda_{\bar\kappa_\ell}$,
  where $\bar\kappa=\kappa\pmod e$. Equivalently, $\Lambda$ is the unique
  element of $P^+_e$ such that
  \begin{equation}\label{E:multicharge}
  (\Lambda,\alpha_i) = \#\set{1\le l\le\ell | \kappa_l\equiv i\pmod e},
     \qquad\text{ for all }i\in I.
  \end{equation}
  All of the bases for the modules and algebras in this paper depend implicitly
  on the choice of~$\charge$ even though the algebras themselves depend only on~$\Lambda$.

  \subsection{Cyclotomic Hecke algebras}
  This section defines the cyclotomic Hecke algebras of type~$A$ and explains
  the connection between these algebras and the degenerate and non-degenerate
  Hecke algebras of type $G(\ell,1,n)$.

  Fix an integral domain $\O$ that contains an invertible element $\xi\in\O^\times$.

\begin{Definition}\label{D:HeckeAlgebras}
  Fix integers $n\ge0$ and $\ell\ge1$. Then the \textbf{cyclotomic Hecke algebra of
  type~$A$} with Hecke parameter $\xi\in\O^\times$ and cyclotomic parameters
  $Q_1,\dots,Q_\ell\in \O$ is the unital associative $\O$-algebra
  $\HH_n=\HH_n(\O,\xi,Q_1,\dots,Q_\ell)$
  with generators $L_1,\dots,L_n$, $T_1,\dots,T_{n-1}$ that are subject to the
  relations
  \begin{align*}
      \prod_{l=1}^\ell(L_1-Q_l)&=0,  &
      (T_r+1)(T_r-\xi )&=0,  \\
      L_rL_t&=L_tL_r, &
    T_rT_s&=T_sT_r &\text{if }|r-s|>1,\\
    T_sT_{s+1}T_s&=T_{s+1}T_sT_{s+1}, &
    T_rL_t&=L_tT_r,&\text{if }t\ne r,r+1,\\
    \span\span L_{r+1}(T_r-\xi+1)=T_rL_r+1,\span\span\span
  \end{align*}
  where $1\le r<n$, $1\le s<n-1$ and $1\le t\le n$.
\end{Definition}

\begin{Remark}\label{R:NewPresentation}
  If $\xi=1$ then, by definition, $\HH_n$ is a degenerate cyclotomic Hecke algebra
  of type $G(\ell,1,n)$. If $\xi\ne1$ then $\HH_n$ is (isomorphic to) an integral
  cyclotomic Hecke algebra of type $G(\ell,1,n)$. To see this define
  $L_k'=(\xi-1)L_k+1$, for $1\le k\le n$, and observe that $\HH_n$ is generated by
  $L_1',T_1,\dots,T_{n-1}$ subject to the usual relations for these algebras as
  originally defined by Ariki and Koike~\cite{AK}. It is now easy to verify our
  claim. For each $1\leq m\leq n$, an eigenvector for $L_m$ of eigenvalue $[k]_{\xi}$
  is the same as an eigenvector for $L'_m$ of eigenvalue $\xi^k$.
  The presentation of~$\HH_n$ in \autoref{D:HeckeAlgebras} unifies the
  definition of the `degenerate' and `non-degenerate' Hecke algebras, which
  corresponds to the cases where $\xi=1$ or $\xi\ne1$, respectively.
\end{Remark}

Let $\Sym_n$ be the \textbf{symmetric group} on~$n$ letters. For $1\le r<n$ let
$s_r=(r,r+1)$ be the corresponding simple transposition. Then
$\{s_1,\dots,s_{n-1}\}$ is the standard set of Coxeter generators for~$\Sym_n$.
A \textbf{reduced expression} for $w\in\Sym_n$ is a word $w=s_{r_1},\dots
s_{r_k}$ with $k$ minimal and $1\le r_j<n$ for $1\le j\le k$. If $w=s_{r_1}\dots
s_{r_k}$ is reduced then set $T_w=T_{r_1}\dots T_{r_k}$. Then
$T_w$ is independent of the choice of reduced expression since the braid
relations hold in~$\HH_n$. It follows arguing as in~\cite[Theorem~3.3]{AK} that
$\HH_n$ is free as an~$\O$-module with basis
$$\set{L_1^{a_1}\dots L_n^{a_n} T_w|0\le a_1,\dots,a_n<\ell\text{ and }w\in\Sym_n}.$$
Consequently, $\HH_n$ is free as an $\O$-module of rank $\ell^n n!$, which is the
order of the complex reflection group of type $G(\ell,1,n)$.

We now restrict our attention to the case of \textit{integral} cyclotomic
parameters. Recall that for any integer $d$ and $t\in\O$ the
\textbf{quantum integer} $[d]_t$ is
  $$[d]_t=\begin{cases}
    \phantom{-}1+t+\dots+t^{d-1},&\text{if }d\ge0,\\
    -(t^{-1}+t^{-2}+\dots+t^d),&\text{if }d<0.
  \end{cases}$$
When $t$ is understood we write $[d]=[d]_t$. Set
$[d]_t^!=[d]^!=[1][2]\dots[d]$ when $d>0$.

An \textbf{integral cyclotomic Hecke algebra} is a cyclotomic Hecke algebra
$\HH_n$ with cyclotomic parameters of the form $Q_r=[\kappa_r]_\xi$, for
$\kappa_1,\dots,\kappa_\ell\in\Z$. The sequence of integers
$\charge=(\kappa_1,\dots,\kappa_\ell)\in\Z^\ell$
is the \textbf{multicharge} of $\HH_n$.

Translating the Morita equivalence theorems of \cite[Theorem~1.1]{DM:Morita} and
\cite[Theorem~5.19]{BK:HigherSchurWeyl} into the current setting, every cyclotomic
Hecke algebras of type~$A$ is Morita equivalent to a direct sum of tensor
products of integral cyclotomic Hecke algebras. Therefore, there is no loss of
generality in restricting our attention to the integral cyclotomic Hecke
algebras of type~$A$.

Recall that $\Lambda\in P^+_e$ and that we have fixed an integer
$e\in\{0,2,3,4,\dots\}$. We assume that $\O$ contains an invertible
element $\xi\in\O^\times$ such that $[e]_\xi=0$ if $e>0$ and
$[f]_\xi\ne0$ for all $f\ge1$ if $e=0$. Hence, either:
\begin{enumerate}
  \item $\xi=1$ and $e$ is prime and equal to the characteristic of~$\O$,
  \item $e>0$ and $\xi$ is a primitive $e$th root of unity, or,
  \item $e=0$ and $\xi$ is not a root of unity.
\end{enumerate}
In addition, fix a multicharge $\charge$ so that $\Lambda=\Lambda_e(\charge)$ as
in \autoref{E:multicharge}.

Let $\H=\H(\O)$ be the integral cyclotomic Hecke algebra
$\HH_n(\O,\xi,\charge)$. Using the definitions it is easy to see that, up to
isomorphism, $\H$ depends only on~$\xi$ and~$\Lambda$. In fact, by
\autoref{T:BKiso} below, it depends only on~$e$ and~$\Lambda$. Nonetheless,
many of the constructions that follow, particularly the definitions of bases,
depend upon the choice of~$\charge$.

\subsection{Graded algebras and cellular bases}\label{S:GradedCellular}

This section recalls the definitions and results from the representation theory
of (graded) cellular algebras that we need.

Let $A$ be a unital associative $\O$-algebra that is free and of finite rank as
an $\O$-module.  In this paper a \textbf{graded module} will always mean a
$\Z$-graded module. That is, an $\O$-module $M$ that has a decomposition
$M=\bigoplus_{n\in\Z}M_d$ as an $\O$-module. If  $m\in M_d$, for $d\in\Z$, then
$m$ is \textbf{homogeneous} of \textbf{degree} $d$ and we set $\deg m=d$. If $M$
is a graded $\O$-module and $s\in\Z$ let $M\<s\>$ be the graded $\O$-module
obtained by shifting the grading on $M$ up by $s$; that is, $M\<s\>_d=M_{d-s}$,
for $d\in\Z$.

Similarly a \textbf{graded algebra} is a  unital associative
$\O$-algebra $A=\bigoplus_{d\in\Z}A_d$ that is a graded
$\O$-module such that $A_dA_e\subseteq A_{d+e}$, for all
$d,e\in\Z$. It follows that $1\in A_0$ and that $A_0$ is a graded
subalgebra of $A$.  A graded (right) $A$-module is a graded
$\O$-module $M$ such that $\underline{M}$ is an
$\underline{A}$-module and $M_dA_e\subseteq M_{d+e}$, for all
$d,e\in\Z$, where $\underline{M}$ and $\underline{A}$ mean forgetting
the $\Z$-grading structures on $M$ and $A$ respectively. Graded submodules,
graded left $A$-modules and so on are all defined in the obvious way.


The following definition extends Graham and Lehrer's~\cite{GL} definition of
cellular algebras to the graded setting.

\begin{Definition}[Graded cellular algebras~\cite{GL,HuMathas:GradedCellular}]
  \label{graded cellular def}\label{D:cellular}
  Suppose that $A$ is an $\O$-algebra that is free of finite rank
  over $\O$. A \textbf{cell datum} for $A$ is an ordered triple
  $(\Pcal,T,C)$, where $(\Pcal,\gdom)$ is the \textbf{weight poset}, $T(\lambda)$
  is a finite set for $\lambda\in\Pcal$, and
  $$C\map{\coprod_{\lambda\in\Pcal}T(\lambda)\times T(\lambda)}A;
     (\s,\t)\mapsto c_{\s\t},
  $$
  is an injective function such that:
  \begin{enumerate}
    \item[(GC$_1$)] $\set{c_{\s\t}|\s,\t\in T(\lambda) \text{ for } \lambda\in\Pcal}$ is an
      $\O$-basis of $A$.
    \item[(GC$_2$)] If $\s,\t\in T(\lambda)$, for some $\lambda\in\Pcal$, and $a\in A$ then
    there exist scalars $r_{\t\v}(a)$, which do not depend on $\s$, such that
      $$c_{\s\t} a=\sum_{\v\in T(\lambda)}r_{\t\v}(a)c_{\s\v}\pmod
      {A^{\gdom\lambda}},$$
      where $A^{\gdom\lambda}$ is the $\O$-submodule of $A$ spanned by
      $\set{c_{\a\b}|\mu\gdom\lambda\text{ and }\a,\b\in T(\mu)}$.
    \item[(GC$_3$)] The $\O$-linear map $\ast\map AA$ determined by
      $(c_{\s\t})^\ast=c_{\t\s}$, for all $\lambda\in\Pcal$ and
      all $\s,\t\in T(\lambda)$, is an anti-isomorphism of $A$.
  \end{enumerate}
  A \textbf{cellular algebra} is an algebra that has a
  cell datum. If $A$ is a cellular algebra with cell datum $(\Pcal,T,C)$ then the basis
  $\set{c_{\s\t}|\lambda\in\Pcal\text{ and } \s,\t\in T(\lambda}$ is a
  \textbf{cellular basis} of~$A$ with $*$ its cellular algebra anti-automorphism.

  If, in addition, $A$ is a $\Z$-graded algebra then a \textbf{graded cell
  datum} for $A$ is a cell datum $(\Pcal,T,C)$ together with a \textit{degree
  function}
  $$\deg\map{\coprod_{\lambda\in\Pcal}T(\lambda)}\Z$$
  such that
  \begin{enumerate}
    \item[(GC$_d$)]the element $c_{\s\t}$ is homogeneous of degree $\deg
    c_{\s\t}=\deg(\s)+\deg(\t)$, for all $\lambda\in\Pcal$ and
    $\s,\t\in T(\lambda)$.
  \end{enumerate}
  In this case, $A$ is a \textbf{graded cellular algebra} with
  \textbf{graded cellular basis} $\{c_{\s\t}\}$.
\end{Definition}

Fix a (graded) cellular algebra $A$ with graded cellular basis
$\{c_{\s\t}\}$. If $\lambda\in\Pcal$ then the graded \textbf{cell
module} is the $\O$-module $C^\lambda$ with basis
$\set{c_\t|\t\in T(\lambda)}$ and with $A$-action
$$c_{\t} a=\sum_{\v\in T(\lambda)}r_{\t\v}(a)c_{\v},$$ where the scalars $r_{\t\v}(a)\in \O$
are the same scalars appearing in (GC$_2)$. One of the key
properties of the graded cell modules is that by
\cite[Lemma~2.7]{HuMathas:GradedCellular} they come equipped with a
homogeneous bilinear form $\<\ ,\ \>$ of degree zero that is
determined by the equation
\begin{equation}\label{E:CellularInnerProduct}
  \<c_\t, c_\u\>c_{\s\v}\equiv c_{\s\t}c_{\u\v}\pmod{A^{\gdom\lambda}},
\end{equation}
for $\s, \t, \u, \v\in T(\lambda)$.
The radical of this form
$$\rad C^\lambda=\set{x\in C^\lambda|\<x,y\>=0
            \text{ for all }y\in C^\lambda}$$
is a graded $A$-submodule of $C^\lambda$ so that $D^\lambda=C^\lambda/\rad
C^\lambda$ is a graded $A$-module. It is shown in
\cite[Theorem~2.10]{HuMathas:GradedCellular} that
$$\set{D^\lambda\<k\>|\lambda\in\Pcal, D^\lambda\ne0 \text{ and }k\in\Z}$$
is a complete set of pairwise non-isomorphic irreducible (graded) $A$-modules
when~$\O$ is a field.

\subsection{Multipartitions and tableaux}\label{S:tableaux}
A \textbf{partition} of~$d$ is a weakly decreasing sequence
$\lambda=(\lambda_1,\lambda_2,\dots)$ of non-negative integers such that
$|\lambda|=\lambda_1+\lambda_2+\dots=d$. An \textbf{$\ell$-multipartition} of~$n$ is an
$\ell$-tuple $\blam=(\lambda^{(1)},\dots,\lambda^{(\ell)})$ of partitions such
that $|\lambda^{(1)}|+\dots+|\lambda^{(\ell)}|=n$. We identify the
multipartition $\blam$ with its \textbf{diagram} that is the set of \textbf{nodes}
$\diag\blam=\set{(l,r,c)|1\le c\le \lambda^{(l)}_r\text{ for }1\le l\le\ell}$,
which we think of as an ordered $\ell$-tuple of arrays of boxes in the plane. For
example, if $\blam=(3,1^2|2,1|3,2)$ then
$$\diag\blam=\TriDiagram(3,1,1|2,1|3,2).$$
In this way we talk of the \textbf{rows}, \textbf{columns} and
\textbf{components} of~$\blam$.

Given two nodes $\alpha=(l,r,c)$ and $\beta=(l',r',c')$ then $\beta$ is
\textbf{below} $\alpha$, or $\alpha$ is \textbf{above} $\beta$, if
$(l,r,c)<(l',r',c')$ in the lexicographic order.

The set of multipartitions of~$n$ becomes a poset ordered by
\textbf{dominance} where~$\blam$ dominates~$\bmu$, or $\blam\gedom\bmu$,
if
$$\sum_{k=1}^{l-1}|\lambda^{(k)}|+\sum_{j=1}^i\lambda^{(l)}_j
     \ge\sum_{k=1}^{l-1}|\mu^{(k)}|+\sum_{j=1}^i\mu^{(l)}_j,$$
for $1\le l\le\ell$ and $i\ge1$. If $\blam\gedom\bmu$ and $\blam\ne\bmu$
then write $\blam\gdom\bmu$. Let $\Parts=(\Parts,\gedom)$ be the poset of
multipartitions of~$n$ ordered by dominance.

Fix a multipartition $\blam$. Then a \textbf{$\blam$-tableau} is a bijective map
$\t\map{\diag\blam}\{1,2,\dots,n\}$, which we identify with a labelling of
$\diag\blam$ by $\{1,2,\dots,n\}$. For example,
$$\Tritab({1,2,3},{4},{5}|{6,7},{8}|{9,10,11},{12,13})\quad\text{and}\quad
  \Tritab({9,12,13},{10},{11}|{6,8},{7}|{1,3,5},{2,4})$$
are both $\blam$-tableaux when $\blam=(3,1^2|2,1|3,2)$ as above. In this way we
speak of the rows, columns and components of tableaux. If $\t$ is a tableau and
$1\le k\le n$ set $\comp_\t(k)=l$ if $k$ appears in the $l$th component of~$\t$.

A $\blam$-tableau is \textbf{standard} if its entries increase along rows and
columns in each component. Both of the tableaux above are standard. Let
$\Std(\blam)$ be the set of standard $\blam$-tableaux and let
$\Std(\Parts)=\bigcup_{\blam\in\Parts}\Std(\blam)$. Similarly set
$\SStd(\blam)=\set{(\s,\t)|\s,\t\in\Std(\blam)}$ and
$\SStd(\Parts)=\set{(\s,\t)|\s,\t\in\Std(\blam)\text{ for some
}\blam\in\Parts}$.

If $\t$ is a $\blam$-tableau set $\Shape(\t)=\blam$ and let $\t_{\downarrow m}$
be the subtableau of~$\t$ that contains the numbers $\{1,2,\dots,m\}$. If~$\t$ is
a standard $\blam$-tableau then $\Shape(\t_{\downarrow m})$ is a multipartition
for all $m\ge0$. We extend the dominance ordering to the set of all standard
tableaux by defining $\s\gedom\t$ if
$$\Shape(\s_{\downarrow m})\gedom\Shape(\t_{\downarrow m}),$$
for $1\le m\le n$. As before, we write $\s\gdom\t$ if $\s\gedom\t$ and
$\s\ne\t$. We extend the dominance ordering to $\SStd(\Parts)$ by declaring that
$(\s,\t)\gedom(\u,\v)$ if $\s\gedom\u$ and $\t\gedom\v$. Similarly,
$(\s,\t)\gdom(\u,\v)$ if $(\s,\t)\gedom(\u,\v)$ and $(\s,\t)\ne(\u,\v)$

It is easy to see that there are unique standard $\blam$-tableaux $\tlam$ and
$\tllam$ such that $\tlam\gedom\t\gedom\tllam$, for all $\t\in\Std(\blam)$.
The tableau $\tlam$ has the numbers $1,2,\dots,n$ entered in order from left to
right along the rows of $\t^{\lambda^{(1)}}$, and then
$\t^{\lambda^{(2)}},\dots,\t^{\lambda^{(\ell)}}$ and similarly, $\tllam$ is the
tableau with the numbers $1,\dots,n$ entered in order down the columns of
$\t^{\lambda^{(\ell)}},\dots,\t^{\lambda^{(2)}},\t^{\lambda^{(1)}}$. When
$\blam=(3,1^2|2,1|3,2)$ then the two $\blam$-tableaux displayed above are
$\tlam$ and $\tllam$.

Given a standard $\blam$-tableau $\t$ define $d(\t)\in\Sym_n$ to be the
permutation such that $\t=\tlam d(\t)$. Let $\le$ be the Bruhat order
on~$\Sym_n$ with the convention that $1\le w$ for all $w\in\Sym_n$. By a
well-known result of Ehresmann and James, if $\s,\t\in\Std(\blam)$ then
$\s\gedom\t$ if and only if $d(\s)\le d(\t)$; see, for example,
\cite[Theorem~3.8]{M:Ulect}.

Recall from \autoref{S:Quiver} that we have fixed a multicharge
$\kappa\in\Z^\ell$. The \textbf{residue} of the node $A=(l,r,c)$ is
$\res(A) = \kappa_l+c-r \pmod e$ (where we adopt the convention that
$i\equiv i\pmod 0$, for $i\in\Z$). Thus, $\res(A)\in I$. A node $A$ is an
\textbf{$i$-node} if $\res(A)=i$. If $\t$ is a $\bmu$-tableaux
and $1\le k\le n$ then the \textbf{residue} of $k$ in~$\t$ is $\res_\t(k)=\res(A)$,
where $A\in\bmu$ is the unique node such that $\t(A)=k$. The
\textbf{residue sequence} of~$\t$ is
$$\res(\t) = (\res_\t(1),\res_\t(2),\dots,\res_\t(n))\in I^n.$$
As an important special case we set $\bi^\mu=\res(\tmu)$, for $\bmu\in\Parts$.

Refine the dominance ordering on the set of standard tableaux by defining
$\s\Gedom\t$ if $\s\gedom\t$ and $\res(\s)=\res(\t)$. Similarly, we write
$(\s,\t)\Gedom(\u,\v)$ if $(\s,\t)\gedom(\u,\v)$, $\res(\s)=\res(\u)$ and
$\res(\t)=\res(\v)$ and $(\s,\t)\Gdom(\u,\v)$ now has the obvious meaning.

Following Brundan, Kleshchev and Wang~\cite[Definition.~3.5]{BKW:GradedSpecht}
we now define the degree of a standard tableau. Suppose that
$\bmu\in\Parts$. A node $A$ is an \textbf{addable node} of~$\bmu$ if
$A\notin\bmu$ and $\bmu\cup\{A\}$ is (the diagram of) a multipartition of~$n+1$.
Similarly, a node $B$ is a \textbf{removable node} of~$\bmu$ if $B\in\bmu$ and
$\bmu\setminus\{B\}$ is a multipartition of~$n-1$. Suppose that $A$ is
an $i$-node and define integers
$$
  d_A(\bmu)=\#\Setbox{addable $i$-nodes of $\bmu$\\ strictly below $A$}
                -\#\Setbox[36]{removable $i$-nodes of $\bmu$\\ strictly below $A$}.
$$
If $\t$ is a standard $\bmu$-tableau define its \textbf{degree} inductively by setting
$\deg_e(\t)=0$, if~$n=0$, and if $n>0$ then
\begin{equation}\label{E:TableauDegree}
   \deg_e(\t)=\deg_e(\t\rest{(n-1)})+d_A(\bmu),
\end{equation}
where $A=\t^{-1}(n)$. When $e$ is understood we write $\deg(\t)$.

The following result shows that the degrees of the standard tableau
are almost completely determined by the Cartan matrix $(c_{ij})$ of
$\Gamma_e$.

\begin{Lemma}[\protect{Brundan, Kleshchev and Wang~\cite[Proposition~3.13]{BKW:GradedSpecht}}]
  \label{L:TableauDegreeRecurrence}
  Suppose that~$\s$ and~$\t$ are standard tableaux such that $\s\gdom\t=\s(r,r+1)$,
  where $1\le r<n$ and $\bi\in I^n$. Let $\bi=\res(\s)$. Then
  $\deg_e(\s)=\deg_e(\t)+c_{i_ri_{r+1}}$.
\end{Lemma}

\subsection{The Murphy basis and cyclotomic Specht modules}\label{S:Murphy}
The cyclotomic Hecke algebra $\H$ is a cellular algebra with several different
cellular bases. This section introduces one of these bases, the Murphy basis,
and uses it to define the Specht modules and simple modules of~$\H$.

Fix a multipartition $\blam\in\Parts$. Following
\cite[Definition~3.14]{DJM:cyc} and~\cite[\Sect6]{AMR}, if $\s,\t\in\Std(\blam)$
define $m_{\s\t}=T_{d(\s)^{-1}}m_\blam T_{d(\t)}$, where
$m_\blam= u_\blam x_\blam$ where
$$u_\blam=\prod_{1\le l<\ell}\prod_{r=1}^{|\lambda^{(1)}|+\dots+|\lambda^{(l)}|}
                  \xi^{-\kappa_{l+1}}(L_r-[\kappa_{l+1}])
\quad\text{and}\quad x_\blam=\sum_{w\in\Sym_\blam}T_w.$$

Let $\ast$ be the unique anti-isomorphism of $\H$ that fixes each of the
generators $T_1,\dots,T_{n-1},L_1,\dots,L_n$ of
\autoref{D:HeckeAlgebras}.

\begin{Theorem}[\protect{\!\cite[Theorem~3.26]{DJM:cyc} and \cite[Theorem~6.3]{AMR}}]
  \label{T:MurphyBasis}
  The cyclotomic Hecke algebra $\H$ is free as an $\O$-module with cellular basis
  $$\set{m_{\s\t}|\s, \t\in\Std(\blam)\text{ for }\blam\in\Parts}$$
  with respect to the weight poset $(\Parts,\gedom)$ and automorphism~$\ast$.
\end{Theorem}

\begin{proof}
  This theorem can be proved uniformly in all cases by modifying the argument of
  \cite[Theorem~3.26]{DJM:cyc}, however, for future reference we explain how to
  deduce this result from the literature for the degenerate and non-degenerate
  algebras.

  First suppose that  $\xi=1$. Then the element $m_\blam$, for $\blam\in\Parts$,
  coincides exactly with the corresponding elements defined for the
  non-degenerate cyclotomic Hecke algebras in \cite[\Sect6]{AMR}. It follows
  that $\set{m_{\s\t}|(\s,\t)\in\Parts}$ is the Murphy basis of the degenerate
  cyclotomic Hecke algebra $\H$ defined in \cite[\Sect6]{AMR} and that the theorem
  is just a restatement of \cite[Theorem~6.3]{AMR} when $\xi=1$.

  Now suppose that $\xi\ne1$ and, as in \autoref{R:NewPresentation}, let
  $L_r'=(\xi-1)L_r+1$ be the `non-degenerate' Jucys-Murphy elements for $\H$, for
  $1\le r\le n$. An application of the definitions shows that if $\kappa\in\Z$
  then
   $$\xi^{-\kappa}(L_r-[\kappa])=\frac{\xi^{-\kappa}}{\xi-1}(L_r'-\xi^{\kappa}).$$
  Therefore, $u_\blam$ is a scalar multiple of the element $u_\blam^+$
  given by \cite[Definition~3.1,3.5]{DJM:cyc}. Consequently, if
  $(\s,\t)\in\SStd(\Parts)$ then $m_{\s\t}$ is a scalar multiple of the
  corresponding Murphy basis element from \cite[Definition~3.14]{DJM:cyc}.
  Hence, the theorem is an immediate consequence of \cite[Theorem~3.26]{DJM:cyc}
  in the non-degenerate case.
\end{proof}

Suppose that $\blam\in\Parts$. The (cyclotomic) \textbf{Specht module}
$\UnS^\blam$ is the cell module associated to $\blam$ using the (ungraded) cellular
basis $\set{m_{\s\t}|(\s,\t)\in\SStd(\Parts)}$. We underline $\UnS^\blam$ to emphasize
that $\UnS^\blam$ is not graded.  When $\O$ is a field let
$\UnD^\blam=\UnS^\blam/\rad\UnS^\blam$ and set
$\Klesh=\set{\blam\in\Parts|\UnD^\blam\ne0}$. Ariki~\cite{Ariki:class} has given a
combinatorial description of the set~$\Klesh$. By the theory of cellular
algebras~\cite{GL}, $\set{\UnD^\bmu|\bmu\in\Klesh}$
is a complete set of pairwise non-isomorphic irreducible $\H$-modules.

The following well-known fact is fundamental to all of the results in this paper.

\begin{Lemma}\label{L:JucysMurphyAction}
  Suppose that $1\le r\leq n$ and that $\s,\t\in\Std(\blam)$, for
  $\blam\in\Parts$. Then
  $$m_{\s\t}L_r \equiv [c_r(\t)]m_{\s\t}
  +\sum_{\substack{\v\gdom\t\\\v\in\Std(\blam)}}r_\v m_{\s\v}\pmod\Hlam,$$
  for some $r_\v\in\O$.
\end{Lemma}

\begin{proof}If $\xi=1$ then this is a restatement of \cite[Lemma~6.6]{AMR}. If
  $\xi\ne1$ then
  $$m_{\s\t}L_r'=\xi^{c_r(\t)}m_{\s\t} +\sum_{\v\gdom\t}r'_\v m_{\s\t}
  					\qquad\pmod\Hlam,$$
  for some $r'_v\in\O$, by \cite[Proposition~3.7]{JM:cyc-Schaper} (and the
  notational translations given in the proof
  of \autoref{T:MurphyBasis}). As $L_r=(L'_r-1)/(\xi-1)$ the result follows.
\end{proof}

\subsection{Cyclotomic quiver Hecke algebras}
Brundan and Kleshchev~\cite{BK:GradedKL} have given a very different
presentation of~$\H$. This presentation is more difficult to work with but it
has the advantage of showing that~$\H$ is a $\Z$-graded algebra.

\begin{Definition}[Brundan-Kleshchev~\cite{BK:GradedKL}]\label{D:QuiverRelations}
  Suppose that $n\ge0$ and $e\in\{0,2,3,4,\dots\}$. The
  \textbf{cyclotomic quiver Hecke algebra}, or \textbf{cyclotomic
  Khovanov-Lauda--Rouquier algebra}, of weight $\Lambda$ and type $\Gamma_e$
  is the unital associative $\O$-algebra $\R=\R(\O)$ with generators
  $$\{\psi_1,\dots,\psi_{n-1}\} \cup
  \{ y_1,\dots,y_n \} \cup \set{e(\bi)|\bi\in I^n}$$
  and relations
  \begin{align*}
    y_1^{(\Lambda,\alpha_{i_1})}e(\bi)&=0,
    & e(\bi) e(\bj) &= \delta_{\bi\bj} e(\bi),
    &{\textstyle\sum_{\bi \in I^n}} e(\bi)&= 1,\\
    y_r e(\bi) &= e(\bi) y_r,
    &\psi_r e(\bi)&= e(s_r{\cdot}\bi) \psi_r,
    &y_r y_s &= y_s y_r,\\[-8mm]
  \end{align*}
  \begin{align}
    \psi_r y_{r+1} e(\bi)&=(y_r\psi_r+\delta_{i_ri_{r+1}})e(\bi),&
    y_{r+1}\psi_re(\bi)&=(\psi_r y_r+\delta_{i_ri_{r+1}})e(\bi),\label{E:ypsi}\\
    \psi_r y_s  &= y_s \psi_r,&&\text{if }s \neq r,r+1,\label{E:ypsiCommute}\\
    \psi_r \psi_s &= \psi_s \psi_r,&&\text{if }|r-s|>1,\notag
\end{align}
\vskip-2em
\begin{align*}
  \psi_r^2e(\bi) &= \begin{cases}
       0,&\text{if }i_r = i_{r+1},\\
       (y_{r}-y_{r+1})e(\bi),&\text{if  }i_r\rightarrow i_{r+1},\\
       (y_{r+1} - y_{r})e(\bi),&\text{if }i_r\leftarrow i_{r+1},\\
       (y_{r+1} - y_{r})(y_{r}-y_{r+1}) e(\bi),
                &\text{if }i_r\rightleftarrows i_{r+1}\\
      e(\bi),&\text{otherwise},\\
\end{cases}\\
\psi_{r}\psi_{r+1} \psi_{r} e(\bi) &= \begin{cases}
    (\psi_{r+1} \psi_{r} \psi_{r+1}-1)e(\bi),\hspace*{7mm}
             &\text{if }i_{r+2}=i_r\rightarrow i_{r+1} ,\\
  (\psi_{r+1} \psi_{r} \psi_{r+1}+1)e(\bi), &\text{if }i_{r+2}=i_r\leftarrow i_{r+1},\\
  \rlap{$\big(\psi_{r+1} \psi_{r} \psi_{r+1} +y_r -2y_{r+1}+y_{r+2}\big)e(\bi)$,}\\
           &\text{if }i_{r+2}=i_r \rightleftarrows i_{r+1},\\
  \psi_{r+1} \psi_{r} \psi_{r+1} e(\bi),&\text{otherwise,}
\end{cases}
\end{align*}
for $\bi,\bj\in I^n$ and all admissible $r$ and $s$. Moreover, $\R$ is naturally
$\Z$-graded with degree function determined by
$$\deg e(\bi)=0,\qquad \deg y_r=2\qquad\text{and}\qquad \deg
  \psi_s e(\bi)=-c_{i_s,i_{s+1}},$$
for $1\le r\le n$, $1\le s<n$ and $\bi\in I^n$.
\end{Definition}

\begin{Remark}
  The presentation of $\R$ given in \autoref{D:QuiverRelations} differs by a
  choice of signs with the definition given in \cite[Theorem~1.1]{BK:GradedKL}.
  The presentation of~$\R$ given above agrees with that used in
  \cite{KMR:UniversalSpecht} as the orientation of the quiver is reversed
  in~\cite{KMR:UniversalSpecht}.
\end{Remark}

The connection between the cyclotomic quiver Hecke algebras of type $\Gamma_e$
and the cyclotomic Hecke algebras of type~$G(\ell,1,n)$ is given by the
following remarkable result of Brundan and Kleshchev.

\begin{Theorem}[\protect{Brundan-Kleshchev's isomorphism theorem~\cite[Theorem~1.1]{BK:GradedKL}}]
\label{T:BKiso}
Suppose that $\O=K$ is a field, $\xi\in K$ as above, and that
$\Lambda=\Lambda(\charge)$. Then there is an isomorphism of algebras
$\UnR\cong\H$.
\end{Theorem}

Rouquier \cite[Corollary~3.20]{Rouq:2KM} has, independently, given
a quick proof of \autoref{T:BKiso}.

We prove a stronger version of \autoref{T:BKiso} in \autoref{T:IntegralKLR}
below. For now we note the following simple corollary of \autoref{T:BKiso}.
Recall that a choice of multicharge $\charge$ determines a dominant weight
$\Lambda_e(\charge)$.

\begin{Corollary}\label{C:Largee}
  Suppose that $n\ge0$, $\charge=(\kappa_1,\dots,\kappa_\ell)\in\Z^\ell$ and that
  $$e>\max\set{n+\kappa_k-\kappa_l|1\le k,l\le\ell}.$$
  Fix invertible scalars $\xi_0\in K$ and $\xi_e\in K$ such that $\xi_0$ is not
  a root of unity and $\xi_e$ is a primitive $e$th root of unity. Then the
  cyclotomic Hecke algebras $\HH^{\Lambda_0(\kappa)}_{K,\xi_0}$ and
  $\HH^{\Lambda_e(\charge)}_{\K,\xi_e}$ are isomorphic $\Z$-graded $K$-algebras.
\end{Corollary}

\begin{proof}Let $\R(0)\cong\HH_n(K,\xi_0,\charge)$ and
  $\R(e)\cong\HH_n(K,\xi_e,\charge)$
  be the corresponding cyclotomic quiver Hecke algebras as in
  \autoref{T:BKiso}. By \cite[Lemma~4.1]{HuMathas:GradedCellular}, $e(\bi)\ne0$
  if and only if $\bi=\res(\t)$, for some standard tableau $\t\in\Std(\Parts)$.
  The definition of $e$ ensures that if $\bi=\bi^{\t}$ then $i_r=i_{r+1}$ or
  $i_r=i_{r+1}\pm 1$ if and only if $i_r\equiv i_{r+1}\pmod{e}$ or $i_r\equiv
  i_{r+1}\pm 1\pmod{e}$. Therefore, $\R(0)\cong\R(e)$ arguing directly from the
  presentations of the cyclotomic quiver Hecke algebras given in
  \autoref{D:QuiverRelations}. Hence, the result follows by
  \autoref{T:BKiso}.
\end{proof}

Therefore, without loss of generality, we may assume that $e>0$. In the appendix
we show how to modify the results and definitions in this paper to cover the
case when $e=0$ directly.


Under the assumptions of the Corollary we note that the algebras
$\HH^\Lambda_{K,\xi_0}$ and $\HH^\Lambda_{K,\xi_e}$ are Morita equivalent by the
main result of \cite{DM:Morita}. That these algebras are actually isomorphic is
another miracle provided by Brundan and Kleshchev's isomorphism theorem.

%
%

\section{Seminormal forms for Hecke algebras}\label{Chap:Seminormal}

In this chapter we develop the theory of seminormal forms in a
slightly more general context than appears in the literature. In
particular, in this paper a seminormal basis will be a basis
for~$\H$ rather than a basis of a Specht module of~$\H$. We also
treat all of the variations of the seminormal bases simultaneously
as this will give us the flexibility to  change seminormal forms
when we use them in the next chapter to study the connections
between $\H$ and the cyclotomic quiver Hecke algebra~$\R$.

\subsection{Content functions and the Gelfand-Zetlin algebra}\label{S:Contents}
Underpinning Brundan and Kleshchev's isomorphism theorem
(\autoref{T:BKiso}) is the decomposition of any $\H$-module into a direct sum of
generalised eigenspaces for the Jucys-Murphy elements $L_1,\dots,L_n$. This
section studies the action of the Jucys-Murphy elements on~$\H$. The results in
this section are well-known, at least to experts, but they are needed in the sequel.

The \textbf{content} of the node $\gamma=(l,r,c)$ is the integer
$$\cont_\gamma=\kappa_l-r+c$$
If $\t\in\Std(\blam)$ is a standard $\blam$-tableau and $1\le k\le n$ then the
\textbf{content} of $k$ in~$\t$ is $\cont_k(\t)=\cont_\gamma$, where
$\t(\gamma)=k$ for $\gamma\in\diag\blam$.

\begin{Definition}\label{D:Separation}
Let $\O$ be a commutative integral domain and suppose that $t\in\O^\times$ is an
invertible element of~$\O$. The pair $(\O,t)$ \textbf{separates} $\Std(\Parts)$ if
$$[n]^!_t\prod_{1\le l<m\le\ell}\prod_{-n<d<n}[\kappa_l-\kappa_m+d]_t\in {\O}^{\times}.$$
\end{Definition}

Fix a multicharge $\charge\in\Z^\ell$ and let $\HO$ be the Hecke algebra defined
over~$\O$ with parameter~$t$. In spite of our notation, note that $\HO$ depends
only on $\charge$ and not directly on $\Lambda=\Lambda_e(\charge)$. Let
$\K$ be a field that contains the field of fractions of~$\O$. Then $\HK=\HO\otimes_\O\K$.

Throughout this chapter we are going to work with the Hecke algebras $\HO$ and
$\HK=\HO\otimes_\O\K$, however, we have in mind the situation of
\autoref{T:BKiso}. By assumption $e>0$, so we can replace the
multicharge~$\charge$ with $(\kappa_1+a_1e,\kappa_2+a_2e,\dots,\kappa_\ell+a_\ell e)$,
for any integers $a_1,\dots,a_\ell\in\Z$, without changing the dominant weight
$\Lambda=\Lambda_e(\charge)$.  In view of \autoref{D:Separation} we therefore
assume that
\begin{equation}\label{E:ChargeSeparation}
    \kappa_l-\kappa_{l+1}\ge n,\qquad\text{ for }1\le l<\ell.
\end{equation}
Until further notice, we fix a multicharge $\charge\in\Z^\ell$ satisfying
\autoref{E:ChargeSeparation} and consider the algebra $\HO$ with parameter~$t$.

Although we do not need this, we remark that it follows from \cite{Ariki:ss} and
\cite[Theorem~6.11]{AMR} that $\H(\K,t)$ is semisimple if and only if $(\K,t)$
separates $\Std(\Parts)$. Our main use of the separation condition is the
following fundamental fact that is easily proved by induction on~$n$; see, for
example, \cite[Lemma~3.12]{JM:cyc-Schaper}.

\begin{Lemma}\label{L:Separation}
  Suppose that $\O$ is an integral domain and $t\in\O^\times$ is invertible.
  Then the following are equivalent:
  \begin{enumerate}
    \item $(\O,t)$ separates $\Std(\Parts)$,
    \item If $\s,\t\in\Std(\Parts)$ then $\s=\t$ if and only if
    $[c_r(\s)]=[c_r(\t)]$, for $1\le r\le n$.
\end{enumerate}
\end{Lemma}

Following~\cite{OkounkovVershik}, define the \textbf{Gelfand-Zetlin subalgebra}
of~$\H$ to be the algebra $\L(\O)=\<L_1,\dots,L_n\>$. The aim of this section is
to understand the semisimple representation theory of~$\L=\L(\O)$. It follows
from \autoref{D:HeckeAlgebras} that~$\L$ is a commutative subalgebra of~$\H$.

If $\O$ is an integral domain then it follows from \autoref{L:JucysMurphyAction}
that, as an $(\L,\L)$-bimodule, $\H(\O)$ has a composition series with
composition factors that are $\O$-free of rank~$1$ upon which $L_r$ acts as
multiplication by $[c_r(\s)]$ from the left and as multiplication by $[c_r(\t)]$
from the right. Obtaining a better description of~$\L$, and of~$\H$ as an
$(\L,\L)$-bimodule, in the non-semisimple case is likely to be important. For
example, the dimension of $\L$ over a field is not known in general.

\begin{Proposition}[\protect{cf. \cite[Proposition~3.17]{AK}}]\label{P:LDecomposition}
  Suppose that $(\K,t)$ separates $\Std(\Parts)$, where $\K$ is a field and
  $0\ne t\in\K$. Then $\H(\K)$ is a semisimple $\(\L,\L)$-bimodule with
  decomposition
  $$
  \H(\K)=\bigoplus_{\substack{\blam\in\Parts\\ \s,\t\in\Std(\blam)}}H_{\s\t},
  $$
  where $H_{\s\t}=\set{h\in\H|L_r h=[c_r(\s)]h\text{ and }hL_r=[c_r(\t)]h,
  \text{ for } 1\le r\le n}$ is one dimensional.
\end{Proposition}

\begin{proof} By \autoref{L:JucysMurphyAction}, the Jucys-Murphy elements
  $L_1,\dots,L_n$ are a family of JM-elements for $\H$ in the sense of
  \cite[Definition~2.4]{M:seminormal}. Therefore, the result is a special case
  of \cite[Theorem~3.7]{M:seminormal}.
\end{proof}

Key to the proof of the results in \cite{M:seminormal} are the following
elements that have their origins in the work of Murphy~\cite{M:Nak}. For
$\t\in\Std(\Parts)$ define
\begin{equation}\label{D:MurphyIdempotents}
    F_\t=\prod_{k=1}^n\prod_{\substack{c\in\Cont\\ [\cont_k(\t)]\ne[c]}}
    \frac{L_k-[c]}{[\cont_k(\t)]-[c]}
\end{equation}
where $\Cont=\set{c_r(\t)|1\le r\le n\text{ and }\t\in\Std(\Parts)}$ is the set of
the possible contents that can appear in a standard tableau of size~$n$. By
definition, $F_\t\in\L(\K)$ and it follows directly from
\autoref{P:LDecomposition} that if $h_{\u\v}\in H_{\u\v}$ then
\begin{equation}\label{E:FInvariance}
F_\s h_{\u\v} F_\t = \delta_{\s\u}\delta_{\v\t}h_{\s\t},
\end{equation}
for all $(\s,\t),(\u,\v)\in\SStd(\Parts)$. Therefore,
$H_{\s\t}=F_\s\H F_\t$.

By \autoref{P:LDecomposition} we can write $1=\sum_{\s,\t}e_{\s\t}$ for unique
$e_{\s\t}\in H_{\s\t}$. Since $F_\t=F_\t^{\bigast}$, the last displayed equation implies that
$F_\t=e_{\t\t}\in H_{\t\t}$ is an idempotent. Consequently,
$$\L(\K)=\bigoplus_{\t\in\Std(\Parts)} H_{\t\t}=\bigoplus_{\t\in\Std(\Parts)}\K F_\t.$$
In particular, $F_\t$ is a primitive idempotent in $\L(\K)$. If follows that
$\L(\K)$ is a split semisimple algebra of dimension $\#\Std(\Parts)$.



\subsection{Seminormal forms}\label{S:SeminormalForm} Seminormal bases for $\H$
are well-known in the literature, having their origins in the work of
Young~\cite{QSAI}.  Many examples of ``seminormal bases'' appear in the
literature. In this section we classify the seminormal bases of~$\H$.  This
characterisation of seminormal forms appears to be new, even in the special case
of the symmetric groups, although some of the details will be familiar to
experts.

Throughout this section we assume that $\K$ is a field, $0\ne t\in\K$ and that
$(\K,t)$ separates $\Std(\Parts)$.  Recall the decomposition
$\H=\bigoplus_{(\s,\t)\in\SStd(\Parts)}H_{\s\t}$ from
\autoref{P:LDecomposition}.

Define an \textbf{anti-involution} on an algebra $A$ to be an algebra
anti-automorphism of~$A$ of order~$2$.

\begin{Definition}\label{D:SeminormalBasis}
  Suppose that $(\K,t)$ separates $\Std(\Parts)$ and let $\iota$ be an anti-involution on $\HK$.
  An \textbf{$\iota$-seminormal basis}
  of $\H(\K)$ is a basis of the form
  $\set{f_{\s\t}|f_{\s\t}=\iota(f_{\t\s})\in H_{\s\t}\text{ for
  }(\s,\t)\in\SStd(\Parts)}$.
\end{Definition}

Recall that $\bigast$ is the unique anti-involution of $\HK$ that fixes
each of the generators $T_1,\dots,T_{n-1},L_1,\dots,L_n$. Then $m_{\s\t}^{\bigast}=m_{\t\s}$, for all
$(\s,\t)\in\SStd(\Parts)$. The assumption that $f_{\s\t}^{\bigast}=f_{\t\s}$ is not
essential for what follows but it is natural because we want to work within the
framework of cellular algebras.

In order to describe the action of $\H$ on its seminormal bases, if
$\t\in\Std(\Parts)$ then define the integers
\begin{equation}\label{E:rho}
  \rho_r(\t)=c_r(\t)-c_{r+1}(\t),\qquad\text{for }1\le r<n.
\end{equation}
Then $\rho_r(\t)$ is the `axial distance' between $r$ and $r+1$ in the
tableau~$\t$.

\begin{Definition}\label{D:SNCS}
A \textbf{$\bigast$-seminormal coefficient system} for $\HKK$ is a set of scalars
$\balpha=\set{\alpha_r(\s)|1\le r<n \text{ and } \s\in\Std(\Parts)}$
in~$\K$ such that if $\t\in\Std(\Parts)$ and $1\le r<n$ then:
\begin{enumerate}
\item
$\alpha_r(\t)\alpha_{r+1}(\t s_r)\alpha_r(\t s_rs_{r+1})
      =\alpha_{r+1}(\t)\alpha_r(\t s_{r+1})\alpha_{r+1}(\t s_{r+1}s_r)$
       if $r<n-1$,
\item
$\alpha_r(\t)\alpha_k(\t s_r)=\alpha_k(\t)\alpha_r(\t s_k)$
if $1\le k<n$ and $|r-k|>1$,
\item if $\v=\t(r,r+1)$ then $\alpha_r(\v)=0$ if $\v\notin\Std(\Parts)$ and
otherwise
\begin{equation*}
   \alpha_r(\t)\alpha_r(\v)=
     \dfrac{[1+\rho_r(\t)][1+\rho_r(\v)]}{[\rho_r(\t)][\rho_r(\v)]}.
\end{equation*}
\end{enumerate}
\end{Definition}

We will see that conditions~(a) and (b) correspond to the braid relations
satisfied by $T_1,\dots,T_{n-1}$ and that~(c) corresponds to the
quadratic relations. Quite surprisingly, as the proof of
\autoref{T:GramDetFactorization} below shows, \autoref{D:SNCS}(c) also
encodes the KLR grading on~$\H$.

Usually, we omit the $\bigast$ and simply call $\balpha$ a seminormal coefficient
system.

\begin{Example} A nice `rational' seminormal coefficient system is given by
$$  
\alpha_r(\t)=\begin{cases}
  \frac{[1+\rho_r(\t)]}{[\rho_r(\t)]}, & \text{if $\t(r,r+1)$ is standard},\\
  0,&\text{otherwise},
\end{cases}
$$ 
for $\t\in\Std(\Parts)$ and $1\le r<n$.
\end{Example}

\begin{Example}\label{Ex:MurphySNC} By \autoref{P:MurphySeminormal} below, the following
  seminormal coefficient system is associated with the Murphy basis of $\H$: if
  $\t\in\Std(\Parts)$ set $\v=\t(r,r+1)$ and define
$$ \alpha_r(\t)=\begin{cases}
         1& \text{if $\v$ is standard and $\t\gdom\v$},\\
         \frac{[1+\rho_r(\t)][1+\rho_r(\v)]}{[\rho_r(\t)][\rho_r(\v)]},
          & \text{if $\v$ is standard and $\v\gdom\t$},\\
  0,&\text{otherwise},
\end{cases}
$$
for $1\le r<n$.
\end{Example}

Another seminormal coefficient system, which is particularly well adapted to
Brundan and Kleshchev's Graded Isomorphism \autoref{T:BKiso}, is given in
\autoref{S:PsiBasis}.

\begin{Lemma}\label{L:SNCExistence}
  Suppose that $(\K,t)$ separates $\Std(\Parts)$ and that $\{f_{\s\t}\}$ is a
  seminormal basis of~$\H$. Then there exists a unique seminormal coefficient
  system $\balpha$ such that if $1\le r<n$ and $(\s,\t)\in\SStd(\Parts)$ then
  $$f_{\s\t}T_r = \alpha_r(\t)f_{\s\v}-\frac{1}{[\rho_r(\t)]}f_{\s\t},$$
  where $\v=\t(r,r+1)$ and $f_{\s\t}=0$ if $(\s,\t)\notin\SStd(\Parts)$.
\end{Lemma}

\begin{proof} The uniqueness statement is automatic, since $\{f_{\s\t}\}$ is a
  basis of~$\HK$, so we need to prove that such a seminormal coefficient system
  $\balpha$ exists.

Fix $(\s,t)\in\SStd(\Parts)$ and $1\le r<n$  and write
$$ f_{\s\t}T_r = \sum_{(\u,\v)\in\SStd(\Parts)}a_{\u\v}f_{\u\v},$$
for some $a_{\u\v}\in\K$. Multiplying on the left by $F_\s$ it follows that
$a_{\u\v}\ne0$ only if~$\u=\s$. If $k\ne r,r+1$ then $L_k$ commutes with $T_r$
so it follows $a_{\s\v}\ne0$ only if $[c_k(\v)]=[c_k(\t)]$, for $k\ne r,r+1$.
Using \autoref{D:Separation}, and arguing as in \autoref{L:Separation}, this
implies that $a_{\s\v}\ne0$ only if $\v\in\{\t,\t(r,r+1)\}$. Therefore, we can
write
$$f_{\s\t}T_r=\alpha_r(\t)f_{\s\v}+\alpha'_r(\t)f_{\s\t},$$
for some $\alpha_r(\t),\alpha'_r(\t)\in\K$, where $\v=\t(r,r+1)$. (Here, and
below, we adopt the convention that $f_{\s\v}=0$ if either of $\s$ or $\v$ is
not standard.) By \autoref{D:HeckeAlgebras}, $T_rL_r=L_{r+1}(T_r-t+1)-1$, so
multiplying both sides of the last displayed equation on the right by~$L_r$ and
comparing the coefficient of $f_{\s\t}$ on both sides shows that
$$[c_{r+1}(\t)]\(\alpha'_r(\t)-t+1\) -1 = \alpha'_r(\t)[c_r(\t)].$$
Hence, $\alpha'_r(\t)=-1/[\rho_r(\t)]$ as claimed. If~$\v$ is not standard then
we set $\alpha_r(\t)=0$. If~$\v$ is standard then comparing the
coefficient of $f_{\s\t}$ on both sides of
$$\Big(\alpha_r(\t)f_{\s\v}-\frac{1}{[\rho_r(\t)]}f_{\s\t}\Big)T_r
     =f_{\s\t}T_r^2=f_{\s\t}\((t-1)T_r+t\)$$
shows that
$\alpha_r(\t)\alpha_r(\v)=\tfrac{[1+\rho_r(\t)][1+\rho_r(\v)]}{[\rho_r(\t)][\rho_r(\v)]}$
in accordance with \autoref{D:SNCS}(c).

If $1\leq r<s-1<n-1$ and $(\s,\t)\in\SStd(\Parts)$ then
$(f_{\s\t}T_r)T_s=f_{\s\t}(T_rT_s)=f_{\s\t}(T_sT_r)=(f_{\s\t}T_s)T_r$. By a direct calculation, we can deduce that
$\alpha_r(\t)\alpha_s(\t(r,r+1))=\alpha_s(\t)\alpha_r(\t(s,s+1))$.

Finally, it remains to show that \autoref{D:SNCS}(a) holds. If $1\le r<n$
then $T_rT_{r+1}T_r=T_{r+1}T_rT_{r+1}$ by \autoref{D:HeckeAlgebras}. On the
other hand, if we set $\t_1=\t(r,r+1)$, $\t_2=\t(r+1,r+2)$, $\t_{12}=\t_1(r+1,r+2)$,
$\t_{21}=\t_2(r,r+1)$ and $\t_{121}=\t_{212}=\t(r,r+2)$ then direct calculation
shows that $0=f_{\s\t}(T_rT_{r+1}T_r-T_{r+1}T_rT_{r+1})$ is equal to
\begin{align*}
  -\Big(\frac1{[\rho_r(\t)]^2[\rho_{r+1}(\t)]}
           -\frac1{[\rho_r(\t)][\rho_{r+1}(\t)]^2}
           +\frac{\alpha_r(\t)\alpha_r(\t_1)}{[\rho_{r+1}(\t_1)]}
           -\frac{\alpha_{r+1}(\t)\alpha_{r+1}(\t_2)}{[\rho_r(\t_2)]}\Big)f_{\s\t}\\
  +\alpha_r(\t)\Big(\frac1{[\rho_r(\t_1)][\rho_{r+1}(\t_1)]}
           +\frac1{[\rho_r(\t)][\rho_{r+1}(\t)]}
           -\frac1{[\rho_{r+1}(\t)][\rho_{r+1}(\t_1)]}\Big)f_{\s\t_1}\\
  -\alpha_{r+1}(\t)\Big(\frac1{[\rho_r(\t_2)][\rho_{r+1}(\t_2)]}
           +\frac1{[\rho_r(\t)][\rho_{r+1}(\t)]}
           -\frac1{[\rho_r(\t)][\rho_r(\t_2)]}\Big)f_{\s\t_2}\\
  -\alpha_r(\t)\alpha_{r+1}(\t_1)\Big(\frac1{[\rho_r(\t_{12})]}
           -\frac1{[\rho_{r+1}(\t)]}\Big)f_{\s\t_{12}}\\
  +\alpha_{r+1}(\t)\alpha_r(\t_2)\Big(\frac1{[\rho_{r+1}(\t_{21})]}
           -\frac1{[\rho_r(\t)]}\Big)f_{\s\t_{21}}\\
  +\big(\alpha_r(\t)\alpha_{r+1}(\t_1)\alpha_r(\t_{12})
  -\alpha_{r+1}(\t)\alpha_r(\t_2)\alpha_{r+1}(\t_{21})\big)f_{\s\t_{121}}.
\end{align*}
By our conventions, if any tableau $\t_?$ is not standard then $f_{\s\t_?}$
and the corresponding $\balpha$-coefficient are both zero.  As the coefficient of
$f_{\s\t_{121}}$ in the last displayed equation is zero it follows that \autoref{D:SNCS}(a) holds.
Consequently, $\balpha=\{\alpha_r(\t)\}$ is a seminormal coefficient system,
completing the proof. (It is not hard to see, using \autoref{D:SNCS} and
identities like $\rho_r(\t_1)=-\rho_r(\t)$ and $\rho_r(\t_{12})=\rho_{r+1}(\t)$,
that the remaining coefficients in the last displayed equation are automatically zero.)
\end{proof}

\autoref{L:SNCExistence} really says that acting from the right on a seminormal
basis determines a seminormal coefficient system. Similarly, the left action on a seminormal
basis determines a seminormal coefficient system. In general, the seminormal
coefficient systems attached to the left and right actions will be different,
however, because we are assuming that our seminormal bases are $\bigast$-invariant
these left and right coefficient systems coincide.  Thus, for
$(\s,\t)\in\SStd(\Parts)$ and $1\le r<n$ we also have $T_rf_{\s\t} =
\alpha_r(\s)f_{\u\t}-\frac{1}{[\rho_r(\s)]}f_{\s\t}$, where $\u=\s(r,r+1)$.

Exactly as eigenvectors are not uniquely determined by their eigenvalues,
seminormal bases are not uniquely determined by seminormal coefficient systems.
We now fully characterize seminormal bases --- and prove a converse to
\autoref{L:SNCExistence}.

Recall that a set of idempotents in an algebra is \textbf{complete} if they sum
to~$1$.

\begin{Theorem}[The Seminormal Basis Theorem]\label{T:Seminormal}
  Suppose that $(\K,t)$ separates $\Std(\Parts)$ and that $\balpha$ is a
  seminormal coefficient system for~$\HKK$. Then $\HKK$ has a $\bigast$-seminormal basis
  $\set{f_{\s\t}|(\s,\t)\in\SStd(\Parts)}$ such that if  $(\s,\t)\in\SStd(\Parts)$ then
  \begin{equation}\label{E:SeminormalForm}
    f_{\s\t}^{\bigast}=f_{\t\s},\quad f_{\s\t}L_k=[c_k(\t)]f_{\s\t}\quad\text{and}\quad
    f_{\s\t}T_r=\alpha_r(\t)f_{\s\v}-\frac{1}{[\rho_r(\t)]}f_{\s\t},
  \end{equation}
  where $\v=\t(r,r+1)$ and $f_{\s\v}=0$ if $\v$ is not standard. Moreover, there
  exist non-zero scalars $\gamma_\t\in\K$, for $\t\in\Std(\Parts)$, such that
  \begin{equation}\label{E:SeminormalGamma}
      F_\u f_{\s\t} F_\v = \delta_{\u\s}\delta_{\t\v}f_{\s\t},\quad
      f_{\s\t}f_{\u\v}=\delta_{\t\u}\gamma_\t f_{\s\v},\quad\text{and}\quad
      F_\t=\frac1{\gamma_\t}f_{\t\t}.
  \end{equation}
   Furthermore, $\set{F_\t|\t\in\Std(\Parts)}$ is a complete set of pairwise
   orthogonal primitive idempotents. In particular, every irreducible
   $\HKK$-module is isomorphic to~$F_\s\HKK$, for some $\s\in\Std(\Parts)$, and
   $F_\s\HKK\cong F_\u\HKK$ if and only if $\Shape(\s)=\Shape(\u)$.

   Finally, the basis
   $\set{f_{\s\t}|\s,\t\in\Std(\blam)\text{ for }\blam\in\Parts}$ is uniquely
   determined by the choice of seminormal coefficient system~$\balpha$ and the
   scalars $\set{\gamma_{\tlam}|\blam\in\Parts}\subseteq\K^\times$.
\end{Theorem}

\begin{proof}
  For each $\blam\in\Parts$ fix an arbitrary pair of tableaux and a non-zero
  element $f_{\s\t}\in H_{\s\t}$. Then $f_{\s\t}$ is a simultaneous eigenvector for all of
  the elements of~$\L$, where they act from the left and from the right.

  Now, suppose that $1\le r<n$ and that $\v=\t(r,r+1)$ is standard. Then
  $\alpha_r(\t)\ne0$ so we can set
  $f_{\s\v}=\frac1{\alpha_r(\t)}f_{\s\t}(T_r+\frac{1}{[\rho_r(\t)]})$.
  Equivalently, $f_{\s\t}T_r=\alpha_r(\t)f_{\s\v}-\frac1{[\rho_r(\t)]}f_{\s\t}$.
  Then using the relations in~$\HKK$ and the defining properties of the
  seminormal coefficient system~$\balpha$, it is straightforward to check that
  $f_{\s\v}L_k=[c_k(\v)]f_{\s\v}$, so that $f_{\s\v}\in H_{\s\v}$. Moreover,
  $f_{\s\v}\ne0$ since
  $f_{\s\t}=\frac1{\alpha_r(\v)}f_{\s\v}(T_r+\frac1{[\rho_r(\v)]})$.

  More generally, it is easy to see that if $\v$ is any $\blam$-tableau then
  there is a sequence of standard tableaux $\v_0=\s, \v_1,\dots,\v_z=\v$ such
  that $\v_{i+1}=\v_i(r_i,r_i+1)$, for some integers $1\le r_i<n$. Therefore,
  continuing in this way it follows that given two tableaux
  $\u,\v\in\Std(\blam)$ we can define non-zero elements $f_{\u\v}\in H_{\u\v}$
  that satisfy \autoref{E:SeminormalForm}. It follows that, once~$f_{\s\t}$ is
  fixed, there is at most one choice of elements
  $\set{f_{\u\v}|\u,\v\in\Std(\blam)}$, such that \autoref{E:SeminormalForm}
  holds.

  To complete the proof that the seminormal coefficient system determines a
  seminormal basis we need to check that the elements $f_{\u\v}$ from the
  last paragraph are well-defined. That is, we need to show that
  $f_{\u\v}$ is independent of the choice of the sequences of simple
  transpositions that link~$\u$ and~$\v$ to~$\s$ and~$\t$, respectively.
  Equivalently, we need to prove that the action of~$\HKK$ given by
  \autoref{E:SeminormalForm} respects the relations of~$\HKK$. Using
  \autoref{E:SeminormalForm}, all of the relations in \autoref{D:HeckeAlgebras}
  are easy to check except for the braid relations of length three that hold by
  virtue of the argument of \autoref{L:SNCExistence}. Hence, by choosing
  elements $f_{\s\t}\in H_{\s\t}$, for $(\s,\t)\in\SStd(\blam)$ and $\blam\in\Parts$,
  the seminormal coefficient system determines a unique seminormal basis.


  Using \autoref{E:FInvariance} it is straightforward to prove
  \autoref{E:SeminormalGamma} so we leave these details to the reader; cf.
  \cite[Theorem~3.16]{M:seminormal}. In particular, this shows that
  $F_\s=\frac1{\gamma_\s}f_{\s\s}$ is an idempotent. To show that $F_\s$ is
  primitive, suppose that~$a$ is a non-zero element of $F_\s\HKK$.
  By~\autoref{E:SeminormalForm}, $a=\sum_{\v\in\Std(\blam)}r_\v f_{\s\v}$, for
  some $r_\v\in\K$. Fix $\t\in\Std(\blam)$ such that $r_\t\ne0$. Then
  $f_{\s\t}=1/r_\t a F_\t\in F_\s\HKK$. Using~\autoref{E:SeminormalForm} we
  deduce that $F_\s\HKK$ has basis $\set{f_{\s\v}|\v\in\Std(\blam)}$.
  Consequently, $a\H=F_\s\HK$, showing that $F_\s\HKK$ is irreducible.
  Therefore, $F_\s$ is a primitive idempotent in~$\HK$.

  The last paragraph, together with \autoref{D:SNCS}(c), implies that if
  $\s,\u\in\Std(\blam)$ then $F_\s\H\cong F_\u\H$ where an isomorphism is given
  by $f_{\s\t}\mapsto f_{\u\t}$, for $\t\in\Std(\blam)$. Consequently, if $\s$
  and $\u$ are standard tableaux of different shape then
  $F_\s\H\not\cong F_\u\H$ because the multiplicity of $S^\blam\cong F_\s\H(\K)$
  in $\H(\K)$ is $\#\Std(\blam)$ by the Wedderburn theorem.

  Finally, it remains to show that the basis $\{f_{\s\t}\}$ is uniquely
  determined by~$\balpha$ and the choice of the $\gamma$-coefficients
  $\set{\gamma_{\tlam}|\blam\in\Parts}$. If
  $\s,\t\in\Std(\blam)$ then we have shown that, once $f_{\s\t}$ is
  fixed, there is a unique seminormal basis
  $\set{f_{\u\v}|\u,\v\in\Std(\blam)}$ satisfying~\autoref{E:SeminormalForm}. In particular,
  taking $\s=\tlam=\t$ and fixing~$f_{\tlam\tlam}$ determines these basis
  elements. By \autoref{E:SeminormalGamma} the choice of~$f_{\tlam\tlam}$ also
  uniquely determines~$\gamma_{\tlam}$. Conversely, by setting
  $f_{\tlam\tlam}=\gamma_{\tlam}F_{\tlam}$ for any choice of non-zero scalars
  $\gamma_{\tlam}\in\K$, for $\blam\in\K$, the seminormal coefficient
  system~$\balpha$ determines a unique seminormal basis.
\end{proof}

The results that follow are independent of the choice of seminormal coefficient
system~$\balpha$, however, the choice of $\gamma$-coefficients will be important
--- and in what follows it will be useful to be able to vary both the seminormal
coefficient system~$\balpha$ and the $\gamma$-coefficients.

The proof of \autoref{T:Seminormal} implies that the choice of $\gamma_{\tlam}$
determines all of the scalars~$\gamma_\s$, for $\s\in\Std(\blam)$. In what
follows we need the following result that makes the relationship
between these coefficients more explicit.

\begin{Corollary}\label{C:GammaRecurrence}
   Suppose that $\t\in\Std(\Parts)$ and
   that $\v=\t(r,r+1)$ is standard, where $1\le r<n$. Then
   $\alpha_r(\v)\gamma_\t=\alpha_r(\t)\gamma_\v$.
\end{Corollary}

\begin{proof}
  Applying \autoref{E:SeminormalForm} and \eqref{E:SeminormalGamma} several times each,
  \begin{align*}
  \gamma_\v f_{\v\v}&=f_{\v\v}f_{\v\v}
     =\frac1{\alpha_r(\t)}f_{\v\t}\Big(T_r+\frac1{[\rho_r(\t)]}\Big)f_{\v\v}
     =\frac1{\alpha_r(\t)}f_{\v\t}T_rf_{\v\v}\\
     &=\frac1{\alpha_r(\t)}f_{\v\t}
    \Big(\alpha_r(\v)f_{\t\v}-\frac{1}{[\rho_r(\v)]}f_{\v\v}\Big)
     =\frac{\alpha_r(\v)}{\alpha_r(\t)}f_{\v\t}f_{\t\v}\\
    &=\frac{\alpha_r(\v)}{\alpha_r(\t)}\gamma_\t f_{\v\v}.
  \end{align*}
  Comparing coefficients, $\alpha_r(\t)\gamma_\v=\alpha_r(\v)\gamma_\t$ as
  required.
\end{proof}

\subsection{Seminormal bases and the Murphy basis}\label{S:GramDet}
In this section we compute the Gram determinant of the Specht modules of $\H$,
with respect to the Murphy basis, as a product of cyclotomic polynomials
when $\xi\ne1$ or as a product of primes when $\xi=1$. These
determinants are already explicitly
known~\cite{AMR,JM:det,JM:Schaper,JM:cyc-Schaper} but all existing
formulas describe them as products of rational functions, or of rational
numbers in the degenerate case.

By \autoref{T:MurphyBasis}, the Murphy basis $\{m_{\s\t}\}$ is a cellular basis
for $\H$ over an arbitrary ring. In this section we continue to work with the
generic Hecke algebra $\H=\HO$ with parameter~$t$ and multicharge $\charge$
satisfying \autoref{E:ChargeSeparation}.

As $(\K,\t)$ separates $\Std(\Parts$), for $\s,\t\in\Std(\blam)$ we can define
$$f_{\s\t}=F_\s m_{\s\t} F_\t.$$
By \autoref{L:JucysMurphyAction},
$f_{\s\t}\equiv m_{\s\t}+\sum r_{\u\v}m_{\u\v}\pmod\Hlam$, for some $r_{\u\v}\in\K$ where
$r_{\u\v}\ne0$ only if $(\u,\v)\gdom(\s,\t)$. It follows that
$\{f_{\s\t}\}$ is a seminormal basis of~$\H(\K)$ in the sense of
\autoref{D:SeminormalBasis}.

For $\blam\in\Parts$ set
$[\blam]_t^!=\prod_{l=1}^\ell\prod_{r\ge1}[\lambda^{(l)}_r]_t^!\in\N[t]$.

\begin{Proposition}\label{P:MurphySeminormal}
  The basis $\set{f_{\s\t}|\s,\t\in\Std(\blam)\text{ for }\blam\in\Parts}$ is
  the $\bigast$-seminormal basis of $\HKK$ determined by the seminormal coefficient system
  defined in \autoref{Ex:MurphySNC} and the choices
  $$\gamma_{\tlam} = [\blam]_t^!\prod_{1\le l<m\le\ell}
  \prod_{(l,r,c)\in[\blam]}[\kappa_l-r+c-\kappa_m],$$
  for $\blam\in\Parts$.
\end{Proposition}

\begin{proof}This is equivalent to \cite[Theorem~2.11]{M:gendeg}
  in the non-degenerate case and to \cite[Proposition~6.8]{AMR} in the
  degenerate case, however, rather than translating the notation from these two
  papers it is easier to prove this directly.

  As noted above, $(\O,t)$ separates $\Std(\Parts)$ and
  $f_{\s\t}\equiv m_{\s\t}+\sum r_{\u\v}m_{\u\v}\pmod\Hlam$, for some $r_{\u\v}\in\K$
  where $r_{\u\v}\ne0$ only if $(\u,\v)\gdom(\s,\t)$. Therefore,
  in view of \autoref{E:SeminormalGamma},
  $\set{f_{\s\t}|(\s,\t)\in\SStd(\Parts)}$ is a $\bigast$-seminormal basis of $\H(\K)$.
  By \autoref{T:Seminormal}, this basis is determined by a seminormal
  coefficient system $\balpha$ and by a choice of scalars
  $\set{\gamma_{\tlam}|\blam\in\Parts}$. If $\t\gdom\v=\t(r,r+1)$
  then, by definition, $m_{\s\t}T_r=m_{\s\v}$. The transition matrix between the
  $\{m_{\s\t}\}$ and $\{f_{\s\t}\}$ is unitriangular so, in view of
  \autoref{T:Seminormal}, $f_{\s\t}T_r=f_{\s\v}-\frac1{[\rho_r(\t)]}f_{\s\t}$.
  Therefore, by \autoref{D:SNCS}(c), the seminormal coefficient system
  corresponding to the basis $\{f_{\s\t}\}$ is the one appearing in
  \autoref{Ex:MurphySNC}.

  It remains to determine the scalars $\set{\gamma_{\tlam}|\blam\in\Parts}$
  corresponding to $\{f_{\s\t}\}$. It is
  well-known, and easy to prove using the relations in~$\H$, that
  $x_\blam^2=[\blam]^!_tx_\blam$.  Therefore, by \autoref{L:JucysMurphyAction},
  $$f_{\tlam\tlam}^2\equiv [\blam]^!_t m_\blam u_\blam
        \equiv [\blam]^!_t\prod_{1\le l<m\le\ell}
           \prod_{(l,r,c)\in[\blam]}[\kappa_l-r+c-\kappa_m]
           \cdot m_\blam\pmod{\Hlam}.$$
  Hence,
  $\gamma_{\tlam}=[\blam]^!_t\prod_{1\le l<m\le\ell}
                   \prod_{(l,r,c)\in[\blam]}[\kappa_l-r+c-\kappa_m]$
  by \autoref{E:SeminormalGamma}.
\end{proof}

As noted after \autoref{T:MurphyBasis}, the Murphy basis
$\set{m_{\s\t}|(\s,\t)\in\SStd(\Parts)}$ of~$\H$ gives a basis
$\set{m_\t|\t\in\Std(\blam)}$ of each Specht module $\UnS^\blam$, for
$\blam\in\Parts$. For example, we can set $m_\t=m_{\tlam\t}+\Hlam$, for
$\t\in\Std(\blam)$. By \autoref{E:CellularInnerProduct}, the cellular basis
equips the Specht module~$\UnS^\blam$ with an inner product $\<\ ,\ \>$. The
matrix
$$\UnGram=\big(\<m_\s,m_\t\>\big)_{\s,\t\in\Std(\blam)}$$
is the \textbf{Gram matrix} of $\UnS^{\blam}$ with respect to the Murphy basis.
Similarly, the seminormal basis yields a second basis
$\set{f_\t|\t\in\Std(\blam)}$ of $\UnS^\blam(\K)$, where $f_\t=m_\t
F_\t=f_{\tlam\t}+\Hlam$, for $\t\in\Std(\blam)$. The transition matrix between
these two bases is unitriangular, so by \autoref{E:SeminormalGamma} we have
\begin{equation}\label{E:GramDetGammaProd}
\det\UnGram = \det\(\<f_\s,f_\t\>\) = \prod_{\t\in\Std(\blam)}\gamma_\t.
\end{equation}
This `classical' formula for $\det\UnGram$ is well-known as it is the
cornerstone used to prove the formula for $\det\UnGram$ as a rational
function in \cite[Theorem~3.35]{JM:cyc-Schaper}. The following definition will
allow us to give an `integral' closed formula for $\det\UnGram$.

\begin{Definition}\label{D:PartitionDegree}
  Suppose that $e\in\{0,2,3,4,\dots\}$, $p$ is a prime integer and that
  $\blam\in\Parts$ is a multipartition of~$n$. Define
  $$\deg_e(\blam)=\sum_{\t\in\Std(\blam)}\deg_e(\t)\qquad\text{and}\qquad
  \Deg_p(\blam)=\sum_{k\ge 1}\deg_{p^k}(\blam).$$
\end{Definition}

By definition, $\deg_e(\blam)$ and $\Deg_p(\blam)$ are integers that, \textit{a
priori}, could be positive, negative or zero. In fact, the next result shows that they
are always non-negative integers, although we do not know of a direct combinatorial
proof of this. By definition, the integers $\deg_e(\blam)$ and $\Deg_p(\blam)$
depend on $\charge$ and $e$. Our definitions ensure that the tableau degrees
$\deg_e(\t)$, for $\t\in\Std(\blam)$, coincide with \autoref{E:TableauDegree}
when $\Lambda=\Lambda_e(\charge)$.

For $k\in\N$, let $\Phi_k=\Phi_k(t)$ be the $k$th cyclotomic polynomial in~$t$. As is
well-known, these polynomials are pairwise distinct irreducible polynomials in~$\Z[t]$ and
\begin{equation}\label{E:CyclotomicPoly}
[n] = \prod_{1<d|n}\Phi_d(t),
\end{equation}
whenever $n\ge1$.

\begin{Theorem}\label{T:GramDetFactorization}
    Suppose that $\kappa_l-\kappa_{l+1}>n$, for $1\le l<\ell$, and that
    $\O=\Z[t,t^{-1}]$. Then
    $$\det\UnGram=t^{\ell(\blam)}\prod_{e\ge2}\Phi_e(t)^{\deg_e(\blam)},$$
    where $\ell(\blam)=\sum_{\t\in\Std(\blam)}\ell(d(\t))$.
\end{Theorem}

\begin{proof} As remarked above, $\det\UnGram=\prod_\t\gamma_\t$. Therefore,
  to prove the theorem it is enough to show that if $\t\in\Std(\blam)$ then
  $$\gamma_\t = t^{\ell(d(\t))}\prod_{e>1}\Phi_e(t)^{\deg_e(\t)}.$$
  We prove this by induction on the dominance ordering.

  Suppose first that $\t=\tlam$. Then \autoref{P:MurphySeminormal} gives an
  explicit formula for~$\gamma_{\tlam}$ and, using \autoref{E:TableauDegree}, it
  is straightforward to check by induction on~$n$ that our claim is true in this case.
  Suppose then that $\tlam\gdom\t$. Then we can write $\t=\s(r,r+1)$ for some
  $\s\in\Std(\blam)$ such that $\s\gdom\t$, and where $1\le r<n$.  Therefore,
  using induction, \autoref{C:GammaRecurrence} and the seminormal coefficient system
  of \autoref{P:MurphySeminormal},
  $$\gamma_\t =t^{\ell(d(\s))}
     \frac{[1+\rho_r(\s)][1+\rho_r(\t)]}{[\rho_r(\s)][\rho_r(\t)]}
     \prod_{e>1}\Phi_e(t)^{\deg_e(\s)}.$$
  By definition, $[k]=-t^k[-k]$, for any $k\in\Z$. Now $\rho_r(\s)=-\rho_r(\t)>0$
  by \autoref{E:ChargeSeparation}, so
  $$\frac{[1+\rho_r(\s)][1+\rho_r(\t)]}{[\rho_r(\s)][\rho_r(\t)]}
  =t\,\frac{[1+\rho_r(\s)][-\rho_r(\t)-1]}{[\rho_r(\s)][-\rho_r(\t)]}
  =t\,\prod_{e>1}\Phi_e(t)^{d_e},$$
  where, according to \autoref{E:CyclotomicPoly}, the integer $d_e$ is given in
  terms of the quiver $\Gamma_e$ by
  $$d_e=\begin{cases}
      -2,&\text{if }i_r=i_{r+1},\\
      2,&\text{if }i_r\leftrightarrows i_{r+1},\\
      1,&\text{if $i_r\leftarrow i_{r+1}$ or }i_r\rightarrow i_{r+1},\\
      0,&\text{otherwise}.
  \end{cases}$$
  Applying \autoref{L:TableauDegreeRecurrence} now completes the proof of
  our claim --- and hence proves the theorem.
\end{proof}

\begin{Remark}We can remove the factor $t^{\ell(\blam)}$ from
  \autoref{T:GramDetFactorization} by rescaling the generators
  $T_1,\dots,T_{n-1}$ so that the quadratic relations in
  \autoref{D:HeckeAlgebras} become $(T_r-t^{\frac12})(T_r+t^{-\frac12})$, for
  $1\le r<n$. Note that the integer $d_e$ in the proof of \autoref{T:Seminormal}
  is equal to the degree of the homogeneous generator $\psi_r e(\bi)$ in the
  cyclotomic KLR algebra~$\R$.
\end{Remark}

Setting $t=1$ gives the degenerate cyclotomic Hecke algebras.  As a special case,
the next result gives an integral closed formula for the Gram determinants of
the Specht modules of the symmetric groups.

\begin{Corollary}
    Suppose that $\kappa_l-\kappa_{l+1}>n$, for $1\le l<\ell$, and that
    $\O=\Z$ and $t=1$. Then
    $$\det\UnGram = \prod_{\substack{0<p\in\Z\\p\text{ prime}}}p^{\Deg_p(\blam)},$$
    for $\blam\in\Parts$.
\end{Corollary}

\begin{proof}This follows by setting $t=1$ in \autoref{T:GramDetFactorization} and using
  the following well-known property of the cyclotomic polynomials:
  $$\Phi_e(1) = \begin{cases}
                  p,&\text{if $e=p^k$ for some $k\ge1$},\\
                  1,&\text{otherwise}.
                \end{cases}
  $$
\end{proof}

\begin{Corollary}\label{C:PositiveDegrees}
  Suppose that $e\in\{0,2,3,4,5,\dots\}$ and that $p>0$ is an integer prime. Then
  $\deg_e(\blam)\ge0$ and
  $\Deg_p(\blam)\ge0$, for all $\blam\in\Parts$.
\end{Corollary}

\begin{proof}
   As the Murphy basis is defined over $\Z[t,t^{-1}]$, the Gram determinant
   $\det\UnGram$ belongs to~$\Z[t,t^{-1}]$. Therefore, $\deg_e(\blam)\ge0$
   whenever $e>1$
   by \autoref{T:GramDetFactorization}. Consequently, $\Deg_p(\blam)\ge0$. Finally,
   if $e\gg0$ then $\deg_0(\t)=\deg_e(\t)$ for any $\t\in\Std(\Parts)$, so
   $\deg_e(\blam)\ge0$ for $e\in\{0,2,3,4,\dots\}$ as claimed.
\end{proof}

The statement of \autoref{C:PositiveDegrees} is purely combinatorial so it
should have a direct combinatorial proof. We now give a second representation
theoretic proof of this result that suggests that a combinatorial proof may be
difficult.

A \textbf{graded set} is a set $D$ equipped with a \textbf{degree} function
$\deg\map D\Z$. Let $q$ be an indeterminate over $\Z$ and define
the \textbf{$q$-cardinality} and \textbf{degree} of~$D$ to be
$$|D|_q=\sum_{d\in D}q^{\deg d}\in\N[q,q^{-1}]
\quad\text{and}\quad
 \deg D=\sum_{d\in D}\deg d\in\Z.$$
If $D$ is a graded set and $z\in\Z$ let $q^zD$ be the graded set with
the same elements
as~$D$ but where the shifted degree function is shifted so that
$d\in D$ now has degree $z+\deg d$. More generally, if $f(q)\in\N[q,q^{-1}]$
let $f(q)D$ be the graded set that is the disjoint union of the appropriate
number of shifted copies of~$D$. For example $(2+q)D=D\sqcup D\sqcup qD$. By definition, $|f(q)D|_q=f(q)|D|_q$.

If $e\in\{0,2,3,4,\dots\}$ let $\Std_e(\blam)$ be the graded set with
elements $\Std(\blam)$ and degree function $\t\mapsto\deg_e(\t)$, for
$\t\in\Std_e(\blam)$.

Fix $e\in\{0,2,3,4,\dots\}$ and consider the Hecke algebra $\H(\C)$ over
$\C$ with Hecke parameter $\xi$, a primitive $e$th root of unity if
$e>0$ or a non-root of unity if $e=0$. Let $S^\blam$ be the graded
Specht module introduced in~\cite{BKW:GradedSpecht} (see
\autoref{S:GradedSpecht}), and let $D^\bmu=S^\bmu/\rad S^\bmu$ be the
graded simple quotient of $S^\bmu$, as in
\cite{HuMathas:GradedCellular}. Let $\Klesh$ be the set of
\textbf{Kleshchev multipartitions} so that
$\set{D^\bmu\<k\>|\bmu\in\Klesh\text{ and }k\in\Z}$ is a complete set of
non-isomorphic graded simple $\H$-modules. As recalled in
\autoref{S:GradedSpecht}, $S^\blam$ comes equipped with a homogeneous
basis $\set{\psi_\t|\t\in\Std_e(\blam)}$.  Let
$d_{\blam\bmu}(q)=[S^\blam{:}D^\bmu]_q$ be the corresponding graded
decomposition number.

Fix a total ordering $\prec$ on $\Std_e(\blam)$ that extends the
dominance ordering, such as the lexicographic ordering. Suppose that
$\bmu\in\Klesh$. By Gaussian elimination, there exists a graded subset
$\DStd(\bmu)$ of $\Std_e(\bmu)$ and a homogeneous basis
$\set{C_\t|\t\in\DStd(\bmu)}$ of~$D^\bmu$ such that
$$C_\t=\psi_\t+\sum_{\v\prec\t}c_{\t\v}\psi_\v+\rad S^\bmu,$$
for some $c_{\t\v}\in\C$ such that $c_{\t\v}\ne0$ only if
$\deg\v=\deg\t$ and $\res(\v)=\res(\t)$.  In particular,
$\Dim D^\blam=|\DStd(\blam)|_q$.  Repeating this argument, with the composition
factors that appear in successive layers of the radical filtration of~$S^\blam$,
shows that there exists a bijection of graded sets
$$\Theta_\blam : \Std_e(\blam)\bijection\bigsqcup_{\bmu\in\Klesh}
               d_{\blam\bmu}(q)\DStd(\bmu).$$
Now if $\bmu\in\Klesh$ then $D^\bmu\cong(D^\bmu)^\circledast$, so that
$\deg\DStd(\bmu)=0$. It follows that
$\deg q^z\DStd(\bmu)=z\dim\UnD^\bmu$, for $z\in\Z$. Therefore, using the
bijection $\Theta_\blam$,
$$\deg_e(\blam) = \deg\Std_e(\blam)
      =\sum_{\bmu\in\Klesh}d_{\blam\bmu}(q)\deg\DStd(\bmu)
      =\sum_{\bmu\in\Klesh}d_{\blam\bmu}'(1)\dim\UnD^\bmu,$$
where $d'_{\blam\bmu}(1)$ is the derivative of the graded decomposition number
$d_{\lambda\bmu}(q)$ evaluated at $q=1$. As we are working with the Hecke
algebra $\H(\C)$ in characteristic zero,
$d_{\blam\bmu}(q)\in\N[q]$ by \cite[Corollary~5.15]{BK:GradedDecomp}.
Consequently, $\deg_e(\blam)\ge0$. Hence, the (deep) fact that
$d_{\blam\bmu}(q)\in\N[q]$ leads to an alternative proof of
\autoref{C:PositiveDegrees}.

In characteristic zero the graded cyclotomic Schur algebras is Koszul by
\cite[Theorem~C]{HuMathas:QuiverSchurI} when $e=0$ and by
\cite{Maksimau:QuiverSchur} and \cite[Proposition
7.8,7.9]{RouquierShanVaragnoloVasserot} in general.  This implies that
the Jantzen and grading filtrations of the graded Weyl modules, and
hence of the graded Specht modules, coincide. Therefore,
\autoref{C:PositiveDegrees} is compatible with this Koszulity Conjecture
via Ryom-Hansen's~\cite[Theorem~1]{RyomHansen:Schaper} description of
the Jantzen sum formula; see also
\cite[Theorem~2.11]{Yvonne:Conjecture}.

The construction of the sets $\DStd(\bmu)$ given above is not unique because it
involves many choices. It natural to ask if there is a canonical choice of basis
for~$S^\blam$ that uniquely determines the sets $\DStd(\bmu)$ and the
bijections $\Theta_\blam$. For level~2 such bijections are implicit in
\cite[\Sect9]{BrundanStroppel:KhovanovIII} when
$e=0$ and \cite[Appendix]{HuMathas:QuiverSchurI} generalizes this to include the
cases when $e>n$. It is interesting to note that the sets $\DStd(\bmu)$,
together with the bijections $\Theta_\blam$, determine the graded decomposition
numbers. More explicitly, if $\s\in\DStd(\bmu)$ then
$$d_{\blam\bmu}(q)= \sum_{\t\in\Theta_{\blam}^{-1}(\s)}q^{\deg\t-\deg\s},$$
where we abuse notation and let $\Theta_\blam^{-1}(\s)$ be the set of tableaux
in $\Std_e(\blam)$ that are mapped onto a (shifted) copy of~$\s$ by~$\Theta_\blam$.
In particular, we can take $\s=\tmu$ because it is easy to see that we must have
$\tmu\in\DStd(\bmu)$ whenever $\bmu\in\Klesh$. Hence, we have shown that
the KLR tableau combinatorics leads to closed combinatorial formulas for
the parabolic Kazhdan-Lusztig polynomials $d_{\blam\bmu}(q)$, and the
graded simple dimensions $\Dim D^\bmu$: both families of
polynomials can be described as the $q$-cardinalities of graded sets of tableaux.

\section{Integral Quiver Hecke algebras}\label{Chap:QuiverHeckeAlgebras}
The Seminormal Basis \autoref{T:Seminormal} compactly describes much of the
semisimple representation theory of~$\HKK$. For symmetric groups,
Murphy~\cite{M:Nak} showed that seminormal bases can also be used to study the
non-semisimple representation theory. Murphy's ideas were extended to the
cyclotomic Hecke algebras in~\cite{M:gendeg,M:seminormal}. In this section we
further extend Murphy's ideas to connect seminormal bases and the KLR grading
on~$\H$.

\subsection{Lifting idempotents}
As \autoref{S:SeminormalForm}, we continue to assume that $\charge$ satisfies
\autoref{E:ChargeSeparation} and that $(\K,t)$ separates $\Std(\Parts)$, where
$\K$ is a field and $0\ne t\in\K$. If $\O$ is a subring of~$\K$ then we identify
$\HO$ with the obvious $\O$-subalgebra of $\H(\K)$ so that
$\H(\K)\cong\HO\otimes_\O\K$ as $\K$-algebras.

Let $J(\O)$ be the \textbf{Jacobson radical} of $\O$, the intersection of all of
the maximal ideals of~$\O$.

\begin{Definition}\label{D:IdempotentSubring}
 Suppose that $\O$ is a subring of $\K$ and $t\in\O^{\times}$. Then
  $(\O,t)$ is an \textbf{$e$-idempotent subring} of $\K$ if the following
  hold:
  \begin{enumerate}
    \item $(\O,t)$ separates $\Std(\Parts)$;
    \item $[k]_t$ is invertible in $\O$ whenever $k\not\equiv0\pmod e$, for  $k\in\Z$; and
    \item $[k]_t\in J(\O)$ whenever $k\in e\Z$.
\end{enumerate}
\end{Definition}

When $e$ and $t$ are understood, we simply call $\O$ an idempotent subring. Note
that $\K$ contains the field of fractions of $\O$, so
\autoref{D:IdempotentSubring}(a) ensures that $\H(\K)$ is semisimple and that it
has a seminormal basis. Until further notice, we fix such a $\bigast$-seminormal
basis $\{f_{\s\t}\}$, together with the corresponding seminormal coefficient
system $\balpha$ and $\gamma$-coefficients.

Let $(\O,t)$ be an $e$-idempotent subring and suppose $c\not\equiv d\pmod e$,
for $c,d\in\Z$.  Then $[c]-[d]=t^d[c-d]$ is invertible in~$\O$. We use this fact
below without mention.

\begin{Examples}{The following local rings are all examples of idempotent subrings.}
  \label{Ex:indeterminate}
  \item Suppose that $\K=\Q$ and $t=1$. Then $(\K,t)$ separates $\Std(\Parts)$
  and $\O=\Z_{(p)}$ is a $p$-idempotent subring of~$\Q$ for any prime~$p$.

  \item Let $K$ be any field and set $\K=K(x)$, where $x$ is an indeterminate
 over~$K$, and $t=x+\xi$, where $\xi$ is a primitive $e$th root of unity in~$K$.
 Then $\O=K[x]_{(x)}$ is an $e$-idempotent subring of $\K$.

  \item Let $\K=\Q(x,\xi)$, where $x$ is an indeterminate over~$\Q$ and
  $\xi=\exp(2\pi i/e)$ is a primitive $e$th root of unity in~$\C$. Let
  $t=x+\xi$. Then $(\K,t)$ separates $\Std(\Parts)$ and $\O=\Z[x,\xi]_{(x)}$ is
  an $e$-idempotent subring of~$\K$.

  \item Maintain the notation of the last example and let $p>1$ be a prime not
  dividing~$e$. Let $\Phi_{e,p}(x)$ be a polynomial in~$\Z[x]$ whose reduction
  modulo~$p$ is the minimum polynomial of a primitive~$e$th root of unity in an
  algebraically closed field of characteristic~$p$. Then
  $\O=\Z[x,\xi]_{(x,p,\Phi_{e,p}(\xi))}$ is an $e$-idempotent subring of $\C(x)$.
\end{Examples}

Suppose that $\bi\in I^n$ and set $\Std(\bi)=\set{\t\in\Std(\Parts)|\res(\t)=\bi}$. Define
the \textbf{residue idempotent} $\fo$ by
\begin{equation}\label{E:ResidueIdempotent}
  \fo =\sum_{\t\in\Std(\bi)}F_\t .
\end{equation}
By \autoref{T:Seminormal}, $\fo$ is an idempotent in $\HKK$. In the
rest of this section, we fix a seminormal basis $\{f_{\s\t}\}$ of
$\HKK$ that is determined by a seminormal coefficient system
$\{\alpha_r(\s)\}$ and a choice of $\gamma_{\tlam}$. Then we have
that $\fo=\sum_{\t\in\Std(\bi)}\frac1{\gamma_\t}f_{\t\t}$.

\begin{Lemma} \label{L:MurphyIdempotent}
  Suppose that $\O$ is an idempotent subring of~$\K$ and that $\bi\in I^n$. Then
  $\fo\in\L(\O)$.
  In particular, $\fo$ is an idempotent in~$\HO$.
\end{Lemma}

\begin{proof}
  This result is proved when~$\O$ is a discrete valuation ring in
  \cite[Lemma~4.2]{M:seminormal}, however, our weaker assumptions necessitate a
  different proof. Motivated, in part, by the proof of \cite[Theorem~2.1]{M:Nak}, if
  $\t\in\Std(\bi)$ define
  $$F_\t' = \prod_{k=1}^n \prod_{\substack{c\in\Cont\\\cont_k(\t)\not\equiv c\pmod e}}
      \frac{L_k-[c]}{[\cont_k(\t)]-[c]}.$$
  Since $\O$ is an $e$-idempotent subring, $F_\t'\in\L(\O)\subset\HO$. By
  \autoref{T:Seminormal}, $\sum_{\s\in\Std(\Parts)}F_\s$ is the identity element
  of~$\HKK$ so, using \autoref{E:SeminormalForm}, we see that
  $$F_\t' = \sum_{\s\in\Std(\Parts)}F_\t' F_\s
          = \sum_{\s\in\Std(\Parts)} a_{\s\t} F_\s,
  $$
  where $a_{\s\t}=\prod_{k,c} ([\cont_k(\s)]-[c])/([\cont_k(\t)]-[c])\in\O$. In
  particular, $a_{\t\t}=1$. If $\s\notin\Std(\bi)$ then there exists an integer
  $k$ such that $\res_k(\s)\ne\res_k(\t)$, so
  $[\cont_k(\s)]-[\cont_k(\t)]\in\O^\times$ and $a_{\s\t}=0$. Therefore,
  $F_\t'=\sum_{\s\in\Std(\bi)}a_{\s\t}F_\s$. Consequently,
  $\fo F_\t' = F_\t' = F_\t'\fo$ by~\autoref{E:SeminormalGamma}. Notice that
  $F_\t'F_\s'=F_\s'F_\t'$ because $\L(\K)$ is a commutative subalgebra
  of~$\HKK$. Therefore,
  $$\prod_{\t\in\Std(\bi)}(\fo-F_\t')
         =\fo\qquad+
         \sum_{\substack{\t_1,\dots,\t_k\in\Std(\bi)\\\text{distinct with }k>0}}
         (-1)^kF_{\t_1}'F_{\t_2}'\dots F_{\t_k}'.$$
On the other hand, since $\fo=\sum_{\s\in\Std(\bi)}F_\s$ and $a_{\t\t}=1$,
  $$\prod_{\t\in\Std(\bi)}(\fo-F_\t')
      =\prod_{\t\in\Std(\bi)}\sum_{\substack{\s\in\Std(\bi)\\\s\ne\t}}
           (1-a_{\s\t})F_\s = 0,$$
because $F_\s F_\t=0$ whenever $\s\ne\t$ by \autoref{E:SeminormalGamma}. Combining the last
two equations,
  $$\fo=\sum_{\substack{\t_1,\dots,\t_k\in\Std(\bi)\\\text{distinct with }k>0}}
          (-1)^{k+1}F_{\t_1}'F_{\t_2}'\dots F_{\t_k}'.$$
  In particular, $\fo\in\L(\O)$ as we wanted to show.
\end{proof}

\begin{Corollary}\label{C:Complete}
    Suppose that $\O$ is an idempotent subring of~$\K$. Then
  $\set{\fo|\bi\in I^n}$ is a complete set of pairwise orthogonal idempotents
  in~$\HO$.
\end{Corollary}

\begin{proof}
  By \autoref{T:Seminormal}, $\set{F_\t|\t\in\Std(\Parts)}$ is a complete set of
  pairwise orthogonal idempotents in $\HKK$. Hence, the result follows from
  \autoref{L:MurphyIdempotent}.
\end{proof}

If $\phi\in\O[X_1,\dots,X_n]$ is a polynomial in indeterminates $X_1,\dots,X_n$
over~$\O$ then set $\phi(L)=\phi(L_1,\dots,L_n)\in\L(\O)$. If $\s$ is a tableau
let $\phi(\s)=\phi([c_1(\s)],\dots,[c_n(\s)])$ be the scalar in~$\O$ obtained by
evaluating the polynomial~$\phi$ on the contents of~$\s$; that is, setting
$X_1=[c_1(\s)], \dots, X_n=[c_n(\s)]$. Then,
$\phi(L)f_{\s\t}=\phi(\s)f_{\s\t}$, for all $(\s,\t)\in\SStd(\Parts)$.

Ultimately, the next result will allow us to `renormalise' intertwiners of
the residue idempotents~$\fo$, for $\bi\in I^n$, so that they depend only on~$e$
rather than on~$\xi$.

\begin{Proposition}\label{P:Invertible}
  Suppose that $\bi\in I^n$ and $\phi\in\O[X_1,\dots,X_n]$ is a polynomial such that
  $\phi(\t)$ is invertible in~$\O$, for all $\t\in\Std(\bi)$. Then
  $$f^\phi_\bi=\sum_{\t\in\Std(\bi)}\frac1{\phi(\t)}F_\t\in\L(\O).$$
  In particular, $f^\phi_\bi\in\HO$.
\end{Proposition}

\begin{proof}By assumption, $\phi(\s)$ is invertible
  in~$\O$ for all $\s\in\Std(\bi)$. In particular,~$f^\phi_\bi$ is a well-defined
  element of~$\L(\K)$. It remains to show that $f^\phi_\bi\in\L(\O)$.

  As in \autoref{L:MurphyIdempotent}, for each $\t\in\Std(\bi)$ define
  $$F_\t'=\prod_{\substack{c\in\Cont\\\cont_k(\t)\not\equiv c\pmod e}}
              \frac{L_k-[c]}{[\cont_k(\t)]-[c]}\in\L(\O),$$
  and write $F_\t'=\sum_{\s\in\Std(\bi)}a_{\s\t}F_\s$ for some $a_{\s\t}\in\O$.
  Recall from the proof of \autoref{L:MurphyIdempotent} that $a_{\t\t}=1$.

  Motivated by the definition of $F_\t'$, set
  $F^\phi_\t=\tfrac{\phi(L)}{\phi(\t)}F'_{\t}$. Then
  $F_\t^\phi\in\L(\O)$ and
  $$F_\t^\phi=\sum_{\s\in\Std(\bi)}a_{\s\t}\frac{\phi(L)}{\phi(\t)}F_\s
  =F_\t+\sum_{\substack{\s\in\Std(\bi)\\\s\ne\t}}\frac{a_{\s\t}\phi(\s)}{\phi(\t)}F_\s$$
  by \autoref{E:SeminormalForm}. Consequently, $F^\phi_\t \fo=F^\phi_\t=\fo F^\phi_\t$.
  The idempotents $\set{F_\s|\s\in\Std(\bi)}$ are pairwise orthogonal, so
  $$f^\phi_\bi F_\t^\phi
    =\Big(\sum_{\s\in\Std(\bi)}\frac1{\phi(\s)}F_\s\Big)
     \Big(\sum_{\s\in\Std(\bi)}\frac{a_{\s\t}\phi(\s)}{\phi(\t)}F_\s\Big)
    =\sum_{\s\in\Std(\bi)}\frac{a_{\s\t}}{\phi(\t)}F_\s
    =\frac1{\phi(\t)}F_\t'.$$
  Therefore, $f_\bi^\phi F_\t^\phi\in\L(\O)$, for all $\t\in\Std(\bi)$. By
  \autoref{E:SeminormalForm}, $f_\bi^\phi\fo=f_\bi^\phi=\fo f_\bi^\phi$, so this
  implies that $f^\phi_\bi(\fo-F^\phi_\t)\equiv f^\phi_\bi\pmod{\L(\O)}$. Hence,
  working modulo $\L(\O)$,
  $$
  f^\phi_\bi \equiv f^\phi_\bi\prod_{\t\in\Std(\bi)}\(\fo-F_\t^\phi)
    =f^\phi_\bi\prod_{\t\in\Std(\bi)}\sum_{\substack{\s\in\Std(\bi)\\\s\ne\t}}
               \frac{a_{\s\t}\phi(\s)}{\phi(\t)}F_\s=0,
  $$
  where the last equality follows using the orthogonality of the idempotents
  $F_\s$ once again. Therefore, $f^\phi_\bi\in\L(\O)$,  completing the proof.
\end{proof}

Let $\phi$ be a polynomial in $\O[X_1,\dots,X_n]$ satisfying the assumptions of
\autoref{P:Invertible}. Then $\phi(L)f^\phi_\bi=\fo=f^\phi_\bi\phi(L)$ by
\autoref{E:SeminormalForm}. Abusing notation, in this situation we write
$$\frac1{\phi(L)}\fo= f^\phi_\bi
         =\sum_{\s\in\Std(\bi)}\frac1{\phi(\s)}F_\s
         =\fo\frac1{\phi(L)}\in\L(\O).$$
Note that, either by direction calculation or because $\L$ is commutative, we are
justified in writing $\fo\tfrac1{\phi(L)}=\tfrac1{\phi(L)}\fo$.

We need the following three special cases of \autoref{P:Invertible}.
For $1\le r<n$ define $M_r=1-L_r+tL_{r+1}$ and $M_r'=1+tL_r-L_{r+1}$, for $1\le
r<n$. Applying the definitions, if $(\s,\t)\in\SStd(\Parts)$ then
\begin{equation}\label{E:Mr}
  M_r f_{\s\t}=t^{c_r(\s)}[1-\rho_r(\s)]f_{\s\t}\quad\text{and}\quad
  M_r'f_{\s\t}= t^{c_{r+1}(\s)}[1+\rho_r(\s)]f_{\s\t}.
\end{equation}
Our main use of \autoref{P:Invertible} is the following application
that corresponds to taking $\phi(L)$ be to $L_r-L_{r+1}$, $M_r$ and $M_r'$,
respectively.

\begin{Corollary}\label{C:Inverting}
  Suppose that $\O$ is an $e$-idempotent subring, $1\le r<n$ and $\bi\in I^n$.
  \begin{enumerate}
    \item If $i_r\ne i_{r+1}$ then
    $\dfrac1{L_r-L_{r+1}}\fo=\Sum_{\t\in\Std(\bi)}\frac{t^{-c_{r+1}(\t)}}{[\rho_r(\t)]}F_\t
          \in\L(\O)$.
    \item If $i_r\ne i_{r+1}+1$ then
    $\dfrac1{M_r}\fo=\Sum_{\t\in\Std(\bi)}\frac{t^{-c_r(\t)}}{[1-\rho_r(\t)]}F_\t\in\L(\O)$.
    \item If $i_r\ne i_{r+1}-1$ then
    $\dfrac1{M_r'}\fo=\Sum_{\t\in\Std(\bi)}\frac{t^{-c_{r+1}(\t)}}{[1+\rho_r(\t)]}F_\t\in\L(\O)$.
\end{enumerate}
\end{Corollary}

\subsection{Intertwiners}
By \autoref{T:BKiso}, if $K$ is a field then the KLR generators of~$\H(K)$ satisfy
$\psi_r e(\bi)=e(s_r\cdot\bi)\psi_r$. This section defines analogous elements
in~$\HO$ that intertwine the residue idempotents $\fo$, for $\bi\in I^n$.

\begin{Lemma}\label{L:IntertwinerEq}
  Suppose that $i_r=i_{r+1}$, for some $\bi\in I^n$ and $1\le r<n$. Then
  $T_r\fo=\fo T_r$.
\end{Lemma}

\begin{proof}This follows directly from the Seminormal Basis
  \autoref{T:Seminormal}. In more detail, note that if $\t\in\Std(\bi)$ then $r$
  and $r+1$ cannot appear in the same row or in the same column of~$\t$.
  Therefore,
  $$T_r\fo-\fo T_r
  =\sum_{\t\in\Std(\bi)}\frac1{\gamma_\t}\Big(T_rf_{\t\t}-f_{\t\t}T_r\Big)
  =\sum_{\substack{\t,\v\in\Std(\bi)\\\v=\t(r,r+1)}}
  \Big(\frac{\alpha_r(\t)}{\gamma_\t}-\frac{\alpha_r(\v)}{\gamma_\v}\Big)f_{\v\t},$$
  by \autoref{E:SeminormalForm}. By \autoref{C:GammaRecurrence}, if $\v=\t(r,r+1)$ then
  $\alpha_r(\t)\gamma_\v=\alpha_r(\v)\gamma_\t$. Hence, $T_r\fo=\fo T_r$ as claimed.
\end{proof}

\begin{Remark}
In the special case of the symmetric groups,
Ryom-Hansen~\cite[\Sect3]{RyomHansen:ModularSeminormal} has proved an analogue of
\autoref{L:IntertwinerEq}.
\end{Remark}

Using \autoref{E:SeminormalForm}, it is easy to verify that $T_r\fo\ne \fo[\bj] T_r$
if $\bj=s_r\cdot\bi\ne\bi$, for $1\le r<n$ and $\bi\in I^n$. The following
elements will allow us to correct for this.

\begin{Lemma}\label{L:LTaction}
  Suppose that $(\s,\t)\in\SStd(\Parts)$ and $1\le r<n$. Let $\u=\s(r,r+1)$. Then
  $(T_rL_r-L_rT_r)f_{\s\t}=\alpha_r(\s)t^{c_{r+1}(\s)}[\rho_r(\s)]f_{\u\t}.$
\end{Lemma}

\begin{proof}Using \autoref{E:SeminormalForm} we obtain
  $$
    (T_rL_r-L_rT_r)f_{\s\t}
        =\alpha_r(\s)\([c_r(\s)]-[c_{r+1}(\s)]\)f_{\u\t}
        =\alpha_r(\s)t^{c_{r+1}(\s)}[\rho_r(\s)]f_{\u\t},
  $$
  where, as usual, we set $f_{\u\t}=0$ if $\u$ is not standard.
\end{proof}

Applying the $\bigast$ anti-involution,
$f_{\s\t}(T_rL_r-L_rT_r)=-\alpha_r(\t)t^{c_{r+1}(\t)}[\rho_r(\t)]f_{\s\v}$,
where $\v=\t(r,r+1)$.

\begin{Lemma}\label{L:InterwinderNe}
  Suppose that $i_r\ne i_{r+1}$, for some $\bi\in I^n$ and $1\le r<n$. Set
  $\bj=s_r\cdot\bi$. Then
  $ (T_rL_r-L_rT_r\)\fo =\fo[\bj](T_rL_r-L_rT_r\).$
\end{Lemma}

\begin{proof} By definition, $\fo=\sum_{\s\in\Std(\bi)}\tfrac1{\gamma_\s}f_{\s\s}$ so,
  by \autoref{L:LTaction},
  \begin{align*}
    (T_rL_r-L_rT_r)\fo &=\sum_{\s\in\Std(\bi)}\frac1{\gamma_\s}(T_rL_r-L_rT_r)f_{\s\s}\\
          &=\sum_{\substack{\s\in\Std(\bi)\\ \u=\s(r,r+1)\in\Std(\Parts)}}
                 \frac{\alpha_r(\s)t^{c_{r+1}(\s)}[\rho_r(\s)]}{\gamma_\s}f_{\u\s}.
  \end{align*}
  Note that if $\s\in\Std(\bi)$ and $\u=\s(r,r+1)$ is standard then
  $\s\in\Std(\bj)$. Similarly,
  $$\fo[\bj](T_rL_r-L_rT_r)=
         \sum_{\substack{\u\in\Std(\bj)\\\s=\u(r,r+1)\in\Std(\bi)}}
         -\frac{\alpha_r(\u)t^{c_{r+1}(\u)}[\rho_r(\u)]}{\gamma_\u}f_{\u\s}.$$

  By \autoref{E:SeminormalForm}, the tableaux in $\Std(\bi)$ and $\Std(\bj)$ that
  have $r$ and $r+1$ in the same row or in the same column do not contribute to
  the right hand sides of either of the last two equations. Moreover, the map
  $\s\mapsto\u=\s(r,r+1)$ defines a bijection from the set of tableaux in
  $\Std(\bi)$ such that $r$ and $r+1$ appear in different rows and columns to
  the set of tableaux in $\Std(\bj)$ that have $r$ and $r+1$ in different rows
  and columns. In particular, $(T_rL_r-L_rT_r)\fo=0$ if and only if
  $\fo[\bj](T_rL_r-L_rT_r)=0$.

  To complete the proof suppose that $\s\in\Std(\bi)$ and
  that $\u=\s(r,r+1)\in\Std(\bj)$.
  Now, $\alpha_r(\u)\gamma_\s=\alpha_r(\s)\gamma_\u$, by \autoref{C:GammaRecurrence},
  and $\rho_r(\u)=-\rho_r(\s)$, by definition. So
  $$ \frac{-\alpha_r(\u)t^{c_{r+1}(\u)}[\rho_r(\u)]}{\gamma_\u}
           =\frac{-\alpha_r(\s)t^{c_{r}(\s)}[-\rho_r(\s)]}{\gamma_\s}
           =\frac{\alpha_r(\s)t^{c_{r+1}(\s)}[\rho_r(\s)]}{\gamma_\s}.$$
  Hence, comparing the equations above,
  $(T_rL_r-L_rT_r\)\fo=\fo[\bj](T_rL_r-L_rT_r\)$ as required.
\end{proof}

Recall the definitions of $M_r$ and $M_r'$ from \autoref{E:Mr}, for $1\le r<n$.
We finish this section by giving the commutation relations for the elements
$M_r$, $M_r'$, $(1+T_r)$ and $(T_rL_r-L_rT_r)$. These will be important later.

\begin{Lemma}\label{L:MrCommute}
  Suppose that $1\le r<n$. Then
  $$(T_rL_r-L_rT_r)M_r=M_r'(T_rL_r-L_rT_r)
  \quad\text{and}\quad (T_r-t)M_r=M_r'(1+T_r).
  $$
\end{Lemma}

\begin{proof}Both formulas can be proved by applying the relations in
  \autoref{D:HeckeAlgebras}. Alternatively, suppose that $(\s,\t)\in\SStd(\Parts)$ and
  set $\v=\t(r,r+1)$. Then, by~\autoref{E:Mr} and \autoref{L:LTaction},
  \begin{align*}
  f_{\s\t}(T_rL_r-L_rT_r)M_r
  &=-\alpha_r(\t)t^{2c_r(\v)}[\rho_r(\t)][1+\rho_r(\t)]f_{\s\v}\\
           &=f_{\s\t}M_r'(T_rL_r-L_rT_r),
  \end{align*}
  where the last equality follows because $c_r(\v)=c_{r+1}(\t)$ and
  $c_{r+1}(\v)=c_r(\t)$.  As the regular representation is a faithful,
  this implies the first formula. The second formula can be proved similarly.
\end{proof}

\subsection{The integral KLR generators}\label{S:KLRGens}
In \autoref{L:IntertwinerEq} and \autoref{L:InterwinderNe}, we have found
elements in $\HO$ that intertwine the residue idempotents~$\fo$.  These
intertwiners are not quite the elements that we need, however, because they still
depend on~$t$, rather than just on~$e$. To remove this dependence on~$t$ we
will use \autoref{P:Invertible} to renormalise these elements.

By \autoref{L:MurphyIdempotent}, if $h\in\HO$ then $h=\sum_{\bi\in
I^n}h\fo$, so that $h$ is completely determined by its projections onto the
spaces $\HO \fo$. We use this observation to define analogues of the KLR
generators in~$\HO$.

Recall from \autoref{E:Mr} that $M_r=1-L_r+tL_{r+1}$. By \autoref{C:Inverting},
if $i_r\ne i_{r+1}+1$ then $M_r$ acts invertibly on $\fo\HO$ so
$\tfrac1{M_r}\fo$ is a well-defined element of~$\HO$.

As in the introduction, define an embedding $I\hookrightarrow \Z; i\mapsto\i$ by
defining $\i$ to be the smallest non-negative integer such that $i=\i+e\Z$, for
$i\in I$.

\begin{Definition}\label{D:KLRLift}
  Suppose that $1\le r<n$. Define elements
  $\psi_r^\O=\sum_{\bi\in I^n}\psi_r^\O \fo$ in~$\HO$ by
  $$\psi_r^\O \fo=\begin{cases}
    (T_r+1)\frac{t^{\i_r}}{M_r}\fo,&\text{if }i_r=i_{r+1},\\
    (T_rL_r-L_rT_r)t^{-\i_r}\fo, &\text{if }i_r=i_{r+1}+1,\\
    (T_rL_r-L_rT_r)\frac1{M_r}\fo, &\text{otherwise}.
  \end{cases}$$
  If $1\le r\le n$ then define $\yo_r=\sum_{\bi\in I^n}t^{-\i_r}(L_r-[\i_r])\fo$.
\end{Definition}

The order of the terms in the definition of $\psio_r$ matters because $M_r$ does
not commute with $T_r+1$ or with $T_rL_r-L_rT_r$ (see \autoref{L:MrCommute}),
although $M_r$ does commute with $\fo$. Notice that $\psio_r$ is independent of
the choice of seminormal coefficient system because the residue idempotents
$\fo$ are independent of this choice.

One subtlety of \autoref{D:KLRLift}, which we will pay for later, is that it
makes use of the embedding $I\hookrightarrow\Z$ in order to give meaning to
expressions like $t^{\pm\i_r}$.

\begin{Remark}
Unravelling the definitions, the element $\psio_r\otimes_\O1_K$ is a scalar
multiple of the choice of KLR generators for~$\HK$ made by Stroppel and
Webster~\cite[(27)]{StroppelWebster:QuiverSchur}. Similarly,
$\yo_r\otimes_\O1_K$ is a multiple of the KLR generator~$y_r$ defined by
Brundan and Kleshchev~\cite[(4.21)]{BK:GradedKL}.
\end{Remark}

\begin{Proposition}\label{P:OGeneration} The algebra $\HO$ is generated by the
  elements
  $$\set{\fo|\bi\in I^n}\cup\set{\psio_r|1\le r<n}\cup\set{\yo_r|1\le r\le n}.$$
\end{Proposition}

\begin{proof}Let $H$ be the $\O$-subalgebra of $\HO$ generated by the elements
  in the statement of the proposition. We need to show that $H=\HO$. Directly
  from the definitions, if $1\le r\le n$ then
  $L_r=\sum_\bi (t^{\i_r}\yo_r+[i_r])\fo\in H$. Therefore, the Gelfand-Zetlin
  algebra $\L(\O)$ is contained in~$H$. Consequently,  $M_r\in H$, for $1\le r<n$.
  By~\autoref{D:HeckeAlgebras}, $L_rT_r-T_rL_r=T_r(L_{r+1}-L_r)-1+(1-t)L_{r+1}$.
  By \autoref{C:Inverting}(a),  if $i_r\ne i_{r+1}$ then
  $\tfrac1{L_r-L_{r+1}}\fo\in\L(\O)\subseteq H$.  Therefore, since
  $M_r$ and $\fo$ commute, we can write
  $$T_r \fo = \begin{cases}
    \(t^{-\i_r}\psio_rM_r-1\)\fo,&\text{if }i_r=i_{r+1},\\
    \(-t^{\i_r}\psio_r+1+(t-1)L_{r+1}\)\frac1{L_{r+1}-L_r}\fo,
         &\text{if }i_r=i_{r+1}+1\\
       \(-\psio_rM_r+1+(t-1)L_{r+1}\)\frac1{L_{r+1}-L_r}\fo,  &\text{otherwise}.
  \end{cases}$$
  by \autoref{D:KLRLift}. Hence, $T_r=\sum_\bi T_r\fo\in H$. As
  $T_1,\dots,T_{n-1},L_1,\dots,L_n$ generate~$\HO$ this implies that $H=\HO$,
  completing the proof.
\end{proof}

We now use the seminormal form to  show that the elements in the statement of
\autoref{P:OGeneration} satisfy most of the relations of
\autoref{D:QuiverRelations}.

\begin{Lemma} \label{L:PsiCommute}
  Suppose that $1\le r<n$ and $\bi\in I^n$. Then
  $\psio_r \fo=\fo[\bj]\psio_r$, where $\bj=s_r\cdot\bi$.
\end{Lemma}

\begin{proof}
By \autoref{L:MurphyIdempotent} and \autoref{P:Invertible}, respectively,
$M_r$ and $\fo$ both belong to~$\L(\O)$, which is a commutative algebra.
Therefore, $\frac1{M_r}\fo$ and $\fo$ commute.
If $i_r=i_{r+1}$ then
  $$\psio_r\fo=(T_r+1)\frac{t^{\i_r}}{M_r}\fo
              =(T_r+1)\fo\frac{t^{\i_r}}{M_r}\fo
              =\fo(T_r+1)\frac{t^{\i_r}}{M_r}\fo
              =\fo\psio_r,$$
  where the third equality comes from \autoref{L:IntertwinerEq}. The remaining
  cases follow similarly using \autoref{L:InterwinderNe}.
\end{proof}

As we will work with right modules we need the right-handed analogue of
\autoref{D:KLRLift}. Note that if $i_r\ne i_{r+1}+1$ then
$\fo\tfrac1{M_r}=\tfrac1{M_r}\fo\in\HO$ by \autoref{P:Invertible}. Similarly, if
$i_r\ne i_{r+1}-1$ then $\fo\frac1{M_r'}=\frac1{M_r'}\fo\in\HO$.  It follows
that all of the expressions in the next lemma make sense.

\begin{Lemma}\label{L:RightPsi}
  Suppose $1\le r<n$ and $\bi\in I^n$. Then
  $$\fo\psio_r = \begin{cases}
    \fo\frac{t^{\i_{r+1}}}{M_r'}(T_r-t),&\text{if }i_i=i_{r+1},\\
    \fo(T_rL_r-L_rT_r)t^{-\i_{r+1}},&\text{if }i_r=i_{r+1}-1,\\
    \fo\frac1{M_r'}(T_rL_r-L_rT_r),&\text{otherwise}.
  \end{cases}
  $$
\end{Lemma}

\begin{proof} By \autoref{L:PsiCommute}, $\fo\psio_r=\fo\psio_r\fo[\bj]$ where
  $\bj=s_r\cdot\bi$.
  Therefore,
  $$\fo\psio_r = \begin{cases}
    \fo(1+T_r)\frac{t^{\i_{r+1}}}{M_r}\fo[\bj],&\text{if }i_i=i_{r+1},\\
   \fo(T_rL_r-L_rT_r)t^{-\i_{r+1}}\fo[\bj], &\text{if }i_r=i_{r+1}-1,\\
   \fo(T_rL_r-L_rT_r)\frac1{M_r}\fo[\bj],  &\text{otherwise}.
  \end{cases}
  $$
  To complete the proof apply \autoref{L:MrCommute}.
\end{proof}

\def\UnSpace{\hspace*{2mm}}
\begin{Lemma}\label{L:fRelations}
  Suppose that $\bi,\bj\in I^n$ and $1\le r,s\le n$. Then
  $$\sum_{\bi\in I^n}\fo=1, \UnSpace
   \fo \fo[\bj]=\delta_{\bi\bj}\fo,\UnSpace
   \yo_r \fo=\fo \yo_r\UnSpace\text{and}\UnSpace \yo_r\yo_s=\yo_s\yo_r.
   $$
   Moreover, if $s\neq r, r+1$ then $\psio_r \yo_s=\yo_s\psio_r$,
   for $1\le r<n$ and $1\le s\le n$.
\end{Lemma}

\begin{proof}
  The elements $\fo$, for $\bi\in I^n$, form a complete set of pairwise
  orthogonal idempotents by \autoref{L:MurphyIdempotent}, which gives the first
  two relations. Since $y_r, \fo\in\L(\O)$ and $\L(\O)$ is a commutative
  algebra, all of the elements $\fo$, $\yo_r$ and~$\yo_s$ commute.

  Now suppose that $s\neq r, r+1$. Then $y_s^{\O}$ commutes with
  $\frac{1}{M_r}f_{\bi}^{\O}$ and with~$T_r$. Hence, $\psio_r
  f_{\bi}^{\O}\yo_s=\yo_s\psio_r f_{\bi}^{\O}$, for any $\bi\in I^n$. Therefore,
  $\psio_r \yo_s=\yo_s\psio_r$.
\end{proof}

\begin{Lemma}\label{L:y1Cyclic}
  Suppose that $\bi\in I^n$. Then
  $$\prod_{\substack{1\le l\le\ell\\\kappa_l\equiv i_1\pmod e}}
  (\yo_1-[\kappa_l-\i_1])\fo=0.$$
\end{Lemma}

\begin{proof} By \autoref{D:HeckeAlgebras}, $\prod_{l=1}^{\ell}(L_1-[\kappa_l])=0$ so that
 $\prod_{l=1}^{\ell}(L_1-[\kappa_l])\fo=0$, for all $\bi\in I$. If
 $\kappa_l\not\equiv i_1\pmod{e}$ then $[\i_1]\ne [\kappa_l]$ so that
 $(L_1-[\kappa_l])$ acts invertibly on~$\fo\H$ by
\autoref{P:Invertible}. Consequently, by \autoref{D:KLRLift},
$$0=\prod_{\substack{1\le l\le\ell\\\kappa_l\equiv i_1\pmod e}}
           (t^{\i_1}\yo_1+[\i_1]-[\kappa_l])\fo
    =t^{\i_1\<\Lambda,\alpha_{i_1}\>}\!\!\!\!
        \prod_{\substack{1\le l\le\ell\\\kappa_l\equiv i_1\pmod e}}
           (\yo_1-[\kappa_l-\i_1])\fo. $$
  As $t$ is invertible in $\O$, the lemma follows.
\end{proof}

Suppose that $\s$ is a standard tableau, $\bi=\res(\s)\in I^n$ and $1\le r<n$.  Define
\begin{equation}\label{E:beta}
\beta_r(\s)= \begin{cases}
  \dfrac{t^{\i_{r}-c_{r}(\s)}\alpha_r(\s)}{[1-\rho_r(\s)]},
      &\text{if }i_r=i_{r+1},\\[1pt]
     t^{c_{r+1}(\s)-\i_r}\alpha_r(\s)[\rho_r(\s)],
         &\text{if }i_r=i_{r+1}+1,\\
     \dfrac{t^{-\rho_r(\s)}\alpha_r(\s)[\rho_r(\s)]}{[1-\rho_r(\s)]},
      &\text{otherwise},
\end{cases}
\end{equation}
and \begin{equation}\label{E:beta2}
\widehat{\beta}_r(\s)= \begin{cases}
  \dfrac{t^{\i_{r+1}-c_{r+1}(\s)}\alpha_r(\s)}{[1+\rho_{r}(\s)]},
      &\text{if }i_r=i_{r+1},\\[1pt]
     -t^{c_{r+1}(\s)-\i_{r+1}}\alpha_r(\s)[\rho_r(\s)],
         &\text{if }i_r=i_{r+1}-1,\\
     -\dfrac{\alpha_r(\s)[\rho_r(\s)]}{[1+\rho_r(\s)]},
      &\text{otherwise}.
\end{cases}
\end{equation}
These scalars describe the action of $\psio_r$ upon the seminormal
basis.

\begin{Lemma}\label{L:PsiExpansion}
  Suppose that $1\le r<n$ and that $(\s,\t)\in\SStd(\Parts)$. Set
  $\bi=\res(\s)$, $\bj=\res(\t)$, $\u=\s(r,r+1)$ and $\v=\t(r,r+1)$. Then
  \begin{align*}
  \psio_rf_{\s\t} &= \beta_r(\s)f_{\u\t}
          -\delta_{i_ri_{r+1}}\frac{t^{\i_{r+1}-c_{r+1}(\s)}}{[\rho_r(\s)]}f_{\s\t},\\
  \intertext{and}
  f_{\s\t}\psio_r &=\widehat{\beta}_r(\t)f_{\s\v}
          -\delta_{j_rj_{r+1}}\frac{t^{\j_{r+1}-c_{r+1}(\t)}}{[\rho_r(\t)]}f_{\s\t}.
  \end{align*}
  Similarly, $\yo_rf_{\s\t}=[c_r(\s)-\i_r]f_{\s\t}$,
  and $f_{\s\t}\yo_r=[c_r(\t)-\j_r]f_{\s\t}$,
  for $1\le r\le n$.
\end{Lemma}

\begin{proof}Applying \autoref{D:KLRLift} and \autoref{E:SeminormalForm},
  $$\yo_r f_{\s\t}=t^{-\i_r}([c_r(\s)]-[\i_r])f_{\s\t}=[c_r(\s)-\i_r]f_{\s\t}.$$
  The proof that $f_{\s\t}\yo_r=[c_r(\t)-\j_r]f_{\s\t}$ is similar. We now consider~$\psio_r$.

  By \autoref{E:SeminormalGamma}, if $\mathbf{k}\in I^n$ then
  $f_{\mathbf{k}}^{\O}f_{\s\t}=\delta_{\bi\mathbf{k}}f_{\s\t}$. We use this observation below
  without mention. By \autoref{L:LTaction},
  $(T_rL_r-L_rT_r)f_{\s\t}=\alpha_r(\s)t^{c_{r+1}(\s)}[\rho_r(\s)]f_{\u\t}$.
  Hence, $\psio_r f_{\s\t}=\beta_r(\s)f_{\u\t}$ when $i_r\ne i_{r+1}$ by \autoref{D:KLRLift}
  and \autoref{E:Mr}. Now suppose that $i_r=i_{r+1}$. Then,
  using \autoref{E:Mr} and \autoref{E:SeminormalForm},
  \begin{align*}
  \psio_rf_{\s\t} &=(1+T_r)\frac{t^{\i_r}}{M_r}f_{\s\t}
     =\frac{t^{\i_r-c_r(\s)}}{[1-\rho_r(\s)]}\Big(\alpha_r(\s)f_{\u\t}
     +\big(1-\frac1{[\rho_r(\s)]}\big)f_{\s\t}\Big)\\
     &=\beta_r(\s)f_{\u\t}-\frac{t^{\i_{r+1}-c_{r+1}(\s)}}{[\rho_r(\s)]}f_{\s\t},
\end{align*}
  as required.  The formula for $f_{\s\t}\psio_r$ is proved similarly using
  \autoref{L:RightPsi} in place of \autoref{D:KLRLift}.
\end{proof}

Note that, in general, $\psio_r f_{\s\t}\ne(f_{\t\s}\psio_r)^{\bigast}$.

The next relation can also be proved using \autoref{L:MrCommute} and
\autoref{L:RightPsi}.

\begin{Corollary}
  Suppose that  $|r-t|>1$, for $1\le r,t<n$. Then
   $\psio_r\psio_t=\psio_t\psio_r$.
\end{Corollary}

\begin{proof}
  It follows easily from \autoref{L:PsiExpansion} that
  $\psi_r\psi_t f_{\s\t}=\psi_t\psi_r f_{\s\t}$, for
  all $(\s,\t)\in\SStd(\Parts)$.  Hence, by \autoref{L:MurphyIdempotent},
  $\psio_r\psio_t\fo=\psio_t\psio_r\fo$, for all $\bi\in I^n$.
\end{proof}

\begin{Lemma}\label{L:PsiYRelation}
  Suppose that $1\le r<n$ and $\bi\in I^n$. Then
  $$\psio_r \yo_{r+1}\fo=(\yo_r\psio_r+\delta_{i_ri_{r+1}})\fo
  \quad\text{and}\quad
    \yo_{r+1}\psio_r\fo=(\psio_r\yo_r+\delta_{i_ri_{r+1}})\fo.$$
\end{Lemma}

\begin{proof}Both formulas can be proved similarly, so we consider only the
  first one. We prove the stronger result that $\psio_r
  \yo_{r+1}f_{\s\t}=(\yo_r\psio_r+\delta_{i_ri_{r+1}})f_{\s\t}$, whenever
  $(\s,\t)\in\SStd(\Parts)$ and $\res(\s)=\bi$.  By
  \autoref{E:ResidueIdempotent} this implies the lemma.

  Suppose first that $i_r=i_{r+1}$. Then, using
  \autoref{L:PsiExpansion},
  \begin{align*}
    \psio_r\yo_{r+1}f_{\s\t}
    &=[c_{r+1}(\s)-\i_{r+1}]\Big(\beta_r(\s)f_{\u\t}
       -\frac{t^{\i_{r+1}-c_{r+1}(\s)}}{[\rho_r(\s)]}f_{\s\t}\Big).
  \end{align*}
  On the other hand, by \autoref{L:PsiExpansion} and \autoref{E:beta},
  \begin{align*}\qquad
  (\yo_r\psio_r+1)f_{\s\t}
    &= [c_r(\u)-\i_{r+1}]\beta_r(\s)f_{\u\t}
    +\Big(1-\frac{t^{\i_r-c_{r+1}(\s)}[c_r(\s)-\i_r]}{[\rho_r(\s)]}\Big)f_{\s\t}\\
    &= [c_r(\u)-\i_{r+1}]\beta_r(\s)f_{\u\t}
       +\frac{[\i_{r+1}-c_{r+1}(\s)]}{[\rho_r(\s)]}f_{\s\t}.
  \end{align*}
  Therefore, $\psio_r\yo_{r+1}f_{\s\t}=(\yo_r\psio_r+1)f_{\s\t}$ since
  $c_r(\u)=c_{r+1}(\s)$ and $i_r=i_{r+1}$.

  If $i_r\ne i_{r+1}$ then the calculation is easier because
  $$\psio_r\yo_{r+1}f_{\s\t} =[c_{r+1}(\s)-\i_{r+1}]\beta_r(\s)f_{\u\t}
               =\yo_r\psio_rf_{\s\t},$$
  where, for the last equality, we again use the fact that $c_r(\u)=c_{r+1}(\s)$.
\end{proof}

%
%
%

The following simple combinatorial identity largely determines both the quadratic
and the (deformed) braid relations for the $\psio_r$, for $1\le r<n$. This
result can be viewed as a graded analogue of the defining property
\autoref{D:SNCS}(c) of a seminormal coefficient system.

\begin{Lemma}\label{L:BetaSquared}
  Suppose that $1\le r<n$ and $\s,\u\in\Std(\blam)$ with $\u=\s(r,r+1)$ and
  $\res(\s)=\bi\in I^n$, for $\blam\in\Parts$. Then
  \begin{align*}
  \beta_r(\s)\beta_r(\u)&=\begin{cases}
    t^{c_r(\s)+c_{r+1}(\s)-\i_r-\i_{r+1}}[1-\rho_r(\s)][1+\rho_r(\s)],
            &\text{if }i_r\leftrightarrows i_{r+1},\\
    t^{c_{r+1}(\s)-\i_{r+1}}[1+\rho_r(\s)],&\text{if }i_r\rightarrow i_{r+1},\\
    t^{c_r(\s)-\i_r}[1-\rho_r(\s)],&\text{if }i_r\leftarrow i_{r+1},\\
    -\frac{t^{2\i_r-2c_{r+1}(\s)}}{[\rho_r(\s)]^2},&\text{if }i_r=i_{r+1},\\
      1,&\text{otherwise}.
    \end{cases}
  \end{align*}
\end{Lemma}

\begin{proof}
  The lemma follows directly from the definition of $\beta_r(\s)$
  using~\autoref{D:SNCS}(c).
\end{proof}


It is time to pay the price for the failure of the embedding
$I\hookrightarrow\Z$ to extend to an embedding of quivers. Together with the
cyclotomic relation, this is place where the KLR grading fails to lift to the
algebra $\HO$.  Recall from \autoref{D:KLRLift} that $\yo_r
\fo=t^{-\i_r}(L_r-[\i_r])\fo$, where $1\le r\le n$ and $\bi\in I^n$.
For $d\in\Z$ define
\begin{equation}\label{E:dotY}
  \dyo[d]_r\fo=t^{d-\i_r}(L_r-[\i_r-d])\fo
  =(t^d\yo_r+[d])\fo.
\end{equation}
In particular, $\dyo[0]_r=\yo_r$ and $\dyo[d]_r\otimes_\O1_K=\yo_r\otimes_\O1_K$
whenever $e$ divides $d\in\Z$,

As a final piece of notation, set $\rho_r(\bi)=\i_r-\i_{r+1}\in\Z$, for $\bi\in I^n$
and $1\le r<n$.

\begin{Proposition}\label{P:PsiSquare}
  Suppose that $1\le r<n$ and  $\bi\in I^n$. Then
  $$(\psio_r)^2\fo = \begin{cases}
    (\dyo[\eps]_r-\yo_{r+1})(\dyo[\spe]_{r+1}-\yo_r)\fo,& \text{if }i_r\leftrightarrows i_{r+1},\\
    (\dyo[\eps]_r-\yo_{r+1})\fo,& \text{if }i_r\rightarrow i_{r+1},\\
    (\dyo[\spe]_{r+1}-\yo_r)\fo,& \text{if }i_r\leftarrow i_{r+1},\\
       0,&\text{if }i_r=i_{r+1},\\
       \fo,&\text{otherwise.}
       \end{cases}$$
\end{Proposition}

\begin{proof}Once again, by \autoref{E:ResidueIdempotent} it is enough to prove the
  corresponding formulas for $(\psio_r)^2 f_{\s\t}$, where
  $(\s,\t)\in\SStd(\Parts)$ and $\bi=\res(\s)$.

  Suppose that $i_r=i_{r+1}$. Let $\u=\s(r,r+1)$ and $\bj=\res(\u)$. By
  \autoref{L:PsiExpansion},
     $$ (\psio_r)^2f_{\s\t}=
       \Big(\frac{t^{2\i_r-2c_{r+1}(\s)}}{[\rho_r(\s)]^2}+\beta_r(\s)\beta_r(\u)\Big)f_{\s\t}
       -\Big(\frac{\beta_r(\s)t^{\i_r-c_r(\s)}}{[\rho_r(\u)]}
       +\frac{\beta_r(\s)t^{\j_r-c_r(\u)}}{[\rho_r(\s)]}\Big)f_{\u\t}.$$
       Note that $\rho_r(\s)=-\rho_r(\u)$ and $i_r=j_r$, so that
       $t^{\j_r-c_r(\u)}[\rho_r(\u)]=-t^{\i_r-c_r(\s)}[\rho_r(\s)]$. Hence, using
       \autoref{L:BetaSquared}, $(\psio_r)^2f_{\s\t}=0$ when $i_r=i_{r+1}$ as
       claimed.

       Now suppose that $i_r\ne i_{r+1}$. Then, by
       \autoref{L:PsiExpansion} and \autoref{L:BetaSquared},
       \begin{align*}
       (\psio_r)^2f_{\s\t}&= \beta_r(\s)\beta_r(\u)f_{\s\t}\\
          &=\begin{cases}
             t^{c_r(\s)+c_{r+1}(\s)-\i_r-\i_{r+1}}[1-\rho_r(\s)][1+\rho_r(\s)]f_{\s\t},
                      &\text{if }i_r\leftrightarrows i_{r+1},\\
             t^{c_{r+1}(\s)-\i_{r+1}}[1+\rho_r(\s)]f_{\s\t},&\text{if }i_r\rightarrow i_{r+1},\\
             t^{c_{r}(\s)-\i_{r}}[1-\rho_r(\s)]f_{\s\t},&\text{if }i_r\leftarrow i_{r+1},\\
             f_{\s\t},&\text{otherwise}.
      \end{cases}
     \end{align*}
     As in \autoref{L:PsiExpansion}, if $d\in\Z$ then
     $\dyo[d]_rf_{\s\t}=[c_r(\s)-\i_r+d]f_{\s\t}$. So, if $i_r\rightarrow i_{r+1}$ then
     \begin{align*}
      (\dyo[\eps]_r-\yo_{r+1})f_{\s\t}
         &=\([c_r(\s)+1-\i_{r+1}]-[c_{r+1}(\s)-\i_{r+1}]\)f_{\s\t}\\
         &=t^{c_{r+1}(\s)-\i_{r+1}}[1+\rho_r(\s)]f_{\s\t}
          =(\psio_r)^2f_{\s\t}.
   \end{align*}
   The cases when $i_r\rightarrow i_{r+1}$ and $i_r\leftrightarrows i_{r+1}$ are similar.
\end{proof}

Set $\OBraid=\psio_r\psio_{r+1}\psio_r-\psio_{r+1}\psio_r\psio_{r+1}$, for
$1\le r<n-1$.

\begin{Proposition}\label{P:PsiBraid}
  Suppose that $1\le r<n$ and $\s,\t\in\Std\(\blam)$, with $\s\in\Std(\bi)$ for
  $\bi\in I^n$.  Then
  $$\OBraid f_{\s\t}=
  \begin{cases}
    (\dyo[\eps]_r+\dyo[\eps]_{r+2}-
    \dyo[\eps]_{r+1}-\dyo[\spe]_{r+1})f_{\s\t},
            &\text{if }i_{r+2}=i_r\rightleftarrows i_{r+1},\\
   -t^{\eps}f_{\s\t}, &\text{if }i_{r+2}=i_r\rightarrow i_{r+1},\\
   f_{\s\t}, &\text{if }i_{r+2}=i_r\leftarrow i_{r+1},\\
    0,&\text{otherwise}.
  \end{cases}$$
\end{Proposition}

\begin{proof}
  We mimic the proof of the  braid relations from \autoref{L:SNCExistence}.

  Define (not necessarily standard) tableaux $\u_1=\s(r,r+1)$,
  $\u_2=\s(r+1,r+2)$, $\u_{12}=\u_1(r+1,r+2)$, $\u_{21}=\u_2(r,r+1)$ and
  $\u_{121}=\u_{12}(r,r+1)=\u_{212}$. To ease notation set $i=i_r$, $j=i_{r+1}$
  and $k=i_{r+2}$. The relationship between these tableaux, and their residues
  $\set{\res_s(\u)|r\le s\le r+2}=\{i,j,k\}$, is illustrated in the following diagram.
  
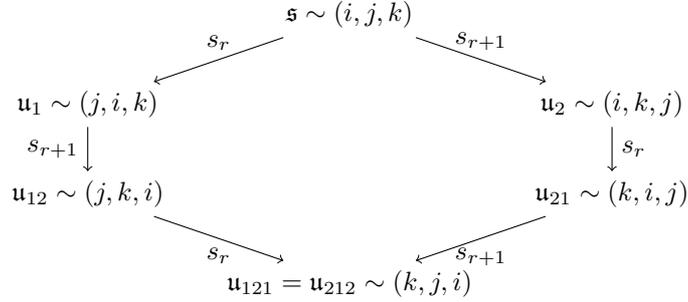
\begin{center}
    \begin{tikzpicture}
     \matrix[matrix of math nodes,row sep=6mm,column sep=6mm]{
     &|(s)|\s\sim(i,j,k)\\
     |(t)|\u_1\sim(j,i,k) &       &|(u)|\u_2\sim(i,k,j)\\
     |(v)|\u_{12}\sim(j,k,i) &       &|(w)|\u_{21}\sim(k,i,j)\\
     &|(x)|\u_{121}=\u_{212}\sim(k,j,i)\\
     };
     \draw[->](s)--node[above]{$s_r$}(t);
     \draw[->](s)--node[above]{$s_{r+1}$}(u);
     \draw[->](t)--node[left]{$s_{r+1}$}(v);
     \draw[->](u)--node[right]{$s_r$}(w);
     \draw[->](v)--node[below]{$s_r$}(x);
     \draw[->](w)--node[below]{$s_{r+1}$}(x);
  \end{tikzpicture}
  \end{center}
  Note that if any tableau $\u\in\{\u_1,\u_2,\u_{12},\u_{21},\u_{121}\}$ is not
  standard then, by definition, $f_{\u\t}=0$ so this term can be ignored in all
  of the calculations below.

  We need to compute $\OBraid f_{\s\t}$. To start with, observe that by
  \autoref{L:PsiExpansion}, the coefficient of $f_{\u_{121}\t}$ in $\OBraid
  f_{\s\t}$ is equal to
  $$\beta_r(\s)\beta_{r+1}(\u_{1})\beta_r(\u_{12})
                -\beta_{r+1}(\s)\beta_r(\u_{2})\beta_{r+1}(\u_{21}).$$
  By definition, the scalars $[\rho_r(\s)]$ and $[1-\rho_r(\s)]$ are determined by
  the positions of~$r$ and~$r+1$ in~$\s$, so it is easy to see that
  \begin{equation}\label{E:identities}
    \begin{aligned}
      \rho_r(\s)&=\rho_{r+1}(\u_{21}),&
      \rho_r(\u_1)&=\rho_{r+1}(\u_{121}),&
      \rho_r(\u_{2})&=\rho_{r+1}(\u_{1}),\\
      \rho_r(\u_{12})&=\rho_{r+1}(\s),&
      \rho_r(\u_{21})&=\rho_{r+1}(\u_{12}),&
      \rho_r(\u_{121})&=\rho_{r+1}(\u_{2}).
    \end{aligned}
  \end{equation}
  Observe that
  $\alpha_r(\s)\alpha_{r+1}(\u_1)\alpha_r(\u_{12})
        =\alpha_{r+1}(\s)\alpha_r(\u_2)\alpha_{r+1}(u_{21})$ by \autoref{D:SNCS}(a).
  Keeping track of the exponent of~$t$, \autoref{E:beta} and
  \autoref{E:identities} now imply that
  $\beta_r(\s)\beta_{r+1}(\u_1)\beta_r(\u_{12})
        =\beta_{r+1}(\s)\beta_r(\u_2)\beta_{r+1}(\u_{21})$. Note
  that~\autoref{D:SNCS}(a) is crucial here.  Therefore, the coefficient of
  $f_{\u_{121}\t}$ in $\OBraid f_{\s\t}$ is zero for any choice of $i,j$ and~$k$.
  As the coefficient of $f_{\u_{121}\t}$ in~$\OBraid f_{\s\t}$ is always zero
  we will omit $f_{\u_{121}\t}$ from most of the calculations that follow.

  There are five cases to consider.

  \Case{$i$, $j$ and $k$ are pairwise distinct}
  By \autoref{L:PsiExpansion} and the last paragraph,
  $$\OBraid f_{\s\t} =\(\beta_r(\s)\beta_{r+1}(\u_1)\beta_r(\u_{12})
  -\beta_{r+1}(\s)\beta_r(\u_2)\beta_{r+1}(\u_{21})\)f_{\u_{121}\t}=0,$$
  as required by the statement of the proposition.

  \Case{$i=j\ne k$}
  In this case, using \autoref{L:PsiExpansion},
  $$ \OBraid f_{\s\t} =\beta_{r+1}(\s)\beta_r(\u_2)\Big(-\frac{t^{\i-c_{r+1}(\s)}}{[\rho_r(\s)]}
    +\frac{t^{\i-c_{r+2}(\u_{21})}}{[\rho_{r+1}(\u_{21})]}\Big)
    f_{\u_{21}\t}.
  $$
  Now $\rho_r(\s)=\rho_{r+1}(\u_{21})$ and $c_{r+1}(\s)=c_{r+2}(\u_{21})$, as in
  \autoref{E:identities}. Hence, $\OBraid f_{\s\t}=0$ when $i=j\ne k$.

  \Case{$i\ne j=k$}
  This is almost identical to Case~2, so we leave the details to the reader.

  \Case{$i=k\ne j$}
  Typographically, it is convenient to set $c=c_r(\s)$, $c'=c_{r+1}(\s)$ and
  $c''=c_{r+2}(\s)$. According to the statement of the proposition, this is the
  only case where $\OBraid f_{\s\t}\ne0$. Using \autoref{L:PsiExpansion}, we see
  that
  \begin{align*}
    \OBraid f_{\s\t} &=\Big(-\beta_r(\s)\frac{t^{\i-c_{r+2}(\u_1)}}{[\rho_{r+1}(\u_1)]}\beta_r(\u_1)
    +\beta_{r+1}(\s)\frac{t^{\i-c_{r+1}(\u_2)}}{[\rho_r(\u_2)]}\beta_{r+1}(\u_2)\Big)f_{\s\t}.\\
    &=\frac{t^{\i-c''}}{[c-c'']}
          \(-\beta_r(\s)\beta_r(\u_1)+\beta_{r+1}(\s)\beta_{r+1}(\u_2)\)f_{\s\t}.
  \end{align*}
  Expanding the last equation using \autoref{L:BetaSquared} shows that
  $$
    \OBraid f_{\s\t} =\begin{cases}
    -\dfrac{t^{c}[1-\rho_r(\s)][1+\rho_r(\s)]
              -t^{c''}[1-\rho_{r+1}(\s)][1+\rho_{r+1}(\s)]}%
              {t^{c''-c'+\j}[c-c'']}f_{\s\t},
                  &\text{if }i\leftrightarrows j,\\
    -\dfrac{[1+\rho_r(\s)]-[1-\rho_{r+1}(\s)]}%
       {t^{c''-\i-c'+\j}[c-c'']}f_{\s\t},
                  &\text{if }i\rightarrow j,\\
    -\dfrac{t^{c}[1-\rho_r(\s)]-t^{c''}[1+\rho_{r+1}(\s)]}
              {t^{c''}[c-c'']}f_{\s\t},
                  &\text{if }i\leftarrow j,\\
             0,&\text{otherwise.}
           \end{cases}
   $$
  (Note that, by assumption, the case $i=j$ does not arise.) If
  $i\leftrightarrows j$ then a straightforward calculation shows that in this
  case
  \begin{align*}
    \OBraid f_{\s\t}&=
   -\Big([c'-\j+2]+[c'-\j]-[c+1-\j]-[c''+1-\j]\Big)f_{\s\t}\\
   &=-\(\dyo[\eps]_{r+1}+\dyo[\spe]_{r+1}-\dyo[\eps]_r-\dyo[\eps]_{r+2}\)f_{\s\t},
 \end{align*}
 where the last equality uses \autoref{L:PsiYRelation} and the
 observation that, because $e=2$, we have $\{1\pm\rho_r(\bi)\}=\{0,2\}$ and
 $\{\i,\j\}=\{0,1\}$.  A similar, but easier, calculation shows that if
 $i\rightarrow j$ then $\OBraid
 f_{\s\t}=-t^{1+\i-\j}f_{\s\t}=-t^{1+\rho_r(\bi)}f_{\s\t}$ and if $ i\leftarrow j$
 then $\OBraid f_{\s\t}=f_{\s\t}$. If $i\ne j$ and $i\noedge j$ then we have
 already seen that  $\OBraid f_{\s\t}=0$, so this completes the proof of Case~4.

  \Case{$i=j=k$}
  We continue to use the notation for $c,c',c''$ from Case~4.  By
  \autoref{L:PsiExpansion} (compare with  the proof of
  \autoref{L:SNCExistence}), $\OBraid f_{\s\t}$ is equal to
  $$\begin{array}{r}
    -\Big(\frac{t^{3\i-2c'-c''}}{[\rho_r(\s)]^2[\rho_{r+1}(\s)]}
    -\frac{t^{3\i-c'-2c''}}{[\rho_{r+1}(\s)]^2[\rho_r(\s)]}
    +\frac{t^{\i-c''}\beta_r(\s)\beta_r(\u_1)}{[\rho_{r+1}(\u_1)]}
    -\frac{t^{\i-c''}\beta_{r+1}(\s)\beta_{r+1}(\u_2)}{[\rho_r(\u_2)]}\Big)f_{\s\t}\\
    +t^{2\i}\beta_r(\s)\Big(\frac{t^{-c''-c}}{[\rho_{r+1}(\u_1)][\rho_r(\u_1)]}
        +\frac{t^{-c'-c''}}{[\rho_r(\s)][\rho_{r+1}(\s)]}
        -\frac{t^{-2c''}}{[\rho_{r+1}(\s)][\rho_{r+1}(\u_1)]}\Big)f_{\u_1\t}\\
    +t^{2\i}\beta_{r+1}(\s)\Big(\frac{t^{-c'-c''}}{[\rho_r(\s)][\rho_r(\u_2)]}
        -\frac{t^{-c''-c'}}{[\rho_r(\u_2)][\rho_{r+1}(\u_2)]}
        -\frac{t^{-c''-c'}}{[\rho_{r+1}(\s)][\rho_r(\s)]} \Big)f_{\u_2\t}\\
    -t^{\i-c''}\beta_r(\s)\beta_{r+1}(\u_1)
    \Big(\frac1{[\rho_r(\u_{12})]}-\frac1{[\rho_{r+1}(\s)]}\Big)f_{\u_{12}\t}\\
    -t^{\i-c'}\beta_{r+1}(\s)\beta_r(\u_2)
           \Big(\frac1{[\rho_r(\s)]}-\frac1{[\rho_{r+1}(\u_{21})]}\Big)f_{\u_{21}\t}.
         \end{array}$$
  Using \autoref{E:identities} it is easy to see that the coefficients of
  $f_{\u_{12}\t}$ and$f_{\u_{21}\t}$ are both zero. On the other hand, if
  $t\ne1$ then the coefficient of $t^{2\i}\beta_r(\s)f_{\u_1\t}$ in $\OBraid f_{\s\t}$ is
  $$ \frac{t-1}{(t^{c'}-t^c)(t^c-t^{c''})}
            +\frac{t-1}{(t^{c}-t^{c'})(t^{c'}-t^{c''})}
            -\frac{t-1}{(t^{c'}-t^{c''})(t^c-t^{c''})}=0.
  $$
  The case when $t=1$ now follows by specialisation.  Similarly,
  the coefficient of~$f_{\u_2\t}$ in $\OBraid f_{\s\t}$ is also zero. Finally, using
  \autoref{L:BetaSquared} and \autoref{E:identities}, the coefficient of
  $f_{\s\t}$ in $\OBraid f_{\s\t}$ is zero as the four terms above, which give the coefficient
  of $f_{\s\t}$ in the displayed equation, cancel out in pairs. Hence,
  $\OBraid f_{\s\t}=0$ when $i=j=k$, as required.

  This completes the proof.
\end{proof}

%
%

\subsection{A deformation of the quiver Hecke algebra}\label{S:KLRDeformation}
Using the results of the last two sections we now describe $\H(\O)$ by
generators and relations using the `$\O$-KLR generators' of~$\HO$.

Let $\RO$ be the abstract algebra defined by the generators and relations in the
statement of \autoref{Thm:KLRDeformation}. We abuse notation and use the same
symbols for the generators of~$\RO$ and the corresponding elements in~$\HO$ that
we defined in \autoref{S:KLRGens}. The previous section shows that there is a
surjection $\RO\surjection\HO$. We want to prove that this map is an
isomorphism.

The next lemma, which is modelled on \cite[Lemma~2.1]{BK:GradedKL}, will be used
to show that $\RO$ is finitely generated as an $\O$-module.

\begin{Lemma}\label{L:yONilpotence}
  Suppose that $1\le r\le n$ and $\bi\in I^n$. Then there exists a
  \textit{multiset} $X_r(\bi)\subseteq e\Z$ such that
  $$\prod_{c\in X_r(\bi)}(\yo_r-[c])\fo=0.$$

\end{Lemma}

\begin{proof}
  We argue by induction on~$r$. If $r=1$ then the relations in~$\RO$ say that we
  can take $X_1(\bi)$ to be the multiset with elements $\kappa_l-\i_1$, where
  $\kappa_l\equiv i_1\pmod e$ for $1\le l\le\ell$. By induction we assume that
  we have proved the result for $\yo_r$ and use this to prove the result for
  $\yo_{r+1}$. There are three cases to consider. Throughout, let
  $\bj=\s_r\cdot\bi$.

  \Case{$i_{r+1}\noedge i_r$}
  Set $X_{r+1}(\bi)=X_r(\bj)$. Then,
  using the relations in~$\RO$,
  \begin{align*}
    \prod_{c\in X_{r+1}(\bi)}(\yo_{r+1}-[c])\fo
       &=\prod_{c\in X_{r+1}(\bi)}(\yo_{r+1}-[c])(\psio_r)^2\fo\\
       &=\psio_r\prod_{c\in X_r(\bj)}(\yo_r-[c])\fo[\bj]\psio_r
        =0,
  \end{align*}
  where the last equality follows by induction.

  \Case{$i_r\rightarrow i_{r+1}$ or $i_r\leftarrow i_{r+1}$}
  We consider only the case when $i_r\rightarrow i_{r+1}$. The case when
  $i_r\leftarrow i_{r+1}$ is similar. Using
  quadratic relation for $\psio_r$ in~$\RO$,
  $$(\yo_{r+1}-[c])\fo
  =t^{1+\rho_i(\bi)}(\yo_r-[c-1-\rho_i(\bi)])\fo-(\psio_r)^2\fo.$$
  Let $X_{r+1}(\bi)$ be the disjoint union
  $X_r(\bj)\sqcup X^+_{r+1}(\bi)$, where
  $X^+_{r+1}(\bi)$ is the multiset $\set{c+1+\rho_r(\bi)|c\in X_r(\bi)}$. If
  $d=c+1+\rho_r(\bi)\in X^+_{r+1}(\bi)$ then, by the last displayed equation,
  \begin{align*}
    \prod_{c'\in X_{r+1}(\bi)}(\yo_{r+1}-[c'])\fo
      &=\prod_{c'\in X_{r+1}(\bi)\setminus\{d\}}(\yo_{r+1}-[c'])\cdot
      \(t^{1+\rho_r(\bi)}(\yo_r-[c]) - (\psio_r)^2\)\fo\\
          &=t^{1+\rho_r(\bi)}(\yo_r-[c])\prod_{c'\in X_{r+1}(\bi)\setminus\{d\}}
               (\yo_{r+1}-[c'])\fo\\
          &\qquad -\psio_r\prod_{c'\in X_{r+1}(\bi)\setminus\{d\}}
               (\yo_r-[c'])\fo[\bj]\psio_r.
  \end{align*}
  The second summand is zero by induction because $X_r(\bj)$ is contained in
  $X_{r+1}(\bi)$. Therefore, arguing this way for every
  $d\in X^+_{r+1}(\bi)$, there exists $N\in\Z$ such that
  $$\prod_{c'\in X_{r+1}(\bi)}(\yo_{r+1}-[c'])\fo
  =t^{Ne} \prod_{c'\in X_r(\bj)}(\yo_{r+1}-[c'])
       \cdot\prod_{c\in X_r(\bi)}(\yo_r-[c])\fo=0,
  $$
  where we again use induction for the last equality. This completes the proof of
  the inductive step when $i_r\rightarrow i_{r+1}$.

  \Case{$i_r\leftrightarrows i_{r+1}$}
  This is similar to Case~2 but slightly more involved. Define
  $X_{r+1}(\bi)=X_r(\bj)\sqcup X^+_{r+1}(\bi)\sqcup X^-_{r+1}(\bi)$, where
  $X^\pm_{r+1}=\set{c\pm 1+\rho_r(\bi)|c\in X_r(\bi)}$. If $c\in X_r(\bi)$ set
  $$
  c^{+}:=c+1+\rho_r(\bi),\quad\, c^{-}:=c-1+\rho_r(\bi).
  $$
  Using the equality
  $$
  (\psio_r)^2\fo=(t^{1+\rho_r(\bi)}\yo_r+[1+\rho_r(\bi)]
         -\yo_{r+1})(t^{1-\rho_r(\bi)}\yo_{r+1}+[1-\rho_r(\bi)]-\yo_r)\fo,
  $$
  we see that
  $$
  (\psio_r)^2\fo=-t^{1-\rho_r(\bi)}(\yo_{r+1}-[c^{+}])(\yo_{r+1}-[c^{-}])\fo
                 +(\yo_r-[c])\fo F_r(y,c),
  $$
  where $F_r(y,c)$ is a polynomial in $\yo_r$ and $\yo_{r+1}$ with coefficients
  in $\Z[t,t^{-1}]$.  Hence,
  \begin{align*}
    &\quad\,\prod_{c'\in X_{r+1}(\bi)}(\yo_{r+1}-[c'])\fo\\
   &=\prod_{c'\in X_{r+1}(\bi)\setminus\{c^{+},c^{-}\}}(\yo_{r+1}-[c'])\cdot
      t^{\rho_r(\bi)-1}\( (\yo_r-[c])\fo F_r(y,c) - (\psio_r)^2\fo\)\\
          &=t^{\rho_r(\bi)-1}(\yo_r-[c])\prod_{c'\in X_{r+1}(\bi)\setminus\{c^{+},c^{-}\}}
               (\yo_{r+1}-[c'])\fo F_r(y,c)\\
          &\qquad -t^{\rho_r(\bi)-1}\psio_r\prod_{c'\in X_{r+1}(\bi)\setminus\{c^{+},c^{-}\}}
               (\yo_r-[c'])\fo[\bj]\psio_r.
  \end{align*}
 The second summand is zero by induction because $X_r(\bj)$ is contained in
  $X_{r+1}(\bi)$. Therefore, arguing this way for every
  $c\in X_{r}(\bi)$, there exists $N\in\Z$ such that
  $$\prod_{c'\in X_{r+1}(\bi)}(\yo_{r+1}-[c'])\fo
  =t^{Ne} \prod_{c'\in X_r(\bj)}(\yo_{r+1}-[c'])
       \cdot\prod_{c\in X_r(\bi)}(\yo_r-[c])\fo F_r(y)=0,
  $$
  where $F_r(y)$ is a polynomial in $\yo_1,\cdots,\yo_{n}$ with
  coefficients in $\Z[t,t^{-1}]$ and we again use induction for the last
  equality. This completes the proof of the inductive step when
  $i_r\leftrightarrows i_{r+1}$.

  \Case{$i_{r+1}=i_r$}
  Let $\phi_r=\psio_r(\yo_r-\yo_{r+1})\fo$.  Then $\phi_r\psio_r=-2\psio_r\fo$,
  so that $(1+\phi_r)^2\fo=\fo$. Moreover, an easy albeit uninspiring
  calculation reveals that
  $$
    (1+\phi_r)\yo_r(1+\phi_r)\fo
     =(\yo_r+\phi_r\yo_r+\yo_r\phi_r+\phi_r\yo_r\phi_r)\fo\\
     =\yo_{r+1}\fo.
  $$
  Therefore, setting $X_{r+1}(\bi)=X_r(\bi)$,
  $$\prod_{c\in X_{r+1}(\bi)}(\yo_{r+1}-[c])\fo
     =(1+\phi_r)\prod_{c\in X_r(\bi)}(\yo_r-[c])\fo(1+\phi_r)=0,$$
  where the last equality follows by induction. This completes the proof.
\end{proof}

Suppose that $(\O,t)$ is an idempotent subring of~$\K$. So far we have not used
the assumption that $[de]\in J(\O)$, for $d\in\Z$. This comes into play in the
next theorem, which is \autoref{Thm:KLRDeformation} from the introduction.

\begin{Theorem}\label{T:IntegralKLR}
  Suppose that $(\O,t)$ is an $e$-idempotent subring of $\K$.
  Then $\RO\cong\HO$ as $\O$-algebras.
\end{Theorem}

\begin{proof}By the results in the last two sections, the elements given in
  \autoref{D:KLRLift} satisfy all of the relations of the corresponding
  generators of $\RO$. Hence, by \autoref{P:OGeneration}, there is a surjective
  $\O$-algebra homomorphism $\theta{:}\RO\surjection\HO$, which maps the
  generators of $\RO$ to the corresponding elements of $\HO$.

  If $w\in\Sym_n$ then set $\psio_w=\psio_{r_1}\dots\psio_{r_k}$, where
  $w=s_{r_1}\dots s_{r_k}$ is a reduced expression for $w$. In general,
  $\psio_w$ will depend upon the choice of reduced expression, however, using
  the relations in $\RO$ it follows that every element in $\RO$ can be written
  as a linear combination of elements of the form $f(y)\psi_we(\bi)$, where
  $f(y)\in\O[\yo_1,\dots,\yo_n]$, $w\in\Sym_n$ and $\bi\in I^n$. Therefore,
  $\RO$ is finitely generated as an $\O$-module by \autoref{L:yONilpotence}.

  Now suppose that $\m$ is a maximal ideal of~$\O$ and let
  $K=\O/\m\cong\O_{\m}/\m\O_{\m}$ and $\zeta=t+\m$. Then
  $1+\zeta+\dots+\zeta^{e-1}=0$ in $K$, since $[e]\in J(\O)\subseteq\m$. Note
  also that $1+\zeta+\dots+\zeta^{k-1}\ne0$ if $k\notin e\Z$ since $\O$ is an
  $e$-idempotent subring. Consequently, $\dyo[de]_r\otimes1_K=\yo_r\otimes 1_K$,
  for all $d\in\Z$. It is easy to see that all of the shifts $1\pm\rho_r(\bi)$
  appearing in the statement of theorem are equal to either~$0$ or to~$e$.
  Therefore,upon base change to $K$ the relations of
  $R_n(\O_{\m})\otimes_{\O_{\m}} K$ coincide with the relations of the quiver
  Hecke algebra $\R(K)$, see \autoref{D:QuiverRelations} and \autoref{T:BKiso}.
  Consequently, $R_n(\O_{\m})\otimes_{\O_{\m}} K\cong\R(K)$, so that
  $\dim R_n(\O_{\m})\otimes_{\O_{\m}} K=\dim\H(K)$ by \cite[Theorem 4.20]{BK:GradedDecomp}.

  By the last paragraph, if $K=\O/\m$, for any maximal ideal $\m$ of~$\O$, then
  $\dim R_n(\O_{\m})\otimes_{\O_{\m}} K=\dim\H(K)=\ell^nn!$.  Moreover, by the
  second paragraph of the proof, $R_n(\O_{\m})$ is a finitely generated
  $\O_{\m}$-algebra.  Therefore,  Nakayama's lemma applies and it implies that
  $R_n(\O_{\m})$ is a free $\O_{\m}$-module of rank~$\ell^nn!$.  Hence, the map
  $\theta_{\m}:R_n(\O_{\m})\bijection\H(\O_{\m})$ is an isomorphism of
  $\O_{\m}$-algebras. It follows that $\theta$ is an isomorphism of
  $\O$-algebras, as required.
\end{proof}

The proof of \autoref{T:IntegralKLR} gives the following.

\begin{Corollary}\label{C:ROField}
  Suppose that $K=\O/\m$, where $\m$ is a maximal ideal of~$\O$. Then
  $$\R(K)\cong\RO\otimes_\O K\cong\H(K).$$
\end{Corollary}

\begin{Remarks}\label{R:OIsomorphism}
  (a) The proof of \autoref{T:IntegralKLR} uses
  \cite[Theorem~4.20]{BK:GradedDecomp} to bound the rank of $\RO$. The
  proof of \cite[Theorem~4.20]{BK:GradedDecomp} does not depend on
  Brundan-Kleshchev's isomorphism \autoref{T:BKiso} (\cite[Theorem
  1.1]{BK:GradedKL}).  Instead, \cite[Theorem~4.20]{BK:GradedDecomp}
  depends on the Ariki-Brundan-Kleshchev Categorification Theorem
  \cite[Theorem~4.18]{BK:GradedDecomp}. Consequently,
  \autoref{T:IntegralKLR} gives a new proof of Brundan and Kleshchev's
  Isomorphism \autoref{T:BKiso}. It should be
  possible to prove \autoref{T:IntegralKLR} directly, without appealing
  to \cite[Theorem 4.18]{BK:GradedDecomp}, by adapting the arguments of
  \cite[Theorem~3.3]{BK:GradedKL}.

  (b) In proving \autoref{T:BKiso}, Brundan and Kleshchev~\cite{BK:GradedKL}
  construct a family of isomorphisms $\R\bijection\HK$ that depend on a choice
  of polynomials $Q_r(\bi)$ that can be varied subject to certain constraints.
  In our setting this amounts to choosing certain `scalars'
  $q_r(\bi)$, which are rational functions in $L_r$ and $L_{r+1}$,  such that $q_r(\bi)\fo\in\HO$
  and then defining
  $$\psi_r^\O \fo=\begin{cases}
      (T_r+1)\frac{t^{\i_r}}{M_r}\fo,&\text{if }i_r=i_{r+1},\\
      (T_rL_r-L_rT_r)q_r(\bi)\fo,
           &\text{otherwise},\\
    \end{cases}$$
  such that the corresponding $\beta$-coefficients still satisfy the constraints
  of \autoref{L:BetaSquared} and \autoref{D:SNCS}(a). To make this more
  precise, as in \autoref{L:PsiExpansion} write
  $$\psio_r f_{\s\t}=\beta'_r(\s)f_{\u\t}
       -\delta_{i_ri_{r+1}}\frac{t^{\i_{r+1}-c_{r+1}(\s)}}{[\rho_r(\s)]},$$
  where $\u=\s(r,r+1)$, $\beta'_r(\s)\in\K$, $\s\in\Std(\bi)$, $\bi\in I^n$
  and $1\le r<n$. (Explicitly,
  $\beta'_r(\s)=t^{c_{r+1}(\s)}[\rho_r(\s)]q_r(\s)$ by \autoref{L:LTaction}, where
  $q_r(\s)\in\K$ is the scalar such that $q_r(\bi)f_{\s\t}=q_r(\s)f_{\s\t}$.)
  Then we require that the scalars $\beta'_r(\s)$ satisfy \autoref{L:BetaSquared}
  and the ``braid relation'' of~\autoref{D:SNCS}(a).  If the $q_r(\bi)$ are
  chosen so that these two identities are satisfied then it is easy to see that
  argument used to prove
  \autoref{T:IntegralKLR} applies, virtually without change, using these more
  general elements. The key point is that \autoref{L:BetaSquared} still holds.
  The corresponding identities in Brundan and Kleshchev's work are~\cite[(3.28),
  (3.29), (4.34) and (4.35)]{BK:GradedKL}.
\end{Remarks}

We end this section by using \autoref{T:IntegralKLR} to give an upper bound for the
nilpotency index of the KLR generators $y_1,\dots,y_n$.
If $1\le r\le n$ and $\bi\in I^n$ set
$$\Dr=\set{c_r(\t)-\i_r|\t\in\Std(\bi)}$$ and
define $d_r(\bi)=\#\Dr$. For example,
$\Dr[1]\subseteq\{\kappa_1-\i_1,\dots,\kappa_\ell-\i_1\}$ so that
$d_1(\bi)=(\Lambda,\alpha_{i_1})$.

Two nodes $\gamma=(l,r,c)$ and $\gamma'=(l',r',c')$ are on the same
\textbf{diagonal} if they have the same content. That is, $\gamma$ and
$\gamma'$ are on the same diagonal if and only if $l=l'$ and
$c-r=c'-r'$. The set of diagonals is indexed by pairs $(l,d)$, with
$1\le l\le\ell$ and $d\in\Z$, and where the corresponding diagonal is
the set of nodes $\mathbb{D}_{l,d}=\set{(l,r,c)|\kappa_l+c-r=d}$.
Hence, $d_r(\bi)=\#\Dr$ counts the number of different diagonals
that~$r$ appears on in~$\Std(\bi)$. More precisely, we have:

\begin{Lemma}\label{L:Explicitdr}
  Suppose that $1\le r\le n$ and $\bi\in I^n$. Then
  $$d_r(\bi)=\#\set{(l,d)|d\equiv i_r\pmod e \text{ and }
         \t^{-1}(r)\in\mathbb{D}_{l,d} \text{ for some }\t\in\Std(\bi)}.$$
  That is, $d_r(\bi)$ is equal to the number of distinct diagonals
  that~$r$ appears on for some tableau $\t\in\Std(\bi)$.
\end{Lemma}

The next result is a stronger version of \autoref{L:yONilpotence}. We do not
know how to prove this result using only the relations in~$\RO$.

\begin{Proposition}\label{P:nilpotent}
  Suppose that $1\le r\le n$ and $\bi\in I^n$.
  Then
  $$\prod_{c\in\Dr}(\yo_r-[c])\fo=0.$$
\end{Proposition}

\begin{proof}By \autoref{L:MurphyIdempotent} and \autoref{L:PsiExpansion},
  \begin{align*}
  \prod_{c\in\Dr}(\yo_r-[c])\fo
    &=\sum_{\t\in\Std(\bi)}\prod_{c\in\Dr}(\yo_r-[c])\frac1{\gamma_\t}f_{\t\t}\\
    &=\sum_{\t\in\Std(\bi)}\frac1{\gamma_\t}\prod_{c\in\Dr}([c_r(\t)-\i_r]-[c])f_{\t\t}
     =0,
\end{align*}
where the last equality follows because $c_r(\t)-\i_r\in\Dr$,
for all $\t\in\Std(\bi)$.
\end{proof}

Even though \autoref{P:nilpotent} is very easy to prove within our framework, it
gives strong information about the nilpotency index of $y_re(\bi)$, for
$\bi\in I^n$ and $1\le e\le n$. By \autoref{P:nilpotent}, and
\autoref{C:ROField}, we have the following.

\begin{Corollary}\label{C:GeneralNilpotence}
  Suppose that $\bi\in I^n$ and $1\le r\le n$.  Then
  $y_r^{d_r(\bi)}e(\bi)=0$ in~$\R$.
\end{Corollary}

When $e=0$ Brundan and Kleshchev \cite[Conjecture~2.3]{BK:GradedKL} conjectured
that $y_r^\ell=0$, for $1\le r\le n$. Hoffnung and Lauda proved this
conjecture as the main result of their paper~\cite{HoffnungLauda:KLRnilpotency}.
Using \autoref{C:GeneralNilpotence} we obtain a quick proof of this result and,
at the same time, a generalization of it to include the cases when $e>n$.

\begin{Corollary}\label{C:HoffnungLauda}
  Suppose that $e=0$ or $e>n$. Then $y_r^{\ell}=0$,  for $1\leq r\leq n$.
\end{Corollary}

\begin{proof}
  If $e=0$ then we may assume that $e\gg0$ by \autoref{C:Largee}, so
  the results of this section apply. Hence, we may assume that $e>n$.

  To prove the corollary it is enough to show that $y_r^\ell e(\bi)=0$,
  whenever $\bi=\res(\t)$ for some standard tableau $\t\in\Std(\Parts)$.
  By \autoref{C:GeneralNilpotence}, this will follow if we show that each
  component contains at most one diagonal with content congruent
  to~$i_r$ upon which~$r$ can appear in any standard tableau $\s$ with
  $\res(\s)=\bi$. Suppose by way of contradiction that there exists a
  standard tableau~$\s$, with $\res(\s)=\bi$, and such that~$r$ appears
  in the same component of~$\s$ and~$\t$ but on different diagonals.
  Then the axial distance between the nodes $\s^{-1}(r)$ and
  $\t^{-1}(r)$ is at least~$e$, so every residue in~$I$ must appear in
  any connected path between these two nodes. As
  $\res(\s)=\bi=\res(\t)$ it follows that $\{i_1,\dots,i_r\}=I$. This
  is a contradiction, however, because $|I|=e>n\ge r$.
\end{proof}

\section{Integral bases for $\HO$}\label{Chap:PsiBases}
  Now that we have proved \autoref{Thm:KLRDeformation}, we begin to use the
  machinery of seminormal forms to study the cyclotomic quiver Hecke
  algebras~$\R$.  In this chapter we reconstruct the `natural' homogeneous bases for the
  cyclotomic Hecke algebras $\H(K)$ and their Specht modules over a field.

\subsection{The $\psi$-basis}\label{S:PsiBasis}
\autoref{T:IntegralKLR} links the KLR grading on $\H\cong\R$ with the semisimple
representation theory of~$\HKK$. We next want to try and understand the graded
Specht modules
of~$\H$~\cite{BKW:GradedSpecht,HuMathas:GradedCellular,KMR:UniversalSpecht} in
terms of the seminormal form.  We start by lifting the homogeneous basis
$\{\psi_{\s\t}\}$ of~$\H$ to~$\HO$.  This turns out to be easier than the
approach taken in \cite{HuMathas:GradedCellular}. Throughout this section, $\O$
is an $e$-idempotent subring of~$\K$.

By \autoref{T:IntegralKLR}, there is a unique anti-isomorphism~$\diamond$
of~$\HO$ such that
$$(\psio_r)^\diamond=\psio_r,\quad (\yo_s)^\diamond =\yo_s\quad\text{and}\quad
(\fo)^\diamond=\fo,$$ for $1\le r<n$, $1\le s\le n$ and $\bi\in I^n$.
\autoref{L:PsiExpansion} shows that, in general, the automorphisms $\bigast$ and
$\diamond$ do not coincide.

Recall from \autoref{D:SeminormalBasis} that a $\diamond$-seminormal basis of $\HK$
is a basis $\{f_{\s\t}\}$ of two-sided eigenvalues for~$\L$ such that
$f_{\s\t}=f_{\t\s}^\diamond$, for all $(\s,\t)\in\SStd(\Parts)$. We define a
\textbf{$\diamond$-seminormal coefficient system} to be a set of scalars
$\{\beta_r(\t)\}$ that satisfy the identity in \autoref{L:BetaSquared} and the
``braid relations'' of \autoref{D:SNCS}(a) (with $\alpha$ replaced
by~$\beta$) as well as the relation \autoref{D:SNCS}(b) (with $\alpha$ replaced
by~$\beta$). The reader may check that the $\diamond$-seminormal coefficients
correspond to the more general setup considered in \autoref{R:OIsomorphism}(b).

The main difference between a $\bigast$-seminormal basis and a $\diamond$-seminormal
basis is that $T_rf_{\s\t}=(f_{\t\s}T_r)^{\bigast}$ for a
$\bigast$-seminormal basis whereas $\psio_r f_{\s\t}=(f_{\t\s}\psio_r)^\diamond$ for a
$\diamond$-seminormal basis.

\begin{Lemma}\label{L:BetaSNCS}
  Suppose that $\{f_{\s\t}\}$ is a $\diamond$-seminormal
  basis of~$\HK$. Then there exists a unique $\diamond$-seminormal coefficient
  system $\{\beta_r(\t)\}$ such that if $1\le r<n$ and $(\s,\t)\in\SStd(\Parts)$
  then
  $$f_{\s\t}\psio_r=\beta_r(\t)f_{\s\v}
       -\delta_{i_ri_{r+1}}\frac{t^{\i_{r+1}-c_{r+1}(\t)}}{[\rho_r(\t)]}f_{\s\t},$$
  where $\v=\t(r,r+1)$ and $\t\in\Std(\bi)$, for $\bi\in I^n$. Conversely, as
  in \autoref{T:Seminormal}, a $\diamond$-seminormal coefficient system, together
  with a choice of scalars $\set{\gamma_{\tlam}|\blam\in\Parts}$, determines a
  unique $\diamond$-seminormal basis.
\end{Lemma}

\begin{proof}By \autoref{E:beta},  a set of scalars $\{\beta_r(\t)\}$ is a
  $\diamond$-seminormal coefficient system if and only if
  $\{\alpha_r(\t)\}$ is a $\bigast$-seminormal coefficient system, where
  $$\alpha_r(\t)=\begin{cases}
       \beta_r(\t)t^{c_{r}(\t)-\i_{r}}[1-\rho_r(\t)],
          &\text{if }i_r=i_{r+1},\\[3mm]
      \dfrac{\beta_r(\t)t^{\i_r-c_{r+1}(\t)}}{[\rho_r(\t)]},
          &\text{if }i_r=i_{r+1}+1,\\[3mm]
       \dfrac{\beta_r(\t)[1-\rho_r(\t)]}{[\rho_r(\t)]},
          &\text{otherwise}.
  \end{cases}$$
  Therefore, as seminormal coefficient systems are determined by the
  action of the corresponding generators of~$\H$ on its right regular
  representation, the result follows from
  \autoref{T:Seminormal} and \autoref{L:PsiExpansion}.
\end{proof}

Henceforth, we will work with $\diamond$-seminormal bases.  \autoref{L:BetaSNCS}
also describes the left action of $\psio_r$ on the $\diamond$-seminormal
basis because $\psio_r f_{\s\t}=(f_{\t\s}\psio_r)^\diamond$.

Exactly as in \autoref{T:Seminormal}, if $\{f_{\s\t}\}$ is a $\diamond$-seminormal
basis then there exists scalars $\gamma_\t\in\K$ such that
$f_{\s\t}f_{\u\v}=\delta_{\u\t}\gamma_\t f_{\s\v}$, for
$(\s,\t),(\u,\v)\in\SStd(\Parts)$.  Repeating the argument of
\autoref{C:GammaRecurrence}, these scalars satisfy the following recurrence
relation.

\begin{Corollary}\label{C:GradedGammaRecurrence}
   Suppose that $\t\in\Std(\Parts)$ and
   that $\v=\t(r,r+1)$ is standard, where $1\le r<n$. Then
   $\beta_r(\v)\gamma_\t=\beta_r(\t)\gamma_\v$.
\end{Corollary}

Motivated by \cite{HuMathas:GradedCellular}, we now define a new basis of $\HO$
that is cellular with respect to the anti-involution~$\diamond$.  Fix $\blam\in\Parts$
and let $\bi^\blam=(i^\blam_1,\dots,i^\blam_n)$, so that
$i^\blam_r=\res_{\tlam}(r)$ for $1\le r\le n$. Following
\cite[Definition~4.7]{HuMathas:GradedCellular}, define
$$\Add_\blam(r)=\SetBox[70]{\alpha}{$\alpha$ is an addable $i^\blam_r$-node of the
multipartition $\Shape(\tlam_{\downarrow r})$ that is \textit{below} $(\tlam)^{-1}(r)$},
$$
for $1\le r\le n$.

Up until now we have worked with an arbitrary seminormal basis of~$\HKK$. In
order to define a `nice' basis of $\HO$ that is compatible with
\autoref{T:IntegralKLR} we now fix the choice of $\gamma$-coefficients by
requiring that
\begin{equation}\label{E:Glam}
  \gamma_{\tlam} = \prod_{r=1}^n\prod_{\alpha\in\Add_\blam(r)}
  [\cont_r(\tlam)-c_\alpha],
\end{equation}
for all $\blam\in\Parts$. Together with a choice of seminormal coefficient system,
this determines $\gamma_\t$ for all $\t\in\Std(\Parts)$ by
\autoref{C:GradedGammaRecurrence}. By definition, $\gamma_{\tlam}$ is typically a
non-invertible element of~$\O$. Nonetheless, if $\bi\in I^n$ then
$\fo=\sum_{\s\in\Std(\bi)}\tfrac1{\gamma_\s}f_{\s\s}$ belongs to~$\HO$ by
\autoref{L:MurphyIdempotent}.

We also fix a choice of seminormal coefficient system by requiring that
$\beta_r(\s)=1$ whenever $\s\gdom\t=\s(r,r+1)$, for $\s\in\Std(\Parts)$ and
$1\le r<n$. More precisely, if $\bi\in I$ and $\s\in\Std(\bi)$ then
we define
\begin{equation}\label{E:NewBetas}
\beta_r(\s)= \begin{cases}
  1,&\text{if $\s\gdom\t$ or $i_r\noedge i_{r+1}$},\\
  -\frac{t^{2\i_r-2c_{r+1}(\s)}}{[\rho_r(\s)]^2},
            &\text{if $\t\gdom\s$ and $i_r=i_{r+1}$},\\
  t^{c_r(\s)+c_{r+1}(\s)-\i_r-\i_{r+1}}[1{-}\rho_r(\s)][1{+}\rho_r(s)],
        &\text{if $\t\gdom\s$ and $i_r\leftrightarrows i_{r+1}$},\\
  t^{c_{r}(\s)-\i_{r}}[1{-}\rho_r(\s)],
            &\text{if $\t\gdom\s$ and $i_r\leftarrow i_{r+1}$},\\
  t^{c_{r+1}(\s)-\i_{r+1}}[1{+}\rho_r(s)],
            &\text{if $\t\gdom\s$ and $i_r\rightarrow i_{r+1}$}.\\
\end{cases}
\end{equation}
where $\s\in\Std(\Parts)$ and $\t=\s(r,r+1)$ is standard, for $1\le r<n$.  The
reader is invited to check that this defines a $\diamond$-seminormal coefficient
system.  As the definition of $\psio_r$ is independent of the choice of
seminormal coefficient system this choice is not strictly necessary for what
follows but it simplifies many of the formulas.

By \autoref{L:BetaSNCS}, this choice of $\diamond$-seminormal coefficient system
and $\gamma$-coefficients determines a unique $\diamond$-seminormal basis $\{f_{\s\t}\}$
of~$\HK$. We will use this basis to define new homogeneous basis of $\H$. The
first step is to define
\begin{align*}
  \ylam\flam &=\prod_{r=1}^n\prod_{\alpha\in\Add_\blam(r)}
           t^{-c_r(\tlam)}(L_r-[c_\alpha])\flam\\
  &=\prod_{r=1}^n\prod_{\alpha\in\Add_\blam(r)}
           t^{\i^\blam_r-c_r(\tlam)}(\yo_r-[c_\alpha-\i^\blam_r])\flam,
\end{align*}
where the second equation follows by rewriting $L_k\fo$ in terms of $y_k\fo$ as
in the proof of \autoref{P:OGeneration}. In particular, these equations show
that $\ylam\flam\otimes_\O1_K$ is a monomial in $y_1,\dots,y_n$ and, further,
that it is (up to a sign) equal to the element~$y^\blam$ defined in
\cite[Definition~4.15]{HuMathas:GradedCellular}.

The next result is a essentially a translation of
\cite[Lemma~4.13]{HuMathas:GradedCellular} into the current setting for the
special case of the tableau~$\tlam$.

\begin{Lemma}\label{L:ylam}
  Suppose that $\blam\in\Parts$. Then there exist scalars $a_\s\in\K$ such that
  $$\ylam\flam=f_{\tlam\tlam}+\sum_{\s\Gdom\tlam}
               a_\s f_{\s\s}.$$
  In particular, $\ylam\flam$ is a non-zero element of $\HO$.
\end{Lemma}

\begin{proof}
  By \autoref{L:MurphyIdempotent}, $\flam=\sum_\s\frac1{\gamma_\s}f_{\s\s}$, so that
  $\ylam\flam=\sum_{\s\in\Std(\ilam)}a_\s f_{\s\s}$, for some $a_\s\in\K$, by
  \autoref{E:SeminormalForm}. It remains to show that $a_{\tlam}=1$ and that
  $a_\s\ne0$ only if $\s\Gedom\tlam$. Using \autoref{E:SeminormalForm}, and recalling
  the definition of $\gamma_{\tlam}$ from \autoref{E:Glam},
  $$\frac1{\gamma_{\tlam}}\ylam f_{\tlam\tlam}
       =\frac1{\gamma_{\tlam}}\prod_{r=1}^n\prod_{\alpha\in\Add_\blam(r)}
       t^{-c_r(\tlam)}([c_r(\tlam)]-[c_\alpha])\cdot f_{\tlam\tlam}
       =f_{\tlam\tlam}.
  $$
  To complete the proof we claim that there exist scalars $a_\s(k)\in\K$, $1\le k\le n$,
  such that
  $$\prod_{r=1}^k\prod_{\alpha\in\Add_\blam(r)} t^{-c_r(\tlam)}(L_r-[c_\alpha])\flam
           =\sum_{\substack{\s\in\Std(\ilam)\\\s_{\downarrow k}\Gedom\tlam_{\downarrow k}}}
                a_\s(k)f_{\s\s}$$
  where $a_{\tlam}(k)=1$. We prove this by induction on~$k$.
  If $k=1$ then the result is immediate from~\autoref{E:SeminormalForm}. Suppose that $k>1$.
  By induction, it is enough to show that
  $$
     (L_k-[c_\alpha])f_{\s\s} =([c_\alpha]-[c_k(\s)])f_{\s\s}=0
  $$
  whenever $\s_{\downarrow(k-1)}\Gedom\tlam_{\downarrow(k-1)}$ and
  $\s_{\downarrow k}\rlap{\hspace*{1mm}$\not$}\Gedom\tlam_{\downarrow k}$, for
  $\s\in\Std(\ilam)$. Fix such a tableau~$\s$. Since
  $\s_{\downarrow(k-1)}\Gedom\tlam_{\downarrow(k-1)}$ we must have
  $(\s_{\downarrow k})^{(l)}=\emptyset$ whenever $l>\comp_{\tlam}(k)$, so the node
  $\alpha=\s^{-1}(k)$ must be below $(\tlam)^{-1}(k)$. Therefore,
  $\alpha\in\Add_\blam(k)$, and $c_k(\s)=\cont_\alpha$ for this $\alpha$,
  and forcing $a_\s(k)=0$ as claimed. This completes the proof.
\end{proof}

For each $w\in\Sym_n$ we now fix a reduced expression $w=s_{r_1}\dots s_{r_k}$
for~$w$, with $1\le r_j<n$ for $1\le j\le k$, and define
$\psio_w=\psio_{r_1}\dots\psio_{r_k}$. By \autoref{T:IntegralKLR} the elements
$\psio_r$ do not satisfy the braid relations so, in general, $\psio_w$ will
depend upon this (fixed) choice of reduced expression.

\begin{Definition}
  Suppose that $\blam\in\Parts$. Define
  $$\psio_{\s\t}=(\psio_{d(\s)})^\diamond \ylam\flam\psio_{d(\t)},$$
for $\s,\t\in\Std(\blam)$.
\end{Definition}

We can now lift the graded cellular basis of \cite[Definitions~5.1]{HuMathas:GradedCellular} to~$\HO$.

\begin{Theorem}\label{T:OIntegralBasis}
  Suppose that $\O$ is an idempotent subring. Then
  $$\set{\psio_{\s\t}|\s,\t\in\Std(\bmu)\text{ for }\bmu\in\Parts}$$
  is a cellular basis of $\HO$ with respect to the
  anti-involution~$\diamond$.
\end{Theorem}

\begin{proof}
  In view of \autoref{E:SeminormalForm} and \autoref{L:PsiExpansion},
  \autoref{L:ylam} implies that
  \begin{equation}\label{E:OTriangular}
    \psio_{\s\t}= f_{\s\t} +\sum_{(\u,\v)\Gdom(\s,\t)}a_{\u\v}f_{\u\v},
  \end{equation}
  for some $a_{\u\v}\in\K$. Therefore,
  $\set{\psio_{\s\t}|(\s,\t)\in\SStd(\Parts)}$ is a basis of
  $\HK$. In fact, these elements are a basis for $\HO$ because if $h\in\HO$ then
  we can write $h=\sum r_{\u\v}f_{\u\v}$, for some $r_{\u\v}\in\K$. Pick
  $(\s,\t)$ to be minimal with respect to dominance such that $r_{\s\t}\ne0$. Then
$r_{\s\t}\in\O$ because $h\in\HO$. Consequently, $h-r_{\s\t}\psio_{\s\t}\in\HO$ so, by continuing in this way, we can write~$h$ as a linear combination of the $\psi$-basis.

  It remains to show that the $\psi$-basis is cellular with respect to the
  anti-involution~$\diamond$.  By definition, if $\blam\in\Parts$ then $\ylam$ and
  $\flam$ commute and they are fixed by the automorphism~$\diamond$. Therefore,
  $(\psio_{\s\t}\)^\diamond=\psio_{\t\s}$, for all $\s,\t\in\Std(\blam)$. By
  \autoref{L:BetaSNCS}, the $\diamond$-seminormal basis $\{f_{\s\t}\}$ is a
  cellular basis with cellular anti-involution~$\diamond$. It remains to verify (GC$_2$)
  from \autoref{D:cellular}. As in \autoref{T:Seminormal}, the seminormal basis
  $\{f_{\u\v}\}$ is cellular. Therefore, if $(\s,\t)\in\SStd(\blam)$
  and $h\in\HO$ then, using \autoref{E:OTriangular} twice,
  \begin{align*}
    \psio_{\s\t}h &=(\psio_{d(\s)})^\diamond\psio_{\tlam\t}
    \equiv(\psio_{d(\s)})^\diamond\Big(f_{\tlam\t}+\sum_{\v\gdom\t}a_\v f_{\tlam\v}\Big)h
    \equiv(\psio_{d(\s)})^\diamond\sum_{\v\in\Std(\blam)}a'_\v f_{\tlam\v}\\
    &\equiv(\psio_{d(\s)})^\diamond \sum_{\v\in\Std(\blam}b_\v\psio_{\tlam\v}\quad
    \equiv \sum_{\v\in\Std(\blam}b_\v\psio_{\s\v}\quad
      \pmod\Hlam,
  \end{align*}
  where $a_\v,a_\v'\in\K$ and $b_\v\in\O$ with the scalars $b_\v$ being
  independent of~$\s$. Hence, (GC$_2$) holds, completing the proof.
\end{proof}

If $K=\O/\m$ for some maximal ideal $\m$ of~$\O$ then
$\H(K)\cong\H(\O)\otimes_\O K$. Set $\psi_{\s\t}=\psio_{\s\t}\otimes1_K$.

\begin{Corollary}[\protect{\cite[Theorem~5.8]{HuMathas:GradedCellular}}]
  \label{C:PsiBasis}
  Suppose that $K=\O/\m$  for some maximal ideal $\m$ of~$\O$. Then
  $\set{\psi_{\s\t}|\s,\t\in\Std(\bmu)\text{ for }\bmu\in\Parts}$
  is a graded cellular basis of $\H(K)$ with $\deg\psi_{\s\t}=\deg\s+\deg\t$,
  for $(\s,\t)\in\SStd(\Parts)$.
\end{Corollary}


\subsection{Graded Specht modules and Gram determinants}\label{S:GradedSpecht}
By \autoref{T:OIntegralBasis}, $\{\psio_{\s\t}\}$ is a cellular basis of $\HO$
so we can use it to define Specht modules for $\HO$ that specialise to the
graded Specht modules in characteristic zero and in positive characteristic.

\begin{Definition}
  Suppose that $\blam\in\Parts$. The Specht module $S^\blam(\O)$ is the right
  $\HO$-module with basis $\set{\psio_\t|\t\in\Std(\blam)}$, where
  $\psio_\t=\psio_{\tlam\t}+\Hlam(\O)$.
\end{Definition}

By \autoref{T:OIntegralBasis} and
\cite[Corollary~5.10]{HuMathas:GradedCellular}, ignoring the grading,
$S^\blam(\O)\otimes_\O K$ can be identified with the graded Specht module
$S^\blam$ of $\H$ defined by Brundan, Kleshchev and
Wang~\cite{BKW:GradedSpecht}. The action of $\HK$ on a graded Specht module is
completely determined by the relations for these modules that are given in
\cite{KMR:UniversalSpecht}. In contrast, in view of \autoref{E:OTriangular} and
\autoref{T:IntegralKLR}, the action of $\HO$ on the Specht module $S^\blam(\O)$
is completely determined by the (choice of) seminormal form.

We now turn to computing the determinant of the Gram matrix
$$\Gram = \(\<\psio_\s,\psio_\t\>\)_{\s,\t\in\Std(\blam)}.$$
\textit{A priori}, it is unclear how the bilinear form on $S^\blam(\O)$ is
related to the usual (ungraded) bilinear from on the Specht module that is
defined using the Murphy basis that we considered in
\autoref{T:GramDetFactorization}. The main problem in relating these two
bilinear forms is that the cellular algebra anti-involutions~$\bigast$ and~$\diamond$, which are used to define these bilinear forms, are different.

Note that the cellular algebra anti-involutions~$\bigast$ and~$\diamond$ on
$\HO$ naturally extend to anti-involutions on the algebra $\HKK$.  The
key point to understanding the graded bilinear form is the following.

\begin{Lemma}\label{L:FtstarInvariance}
  Suppose that $\t\in\Std(\Parts)$. Then $(F_\t)^\diamond=F_\t$.
\end{Lemma}

\begin{proof}
  By definition, $F_\t$ is a linear combination of products of Jucys-Murphy elements, so it can also be
  written as a polynomial, with coefficients in~$\K$, in $\yo_r$,
  $f_{\bi}^{\O}$, for $1\le r\le n$ and $\bi\in I^n$. As $(\yo_r)^\diamond=\yo_r$,
  $(f_{\bi}^{\O})^{\diamond}=f_{\bi}^{\O}$, for $1\le r\le n$ and $\bi\in I^n$, the
  result follows.
\end{proof}

Recall that if $\t\in\Std(\blam)$ then $\psio_\t=\psi^\O_{\tlam\t}+\Hlam$ is a
basis element of the Specht module $S^\blam(\O)$. In order to compute
$\det\Gram$, set $f_\t=\psio_\t F_\t$, for $\t\in\Std(\blam)$.
Recall that $S^\blam(\K)=S^\blam(\O)\otimes_{\O}\K$.

\begin{Lemma}\label{L:GramDetTransfer}
  Suppose that $\blam\in\Parts$. Then $\set{f_\t|\t\in\Std(\blam)}$ is a basis
  of~$S^\blam(\K)$. Moreover,
  $\det\Gram=\det\(\<f_\s,f_\t\>\)=\prod_{\s\in\Std(\blam)}\gamma_{\s}$.
\end{Lemma}

\begin{proof}
  By definition, $f_\t=f_{\tlam\t}+(\HKK)^{\rhd\blam}$. Therefore, $f_\t\in
  S^\blam(\K)$ and $f_\t=\psio_\t+\sum_{\v\gdom\t}r_{\t\v} \psio_\v$ by
  \autoref{E:OTriangular}, for some scalars $r_{\t\v}\in\K$. Set $r_{\t\t}=1$
  and $U=\(r_{\t\v}\)$. Then $\set{f_\t|\t\in\Std(\blam)}$ is a $\K$-basis
  of~$S^\blam(\K)$ and $\Gram=(U^{-1})^{tr}\(\<f_\s,f_\t\>\)U^{-1}$ Taking
  determinants shows that $\det\Gram=\deg\(\<f_\s,f_\t\)\)$ since $U$ is
  unitriangular.  To complete the proof observe that
  $\<f_\s,f_\t\>f_{\tlam\tlam}\equiv
  f_{\tlam\s}f_{\t\tlam}=\delta_{\s\t}\gamma_\s f_{\tlam\tlam}\pmod{\Hlam}$,
  where we are implicitly using \autoref{L:FtstarInvariance}. The result
  follows.
\end{proof}

\autoref{L:GramDetTransfer} is subtly different from
\autoref{E:GramDetGammaProd} because, in spite of our notation,
the~$\gamma_\t$'s appearing in the two formulas satisfy different
recurrence relations. It is not hard to show that the quotient of
$\gamma_\t$, as defined in this section, by the~$\gamma_\t$ defined in
\autoref{S:GramDet} in a unit in~$\O$, for all $\t\in\Std(\Parts)$.

\begin{Lemma}\label{L:KLRGamma}
  Suppose that $\t\in\Std(\blam)$, for $\blam\in\Parts$. Then
  $\gamma_\t=u_\t\Phi_e(t)^{\deg_e(\t)}$, for some unit $u_\t\in\O^\times$.
\end{Lemma}

\begin{proof} We argue by induction on the dominance order on $\Std(\blam)$. If
  $\t=\tlam$ then \autoref{E:Glam} ensures that
  $\gamma_{\tlam}=u_{\tlam}\Phi_e(t)^{\deg_e(\tlam)}$, for some unit
  $u_{\tlam}\in\O$. Now suppose that $\tlam\gdom\t$. Then there exists a
  standard tableau $\s\in\Std(\blam)$ such that $\s\gdom\t$ and $\t=\s(r,r+1)$,
  where $1\le r<n$. Arguing exactly as in \autoref{C:GradedGammaRecurrence}
  shows that $\beta_r(\s)\gamma_\t=\beta_r(\t)\gamma_\s$. Therefore,
  $\gamma_\t=\tfrac{\beta_r(\t)}{\beta_r(\s)}\gamma_\s=\beta_r(\t)\gamma_\s$.
  Hence, the lemma follows by induction exactly as in the proof of
  \autoref{T:GramDetFactorization}.
\end{proof}

\begin{Remark}
  Looking at the definition of a $\diamond$-seminormal coefficient system shows
  that the quantities $\tfrac{\beta_r(\t)}{\beta_r(\s)}$, which are used in the
  proof of \autoref{L:KLRGamma}, are independent of the choice of
  $\diamond$-seminormal coefficient system. This shows that the choice of
  $\diamond$-seminormal coefficient system made in \autoref{E:NewBetas} really is
  only for convenience.
\end{Remark}

By general nonsense, the determinants of $\Gram$ and $\UnGram$ differ by a scalar
in~$\K$. The last two results readily imply the next theorem, the real content
of which is that this scalar is a unit in~$\O$.

\begin{Theorem}\label{T:eGramDet}
  Suppose that $\blam\in\Parts$. Then $\det\Gram=u\Phi_e(t)^{\deg_e(\blam)},$
  for some unit $u\in\O^\times$. Consequently,
  $\det\Gram=u'\det\UnGram$, for some unit $u'\in\O^\times$.
\end{Theorem}

If $\bi\in I^n$ and $\blam\in\Parts$ let
$\Std_\bi(\blam)=\set{\t\in\Std(\blam)|\res(\t)=\bi}$.

The Specht module $S^\blam$ over $\O$ decomposes as a direct sum of generalised
eigenspaces as an $\L(\O)$-module: $S^\blam=\bigoplus_{\bi\in I^n} S_\bi^\blam$,
where $S^\blam_\bi=S^\blam\fo$.  The weight space $S^\blam_\bi$ has basis
$\set{\psio_\t|\t\in\Std_\bi(\blam)}$ and the bilinear linear
form $\<\ ,\ \>$ on $S^\blam$ respects the weight space decomposition of
$S^\blam$. Set
$$\deg_{e,\bi}(\blam)=\sum_{\t\in\Std_\bi(\blam)} \deg\t.$$
and let $\Gram_\bi$ be restriction of the Gram matrix of $S^\blam$ to
$S^\blam_\bi$, for $\bi\in I^n$. Then we have the following refinement of
\autoref{T:eGramDet} (and \autoref{T:GramDetFactorization}).

\begin{Corollary}
  Suppose that $\blam\in\Parts$ and $\bi\in I^n$. Then
  $\deg\Gram_\bi=u_\bi\Phi_e(t)^{\deg_{e,\bi}(\blam)}$, for some unit $u_\bi\in\O^{\times}$.
  Moreover, $\deg_{e,\bi}(\blam)\ge0$.
\end{Corollary}

\section{A distinguished homogeneous basis for $\H$}\label{Chap:BBasis}
The $\psi$-basis of $\H(\O)$, the homogeneous bases of $\H$ constructed in
\cite{HuMathas:GradedCellular}, and the homogeneous basis of the graded Specht
modules given by Brundan, Kleshchev and Wang~\cite{BKW:GradedSpecht}, are all
indexed by pairs of standard tableaux. Unfortunately, unlike in the ungraded case,
these basis elements depend upon choices of reduced expressions for the
permutations corresponding to these tableaux. In this section we construct new
bases for these modules that depend only on the corresponding tableaux.

\subsection{A new basis of $\H(\O)$} To construct our new basis for
$\H$ we need to work over a complete discrete valuation ring. We start by
setting up the necessary machinery.

Recall that the algebra $\H$ is defined over the field~$K$ with parameter~$\xi$
and that $e>1$ is minimal such that $[e]_\xi=0$. Let $x$ be an indeterminate
over~$K$ and let $\O=K[x]_{(x)}$ and $t=x+\xi$. Then $(\O,t)$ is an
idempotent subring by \autoref{Ex:indeterminate}(b) and $K(x)$ is the field of
fractions of~$\O$. Note that $\O$ is a local ring with maximal ideal $\m=x\O$.

Let $\OO$ be the $\m$-adic completion of $\O$. Then $\OO$ is a complete discrete
valuation ring with field of fractions $K((x))$
Let $\KK=K((x))$ be the $\m$-adic completion of~$K(x)$. Then $\OO$ is an
idempotent subring of $\KK$.

Define a valuation on $\KK^\times$ by setting $\nu_x(a)=n$ if $a=ux^n$, where
$n\in\Z$ and $u\in\OO^\times$ is a unit in~$\OO$.  We need to work with a complete
discrete valuation ring because of the following fundamental but elementary fact
that is proved, for example, as \cite[Proposition~II.5]{Serre:LocalFields}.

\begin{Lemma}\label{L:madicExpansion}
  Suppose that $a\in\KK$. Then $a$ can be written uniquely as a convergent series
  $$a=\sum_{n\in\Z}a_nx^n,\qquad\text{with }a_n\in K,$$
  such that if $a\ne0$ then $a_n\ne0$ only if $n\ge\nu_x(a)$. Moreover, $a\in\OO$
  if and only if~$a_n=0$ for all $n<0$.
\end{Lemma}

In particular, $x^{-1}K[x^{-1}]\cap\OO=0$, where we embed $x^{-1}K[x^{-1}]$
into~$\KK$ in the obvious way.

\begin{Theorem}\label{T:KLbasis}
    Suppose that $(\s,\t)\in\SStd(\Parts)$. There exists a unique
  element $\Boo_{\s\t}\in\H(\OO)$ such that
  $$\Boo_{\s\t}=f_{\s\t}+\sum_{\substack{(\u,\v)\in\SStd(\Parts)\\(\u,\v)\Gdom(\s,\t)}}
  p_{\u\v}^{\s\t}(x^{-1})f_{\u\v},$$
  where $p_{\u\v}^{\s\t}(x)\in xK[x]$. Moreover,
  $\set{\Boo_{\s\t}|(\s,\t)\in\SStd(\Parts)}$
  is a cellular basis of~$\H(\OO)$.
\end{Theorem}

\begin{proof}The existence of an element $\Boo_{\s\t}$ with the required
  properties follows directly from \autoref{E:OTriangular} and
  \autoref{L:madicExpansion} using Gaussian elimination. (See the proof of \autoref{P:KLbasis},
  below, which proves a stronger result in characteristic zero.)  To prove uniqueness of the element
  $\Boo_{\s\t}$, suppose, by way of contradiction, that there exist two elements
  $\Boo_{\s\t}$ and $B_{\s\t}'$ in~$\H(\OO)$ with the required properties. Then
  $\Boo_{\s\t}-B'_{\s\t}=\sum r_{\u\v} f_{\u\v}\in\H(\OO)$ and, by assumption,
  $r_{\u\v}\in x^{-1}K[x^{-1}]$ with $r_{\u\v}\ne0$ only if
  $(\u,\v)\Gdom(\s,\t)$. Pick $(\a,\b)$ minimal with respect to dominance such
  that $r_{\a\b}\ne0$. Then, by \autoref{T:OIntegralBasis}, if we write
  $\Boo_{\s\t}-B'_{\s\t}$ as a linear combination of $\psi$-basis elements then
  $\psi_{\a\b}^{\OO}$ appears with coefficient $r_{\a\b}$. Therefore,
  $r_{\a\b}\in x^{-1}K[x^{-1}]\cap\OO=0$, a contradiction. Hence,
  $\Boo_{\s\t}=B'_{\s\t}$ as claimed.

  By \autoref{E:OTriangular}, the transition matrix between the $B$-basis and
  the $\psi$-basis is unitriangular, so $\{\Boo_{\s\t}\}$ is a basis of
  $\H(\OO)$. To show that the $B$-basis is cellular we need to check properties
  (GC$_1$)--(GC$_3$) from \autoref{D:cellular}. We have already verified
  (GC$_1)$ Moreover, (GC$_3$) holds because $(\Boo_{\s\t})^\diamond=\Boo_{\t\s}$ by
  the uniqueness of $\Boo_{\t\s}$ since $\{f_{\u\v}\}$ is $\diamond$-seminormal
  basis.  It remains to prove (GC$_2$), which we do in three steps.

  \Step 1. We claim that if $h\in\H(\OO)$ and  $\t\in\Std(\blam)$ then
  $$\Boo_{\tlam\t}h \equiv \sum_{\v\in\Std(\blam)} b_\v \Boo_{\tlam\v}
     \quad\pmod\Hlam,$$
  for some scalars $b_\v\in\OO$ that depend only on $\t$, $\v$ and $h$ (and not
  on $\tlam$).

  To see this first note that
  $\psio_{\tlam\t}=f_{\tlam\t}+\sum_{\v\Gdom\t}a_\v f_{\tlam\v}$ by
  \autoref{E:OTriangular}, for some $a_\v\in K(x)$. Therefore, it follows by
  induction on the dominance order that if $\t\in\Std(\blam)$ then
  $$\Boo_{\tlam\t}=f_{\tlam\t}+\sum_{\v\Gdom\t}p_{\t\v}f_{\tlam\v}\quad\pmod\Hlam,$$
  for some $p_{\t\v}\in x^{-1}K[x^{-1}]$. As  the seminormal basis is cellular,
  and the transition matrix between the seminormal basis and the $B$-basis is
  unitriangular, our claim now follows.

  \Step 2. As the Specht module $S^\blam$ is cyclic there exists an element
  $\Do_\t\in\H(\OO)$ such that $\Boo_{\tlam\t}\equiv \Boo_{\tlam\tlam}\Do_\t\pmod\Hlam$.
  We claim that
  $$\Boo_{\s\t}\equiv (\Do_s)^\diamond \Boo_{\tlam\tlam} \Do_\t\pmod\Hlam,$$
  for all $\s,\t\in\Std(\blam)$.

  To prove this claim, embed $\H(\OO)$ in $\H(\KK)$. Note that
  $f_{\tlam\tlam}f_{\u\v}=0$ if $\u\ne\tlam$, so we may assume that
  $\Do_\t\equiv\sum_\v q_{\t\v}f_{\tlam\v}\pmod{\Hlam}$, for some $q_{\t\v}\in\KK$.  Then
  $$\Boo_{\tlam\t}\equiv
  \Boo_{\tlam\tlam}\Do_\t=\sum_{\v\in\Std(\blam)}\gamma_{\tlam}q_{\t\v}f_{\tlam\v}\quad\pmod\Hlam.$$
  Therefore, $q_{\t\v}=\tfrac1{\gamma_{\tlam}}p_{\t\v}$, where $p_{\t\v}\in
  \delta_{\t\v}+x^{-1}K[x^{-1}]$ is as in Step~1. In particular,
  $q_{\t\t}=\tfrac{1}{\gamma_{\tlam}}$ and $q_{\t\v}\ne0$ only if $\v\Gedom\t$. Consequently,
  \begin{align*}
  (\Do_s)^\diamond \Boo_{\tlam\tlam}\Do_\t
  &\equiv\sum_{\substack{(\u,\v)\Gedom(\s,\t)\\\u,\v\in\Std(\blam)}}
              q_{\s\u}q_{\t\v}f_{\u\tlam}f_{\tlam\tlam}f_{\tlam\v}
         =\sum_{(\u,\v)\Gedom(\s,\t)}\gamma_{\tlam}^2q_{\s\u}q_{\t\v}f_{\u\v}\\
        &=f_{\s\t}+\sum_{(\u,\v)\Gdom(\s,\t)}p_{\s\u}p_{\t\v}f_{\u\v}
           \quad\pmod\Hlam.
  \end{align*}
  By construction, $(\Do_s)^\diamond \Boo_{\tlam\tlam}\Do_\t\in\H(\OO)$. Consequently, our claim
  now follows using the uniqueness property of $\Boo_{\s\t}$ since
  $p_{\s\u}p_{\t\v}\in x^{-1}K[x^{-1}]$ when $\s\ne\u$ or $\t\ne\v$.

  \Step 3. We can now verify (GC$_2$). If $h\in\H(\OO)$ then, using steps~1 and~2,
  $$\Boo_{\s\t}h\equiv (\Do_s)^\diamond \Boo_{\tlam\t}h
         \equiv \sum_{\v\in\Std(\blam)}b_\v (\Do_s)^\diamond \Boo_{\tlam\v}
         \equiv \sum_{\v\in\Std(\blam)}b_\v \Boo_{\s\v} \quad\pmod\Hlam,
  $$
  where $b_\v$ depends only on~$\t$, $\v$ and $h$ and not on~$\s$. Hence, the
  $B$-basis satisfies all of the cellular basis axioms and the theorem is proved.
\end{proof}

By \autoref{T:KLbasis}, if $(\s,\t)\in\SStd(\Parts)$ then
$\Boo_{\s\t}\in\H(\OO)$, however, our notation suggests that
$\Boo_{\s\t}\in\H(\O)$, where $\O=K[x]_{(x)}$. The next result justifies our
notation and shows that we can always work over the ring~$\O$.

\begin{Corollary}\label{C:OForm}
  Let $\O=K[x]_{(x)}$. Then $\set{\Boo_{\s\t}|(\s,\t)\in\SStd(\Parts)}$ is a
  graded cellular basis of $\H(\O)$.
\end{Corollary}

\begin{proof} Fix $(\s,\t)\in\SStd(\Parts)$. Then
  it is enough to prove that $\Boo_{\s\t}\in\H(\O)$. First note that
  by construction the $\diamond$-seminormal basis is defined over the rational function
  field $K(x)$, so  $\Boo_{\s\t}$ is defined over the ring $R=K(x)\cap\OO$
  since if $(\u,\v)\in\SStd(\Parts)$ then $p^{\s\t}_{\u\v}(x^{-1})\in K[x^{-1}]\subset K(x)$ by
  \autoref{T:KLbasis}. Every element of $K(x)$
  can be written in the form $f(x)/g(x)$, for $f(x),g(x)\in K[x]$ with
  $\gcd(f,g)=1$. Expanding $f/g$ into a power series, as in
  \autoref{L:madicExpansion}, it is not difficult to see that if $f/g\in\OO$ then
  $g(0)\ne0$. Therefore, $R\subseteq\O$ so
  that $\Boo_{\s\t}$ is defined over $\O$ as claimed.
\end{proof}

By similar arguments, $\Do_\t\in\HO$, for all $\t\in\Std(\Parts)$.

If  $K$ is a field of characteristic zero then we can determine
the degree of the polynomials $p^{\s\t}_{\u\v}\ne0$, for
$(\u,\v)\Gedom(\s,\t)\in\SStd(\Parts)$.

\begin{Proposition}\label{P:KLbasis}
  Suppose that $K$ is a field of characteristic zero. Suppose that $(\u,\v)\Gdom(\s,\t)$ for
  $(\s,\t),(\u,\v)\in\SStd(\Parts)$. Then $p^{\s\t}_{\u\v}(x)\in xK[x]$ and
  $$\deg p^{\s\t}_{\u\v}(x)\le\tfrac12(\deg\u-\deg\s+\deg\v-\deg\t).$$
  In particular, $p^{\s\t}_{\u\v}(x)\ne0$ only if
  $\deg\u+\deg\v>\deg\s+\deg\t$.
\end{Proposition}

\begin{proof}
  We argue by induction on the dominance orders on $\Parts$ and $\Std(\Parts)$.
  Note that $\deg p(x)=d$ if and only if $\nu_x\(p(x^{-1})\)=-d$. For convenience,
  throughout the proof given two tableaux $\s,\u\in\SStd(\Parts)$ set
  $\deg(\s,\u)=\deg\s-\deg\u$. Therefore, the proposition is equivalent to the
  claim that $\nu_x\(p^{\s\t}_{\u\v}(x^{-1})\)\ge\frac12\(\deg(\s,\u)+\deg(\t,\v)\)$.

  Suppose first that $\blam=(n|0|\dots|0)$. Then $\s=\tlam=\t$ and
  $\psio_{\tlam\tlam}=f_{\tlam\tlam}$ so there is nothing to prove.
  Hence, we may assume that $\blam\ne(n|0|\dots|0)$ and that
  the proposition holds for all more dominant shapes.

  Next, consider the case when $\s=\tlam=\t$. By the proof of \autoref{L:ylam},
  if $\s\in\Std(\ilam)$ and $\s\Gedom\tlam$ then
  $\ylam f_{\s\s}=u_\s'\gamma_{\tlam}f_{\s\s}$ for some unit $u'_\s\in\O^\times$.
  Therefore, by \autoref{L:KLRGamma}, there exist units $u_\s\in\O^\times$ so
  that in $\H(K(x))$
  $$\psio_{\tlam\tlam}=\sum_{\s\Gedom\tlam}\frac{u_\s'\gamma_{\tlam}}{\gamma_\s}f_{\s\s}
                  =f_{\tlam\tlam}+\sum_{\s\Gdom\tlam}u_\s
                  \Phi_e(t)^{\deg(\tlam,\s)}f_{\s\s}.$$
  Since $t=x+\xi$, the constant term of $\Phi_e(t)$ is $\Phi_e(\xi)=0$, so $x$
  divides $\Phi_e(t)$ and $\nu_x(\Phi_e(t)^{\deg(\tlam,\s)})=\deg(\tlam,\s)$
  since the coefficient of $x$ in $\Phi_e(t)$ is non-zero.
  (If~$K$ is field of positive characteristic this may not be true.) Expanding each unit
  $u_\s$ into a power series, as in \autoref{L:madicExpansion}, the coefficient
  of $f_{\s\s}$
  can be written as $b_\s+c_\s$ where $b_\s\in x^{-1}K[x^{-1}]$ and $c_\s\in\O$.
  In particular, if $b_\s\ne0$ and $c_\s\ne0$ then $\nu_x(c_\s)\ge0>\nu_x(b_\s)$
  and $\nu_x(c_\s)>\nu_x(b_\s)\ge\deg(\tlam,\s)$.  Pick $\t$ minimal with
  respect to dominance such that $c_\t\ne0$. Note that
  $\nu_x(c_\t)\ge\deg(\tlam,\t)$, with equality only if $b_\t=0$. Using induction, replace
  $\psio_{\tlam\tlam}$ with the element $A_{\tlam\tlam}=\psio_{\tlam\tlam}-c_\t\Boo_{\t\t}$. By
  construction $A_{\tlam\tlam}\in\H(\O)$ and, by \autoref{E:OTriangular}, the coefficient
  of~$f_{\t\t}$ in $A_{\tlam\tlam}$ is~$b_\t\in x^{-1}K[x^{-1}]$.
  If $(\u,\v)\Gedom\(\t,\t)$ then,
  $f_{\u\v}$ appears in $\Boo_{\t\t}$ with coefficient $p^{\t\t}_{\u\v}(x^{-1})$
  and, by induction,
  $\nu_x\(p^{\t\t}_{\u\v}(x^{-1})\)\ge\tfrac12\(\deg(\t,\u)+\deg(\t,\v)\)$. Therefore,
  \begin{align*}
    \nu_x\(c_\t p^{\t\t}_{\u\v}(x^{-1})\)&=\nu_x(c_\t)+\nu_x\(p^{\t\t}_{\u\v}(x^{-1})\)
             \ge\deg(\tlam,\t)+\tfrac12(\deg(\t,\u)+\deg(\t,\v\))\\
            &=\tfrac12\(\deg(\tlam,\u)+\deg(\tlam,\v)\).
  \end{align*}
  It follows that if $f_{\u\v}$ appears in $A_{\tlam\tlam}$ with non-zero
  coefficient $a_{\u\v}$ then $\nu_x(a_{\u\v})\ge\tfrac12\(\deg(\tlam,\u)+\deg(\tlam,\v)\)$.
  If $A_{\tlam\tlam}$ now has the required properties then we can set
  $B_{\tlam\tlam}=A_{\tlam\tlam}$. Otherwise, let $(\s,\t)$ be a pair of tableau
  that is minimal with respect to dominance such that the coefficient of
  $f_{\s\t}$ in $A_{\tlam\tlam}$ is of the form $b_{\s\t}+c_{\s\t}$ with
  $c_{\s\t}\ne0$, $\nu_x(c_{\s\t})\ge0$, $b_{\s\t}\in x^{-1}K[x^{-1}]$ and
  $\nu_x(b_{\s\t})\ge\tfrac12\(\deg(\tlam,\s)+\deg(\tlam,\t)\)$.  Replacing
  $A_{\tlam\tlam}$ with $A_{\tlam\tlam}-c_{\s\t}\Boo_{\s\t}$ and continuing in
  this way we will, in a finite number of steps, construct an
  element~$B'_{\tlam\tlam}$ with all of the required properties. By the uniqueness statement
  in \autoref{T:KLbasis}, $\Boo_{\tlam\tlam}=B'_{\tlam\tlam}$ so this proves the
  proposition for the polynomials $p^{\tlam\tlam}_{\u\v}(x^{-1})$.

  Finally, suppose that $(\s,\t)\in\SStd(\blam)$ with
  $(\tlam,\tlam)\gdom(\s,\t)$.  Without loss of generality, suppose that
  $\s=\a(r,r+1)$ where $\a\in\Std(\bi)$, for $\bi\in I^n$, and $\a\gdom\s$. Using
  \autoref{L:PsiExpansion},
  \begin{align*}
  \psio_r\Boo_{\a\t}
  & = \sum_{(\u,\v)\Gedom(\a,\t)}p^{\a\t}_{\u\v}(x^{-1})\psio_rf_{\u\v}\\
  & = \sum_{(\u,\v)\Gedom(\a,\t)}p^{\a\t}_{\u\v}(x^{-1})\Big(\beta_r(\u)f_{\u(r,r+1),\v}
        -\delta_{i_ri_{r+1}}\frac{t^{\i_{r+1}-c_{r+1}(\u)}}{[\rho_r(\u)]}f_{\u\v}\Big).
  \end{align*}
  By induction, $\nu_x(p^{\a\t}_{\u\v})\ge\frac{1}{2}(\deg(\a,\u)+\deg(\t,\v))$. Therefore,
  using \autoref{L:KLRGamma} (as in the proof of \autoref{T:GramDetFactorization}), it follows that
  if $c_{\u\v}\ne0$ is the coefficient of $f_{\u\v}$ in the last equation then
  $\nu_x(c_{\u\v})\ge\tfrac12\(\deg(\s,\u)+\deg(\t,\v)\)$. Hence, the
  proposition follows by repeating the argument of the last paragraph.
\end{proof}

\subsection{A distinguished homogeneous basis of $\H(K)$}
This section uses \autoref{T:KLbasis} to construct a new graded cellular basis
of~$\H(K)$. The existence of such a basis is not automatically guaranteed by
\autoref{T:KLbasis} because the elements $\Boo_{\s\t}\otimes1_K$, for
$(\s,\t)\in\SStd(\Parts)$, are not necessarily homogeneous.

The isomorphisms $K\cong\O/x\O\cong\OO/x\OO$ extend to $K$-algebra
isomorphisms
$$\H(K)\cong\H(\O)\otimes_\O K\cong\H(\OO)\otimes_{\OO}1_K.$$
We identify these three $K$-algebras.

\begin{Lemma}\label{L:BstApproximation}
  Suppose that $(\s,\t)\in\SStd(\Parts)$. Then
  $$\Boo_{\s\t}\otimes1_K=\psi_{\s\t}+\sum_{(\u,\v)\Gdom(\s,\t)}a_{\u\v}\psi_{\u\v},$$
  for some $a_{\u\b}\in K$. In particular, the homogeneous component of
  $\Boo_{\s\t}\otimes1_K$ of degree $\deg\s+\deg\t$ is non-zero.
\end{Lemma}

\begin{proof}
  This is immediate from \autoref{T:KLbasis}, \autoref{E:OTriangular} and \autoref{C:PsiBasis}.
\end{proof}

Recall from Step~2 in the proof of \autoref{T:KLbasis} that for each
$\v\in\Std(\blam)$ there exists an element
$\Do_\v\in\H(\O)$ such that $\Boo_{\s\t}\equiv(\Do_\s)^\diamond
B_{\tlam\tlam}\Do_\t\pmod\Hlam$.

\begin{Definition}\label{D:BstBasis}
  Suppose that $\blam\in\Parts$.
  \begin{enumerate}
    \item If $\v\in\Std(\blam)$ let $D_\v$ be the homogeneous component of
    $\Do_\v\otimes1_K$ of degree $\deg\v-\deg\tlam$.
    \item Define $B_{\tlam\tlam}$ to be the homogeneous component of
    $\Boo_{\tlam\tlam}\otimes1_K$ of degree $2\deg\tlam$. More generally, if
    $\s,\t\in\Std(\blam)$ define
    $B_{\s\t} = D_\s^\diamond B_{\tlam\tlam} D_\t.$
\end{enumerate}
\end{Definition}

By \autoref{T:KLbasis}, $(B_{\tlam\tlam}^{\O})^\diamond=B_{\tlam\tlam}^{\O}$ which
implies that $B_{\tlam\tlam}^\diamond=B_{\tlam\tlam}$.  Consequently, if
$\s,\t\in\Std(\blam)$ then $B_{\s\t}^\diamond=B_{\t\s}$. If $B_{\s\t}\ne0$ then, by
construction, $B_{\s\t}$ is homogeneous of degree $\deg\s+\deg\t$.
Unfortunately, it is not clear from the definitions that $B_{\s\t}$ is non-zero.

\begin{Proposition}\label{P:BBasisNonZero}
  Suppose that $(\s,\t)\in\SStd(\Parts)$. Then
  $$B_{\s\t}\equiv\psi_{\s\t}+\sum_{(\u,\v)\Gdom(\s,\t)}b_{\u\v}\psi_{\u\v}
     \qquad\pmod\Hlam,$$
  for some $b_{\u\v}\in K$. In particular, $B_{\s\t}\ne0$.
\end{Proposition}

\begin{proof}Fix $\blam\in\Parts$ and suppose that $\s,\t\in\Std(\blam)$. If
  $\s=\t=\tlam$ then $B_{\tlam\tlam}$ is the homogeneous component of
  $\Boo_{\tlam\tlam}\otimes 1_{K}$ of degree $2\deg\tlam$, so the result is just
  \autoref{L:BstApproximation} in this case. Now consider the case when
  $\s=\tlam$ and $\t$ is an arbitrary standard $\blam$-tableau. Then, since
  $B_{\tlam\tlam}^{\O}\equiv\psi_{\tlam\tlam}^{\O}\pmod{\Hlam}$,
  $$\Boo_{\tlam\t}\otimes 1_{K}\equiv(\psi_{\tlam\tlam}^{\O}\otimes1_{K})(\Do_\t\otimes1_K)\pmod\Hlam.$$
  Looking at the homogeneous component of degree $\deg\tlam+\deg\t$ shows that
  \begin{align*}
    B_{\tlam\t}&= B_{\tlam\tlam}D_\t
      \equiv\psi_{\tlam\t}+\sum_{\v\Gdom\t}a_{\tlam\v}\psi_{\tlam\v}\qquad\pmod\Hlam,\\
   \intertext{by \autoref{L:BstApproximation}.
   Set $b_{\tlam\v}=a_{\tlam\v}$ with $b_{\tlam\t}=1$.
   Similarly,}
    D_\s^\diamond \psi_{\tlam\tlam}&\equiv D_\s^\diamond B_{\tlam\tlam}=B_{\s\tlam}
      \equiv\sum_{\u\Gedom\s}b_{\u\tlam}\psi_{\u\tlam}\qquad\pmod\Hlam,\\
 \intertext{where $b_{\u\tlam}=a_{\tlam\u}$ with $b_{\s\tlam}=1$.
  By \autoref{C:PsiBasis}, $\{\psi_{\u\v}\}$ is a graded cellular basis of
  $\H(K)$ so, working modulo $\Hlam$,}
    B_{\s\t} &= D_\s^\diamond B_{\tlam\tlam} D_\t
    \equiv  \sum_{\v\Gedom\t} b_{\tlam\v}D_\s^\diamond\psi_{\tlam\v}
    \equiv \sum_{\v\Gedom\t} \sum_{\u\Gedom\s}
                   b_{\tlam\v} b_{\u\tlam}\psi_{\u\v}\\
    &=\psi_{\s\t} + \sum_{(\u,\v)\Gdom(\s,\t)}
    b_{\u\tlam}b_{\tlam\v}\psi_{\u\v}\qquad\pmod\Hlam.
  \end{align*}
  Setting $b_{\u\v}=\b_{\s\tlam}b_{\tlam\v}$ completes the proof.
\end{proof}

Combining these results gives us a new graded cellular basis of~$\H$.

\begin{Theorem} \label{T:CellularBBasis}
  Suppose that $K$ is a field. Then
  $\set{B_{\s\t}|(\s,\t)\in\SStd(\Parts)}$ is a graded cellular basis
  of~$\H(K)$ with cellular algebra automorphism~$\diamond$.
\end{Theorem}

\begin{proof} By \autoref{P:BBasisNonZero} and \autoref{C:PsiBasis},
  $\set{B_{\s\t}|(\s,\t\in\SStd(\Parts)}$ is a basis of $\H(K)$. By definition,
  if $(\s,\t)\in\SStd(\Parts)$ then $B_{\s\t}$ is homogeneous of degree
  $\deg\s+\deg\t$ and $B_{\s\t}^\diamond=B_{\t\s}$. Therefore, the basis
  $\{B_{\s\t}\}$  satisfies (GC$_1$), (GC$_3$) and (GC$_d$) from
  \autoref{D:cellular}. Finally, since
  $B_{\s\t}\equiv D_\s^\diamond B_{\tlam\tlam} D_\t\pmod\Hlam$, (GC$_2$) follows by
  repeating the argument from Step~3 in the proof of \autoref{T:KLbasis}.
\end{proof}

The graded cellular basis $\set{B_{\s\t}|(\s,\t)\in\SStd(\Parts)}$ of $\H(K)$ is
distinguished in the sense that, unlike~$\psi_{\s\t}$, the element $B_{\s\t}$
depends only on $(\s,\t)\in\SStd(\Parts)$ and not on a choice of reduced
expressions for the permutations $d(\s)$ and $d(\t)$.

\begin{Example}
  We give an example to show what $B$-basis elements look like.  Suppose that
  $K$ is a field of characteristic zero, that $e>2$, and let
  $\Lambda=2\Lambda_0+\Lambda_1$. Fix a multicharge
  $\charge=(\kappa_1,\kappa_2,\kappa_3)$ such that
  $\charge\equiv(0,1,0) \pmod e$ and $\charge$ satisfies
  \autoref{E:ChargeSeparation}. We use the notation of the last two
  sections, so we work over the rings $(\K,\OO,K)$ and $t=x+\xi\in\O$.

  Let $\blam=(1|1|1)$ and set $\t=\tritab({3}|{2}|{1})$.
  The permutation $d(\t)$ has two reduced expressions: $s_1s_2s_1$ and
  $s_2s_1s_2$. Let $\psi_{\tlam\t}$ and $\hat\psi_{\tlam\t}$,
  respectively, be the $\psi$-basis elements corresponding to these two
  reduced expressions. By \autoref{D:QuiverRelations},
  $\psi_1\psi_2\psi_2e(0,1,0)=\psi_2\psi_1\psi_2e(0,1,0)-e(0,1,0),$
  so $\psi_{\tlam\t}=\hat\psi_{\tlam\t}-\psi_{\tlam\tlam}\neq \hat\psi_{\tlam\t}$,
  where $0\neq\psi_{\tlam\tlam}=y_1e(0,1,0)$. The set of standard tableau with residue sequence $(0,1,0)$ and which are dominant than or equal to $\t$ is $\{\t^\bnu, \t^\bnu,\t^\blam,\t\}$, where $\bmu=(2|-|1)$ and $\bnu=(1|1^2|-)$. Of these tableaux, only~$\t$ and~$\tlam$
  have degree~$1$, so it follows that $B_{\tlam\tlam}=\psi_{\tlam\tlam}$ and
  $B_{\tlam\t}=\psi_{\tlam\t}+c\psi_{\tlam\tlam}$, for some $c\in K$. To
  compute~$c$ it is enough to work with the seminormal basis
  $\set{f_\s|\s\in\Std(\blam)}$ of the Specht module $S^\blam$ over~$\O$. Using
  \autoref{L:BetaSNCS} and \autoref{E:NewBetas},
  $$\psio_\t=f_{\tlam}\psio_1\psio_2\psio_1
  =f_{\t} -\frac{t^{\kappa_2-1-\kappa_3}[1+\kappa_1-\kappa_2]_t}
                {[\kappa_1-\kappa_3]_t}f_{\tlam}.$$
  Since $\charge\equiv(0,1,0)$ we can write
  $1+\kappa_1-\kappa_2=ae$ and $\kappa_1-\kappa_3=be$, for some~$a,b\in\Z$.
  Moreover, $a,b\ne0$ by \autoref{E:ChargeSeparation}. It is
  straightforward to check that when $x=0$ the coefficient of $f_{\tlam}$
  above is equal to~$-\frac ab$, so this coefficient is invertible
  in~$\OO$. Hence,
  $\Boo_{\tlam\t}=f_{\tlam\t}+c_1f_{\t^{\bmu}\t^{\bmu}}
                             +c_2f_{\t^{\bnu}\t^{\bnu}}\in\HO$,
  for some $c_1,c_2\in x^{-1}K[x^{-1}]$.
  Since $\deg\t^{\bmu}=\deg\t^{\bnu}=2>1$, we conclude that
  $B_{\tlam\t}=\psi_{\tlam\t}+\frac{a}{b}\psi_{\tlam\tlam}$.
\end{Example}

We have not yet proved \autoref{Thm:BBasis} from the introduction
because it is not clear that $B_{\s\t}$ is the homogeneous
component of $\Boo_{\s\t}\otimes 1_K$ of degree~$\deg\s+\deg\t$. In fact,
there is no reason why this should be true.

As in \autoref{Thm:BBasis}, suppose that $K$ is a field of
characteristic zero. Using \autoref{P:KLbasis}, it follows by
induction on the dominance ordering that the homogeneous components of
$\Boo_{\s\t}\otimes1_K$ have degree greater than or equal to
$\deg\s+\deg\t$. Moreover, if $B_{\s\t}'$ is the homogeneous component
of $\Boo_{\s\t}\otimes1_K$ of degree $\deg\s+\deg\t$ then
$$B_{\s\t}'\equiv B_{\s\t}\pmod\Hlam$$
by the proof of \autoref{T:KLbasis} (specifically the definition of
$D^\O_\s$ and $D^\O_\t$). In particular, $B_{\s\t}'\ne0$. As
$B_{\s\t}'$ is the minimial homogeneous component of $\Boo_{\s\t}\otimes1_K$,
\autoref{T:KLbasis} readily implies that
$\set{B_{\s\t}'|(\s,\t)\in\SStd(\Parts)}$ is a graded
cellular basis of~$\H$. Hence, all of the claims in
\autoref{Thm:BBasis} now follow.

If $K$ is a field of positive characteristic it is not clear if
$\Boo_{\s\t}\otimes1_K$ has homogeneous components of degree less
than~$\deg\s+\deg\t$. It is precisely for this reason that we need the
elements $D_\s$ and $D_\t$ in \autoref{D:BstBasis}.


\appendix
\def\theequation{\thesection\arabic{equation}}

\section{Seminormal forms for the linear quiver}\label{Chap:e=0}
  In this appendix we show how the results in this paper work when $e=0$ so that
  $\xi\in K$ is either not a root of unity or $\xi=1$ and $K$ is a field of
  characteristic zero.
  In order to define a modular system we have to leave the case where the cyclotomic
  parameters $Q_1,\dots,Q_\ell$ are \textit{integral}, that is, when
  $Q_l=[\kappa_l]$ for $1\le l\le\ell$. This causes quite a few notational
  inconveniences, but otherwise the story is much the same as for the case
  when~$e>0$.  We do not develop the full theory of ``$0$-idempotent subrings''
  here. Rather, we show just one way of proving the results in this paper when
  $e=0$.

  Fix a field $K$ and $0\ne\xi\in K$ of quantum characteristic~$e$. That is,
  either $\xi=1$ and $K$ is a field of characteristic zero or $\xi^d\ne1$ for
  $d\in\Z$.  The multicharge $\charge\in\Z^\ell$ is arbitrary.

  Let $\O=\Z[x,\xi]_{(x)}$ be the localisation of $\Z[x,\xi]$ at the principal
  ideal generated by~$x$. Let $\K=\Q(x,\xi)$ be the field of fractions of~$\O$.
  Define $\HO$ to be the cyclotomic Hecke algebra of type~$A$
  with Hecke parameter $t=\xi$, a unit in~$\O$, and cyclotomic parameters
  $$Q_l=x^l+[\kappa_l],\qquad\text{for }1\le l\le\ell,$$
  where, as before, $[k]=[k]_t$ for $k\in\Z$.  Then $\H(\K)=\HO\otimes_\O\K$ is
  split semisimple in view of Ariki's semisimplicity condition~\cite{Ariki:ss}.
  Moreover, by definition, $\H(K)\cong\HO\otimes_\O K$, where we consider $K$ as
  an $\O$-module by setting $x$ act on~$K$ as multiplication by zero.

  Define a new \textbf{content} function for $\HO$ by setting
  $$C_\gamma =t^{c-r}x^{l}+[\kappa_l+c-r],$$
  for a node $\gamma=(l,r,c)$. We will also need the previous definition of
  contents below. If $\t\in\Std(\Parts)$ is a tableau and $1\le k\le n$ then set
  $C_k(\t)=C_\gamma$, where $\gamma$ is the unique node such that
  $\t(\gamma)=k$.

  As in \autoref{S:Murphy}, let $\set{m_{\s\t}|(\s,\t)\in\SStd(\Parts)}$ be the
  Murphy basis of $\HO$.  Then the analogue of \autoref{L:JucysMurphyAction} is
  that if $1\le r\le n$ then
  $$m_{\s\t}L_r = C_r(\t)m_{\s\t}+\sum_{(\u,\v)\gdom(\s,\t)}r_{\u\v}m_{\u\v},$$
  for some $r_{\u\v}\in\O$. As in \autoref{S:Contents} define a $\bigast$-seminormal basis
  of $\H(\K)$ to be a basis $\{f_{\s\t}\}$ of simultaneous two-sided
  eigenvectors for $L_1,\dots,L_n$ such that $f_{\s\t}^{\bigast}=f_{\t\s}$.

  Define a \textbf{seminormal coefficient system} for $\HO$ to be a set
  of scalars $\balpha=\{\alpha_r(\s)\}$ satisfying \autoref{D:SNCS}(a),
  \autoref{D:SNCS}(b) and such that if $\s\in\Std(\Parts)$ and
  $\u=\s(r,r+1)\in\Std(\Parts)$ then
  \begin{equation}\label{EA:SNCS}
    \alpha_r(\s)\alpha_r(\u)
        = \frac{(1-C_r(\s)+tC_r(\u))(1+tC_r(\s)-C_r(\u))}{P_r(\s)P_r(\u)},
  \end{equation}
  where $P_r(\s)=C_r(\u)-C_r(\s)$, and where $\alpha_r(\s)=0$ if
  $\u\notin\Std(\Parts)$.

  As in \autoref{T:Seminormal}, each seminormal basis of $\HK$ is determined by
  a seminormal coefficient system $\balpha=\{\alpha_r(\s)\}$, such that
  $$T_rf_{\s\t} =
  \alpha_r(\s)f_{\u\t}+\frac{1+(t-1)C_{r+1}(\s)}{P_r(\s)}f_{\s\t},
  \qquad\text{where }\u=\s(r,r+1),
  $$
  together with a set of scalars $\set{\gamma_{\tlam}|\blam\in\Parts}$.  Notice
  that $I=\Z$, since $e=0$, so if $\bi\in I^n$ then $\t\in\Std(\bi)$ if and only
  if $c_r(\t)=i_r$ and, in turn, this is equivalent to the constant term of
  $C_r(\t)$ being equal to $[i_r]$, for $1\le r\le n$. Arguing as in
  \autoref{L:MurphyIdempotent},
  $$\fo = \sum_{\t\in\Std(\bi)}\frac1{\gamma_\t}f_{\t\t}\in\HO.$$

  With these definitions in place all of the arguments in
  \autoref{Chap:QuiverHeckeAlgebras} go through with only minor changes. In
  particular, if $1\le r\le n$ and $\bi\in I^n$ then \autoref{D:KLRLift} should be
  replaced by
  $$\psio_r\fo=\begin{cases}
    (T_r+1)\frac1{M_r}\fo,&\text{if }i_r=i_{r+1}\\
    (T_rL_r-L_rT_r)\fo,&\text{if }i_r=i_{r+1}+1,\\
    (T_rL_r-L_rT_r)\frac1{M_r}\fo,&\text{otherwise,}\\
  \end{cases}$$
  and $\yo_r\fo=\(L_r-C_r(\t)\)\fo$ where, as before, $M_r=1-L_r+tL_{r+1}$.
  With these new definitions, if $\s\in\Std(\bi)$, for $\bi\in I^m$, and $1\le
  r\le n$ then \autoref{L:PsiExpansion} becomes
  $$\psio_r f_{\s\t}=B_r(\s)f_{\s\t}+\frac{\delta_{i_ri_{r+1}}}{P_r(\s)}f_{\u\t},$$
  where $\u=\s(r,r+1)$ and
  $$B_r(\s)=\begin{cases}
    \frac{\alpha_r(\s)}{1-C_r(\s)+tC_{r+1}(\s)},&\text{if }i_r=i_{r+1},\\
    \alpha_r(\s)P_r(\s),&\text{if }i_r=i_{r+1}+1,\\
    \frac{\alpha_r(s)P_r(\s)}{1-C_r(\s)+tC_{r+1}(\s)},&\text{otherwise.}
  \end{cases}$$
  Observe that if $\u=\s(r,r+1)$ is a standard tableau then, using
  \autoref{EA:SNCS}, the definitions imply that
  $$B_r(\s)B_r(\u) = \begin{cases}
    \frac1{P_r(\s)P_r(\u)},&\text{if }i_r=i_{r+1},\\
    (1-C_r(\s)+tC_r(\u))(1+tC_r(\s)-C_r(\u)),&\text{if }i_r\leftrightarrows i_{r+1},\\
    (1+tC_r(\s)-C_r(\u)),&\text{if }i_r\rightarrow i_{r+1},\\
    (1-C_r(\s)+tC_r(\u)),&\text{if }i_r\leftarrow i_{r+1},\\
    1,&\text{otherwise.}
  \end{cases}$$
  Comparing this with \autoref{L:BetaSquared}, it is now easy to see that
  analogues of \autoref{P:PsiSquare} and \autoref{P:PsiBraid} both hold in this
  situation. Hence, repeating the arguments of \autoref{S:KLRDeformation},
  a suitable modification of \autoref{Thm:KLRDeformation} also holds. Similarly,
  the construction of the bases in \autoref{Chap:PsiBases} and
  \autoref{Chap:BBasis} now goes though largely without change.

\section*{Acknowledgments}

Both authors were supported by the Australian Research Council. The first author
was also supported by the National Natural Science Foundation of China.

\bibliographystyle{andrew}

\end{document}